\documentclass[11pt,oneside,english]{article}
\usepackage{amssymb}

\usepackage{graphicx}
\usepackage[all]{xy}
\usepackage[T1]{fontenc}
\usepackage{amsmath}

\usepackage[french]{babel}

\def \S{{\mathcal S}}

\newtheorem{theo}{Theorem}[section]

\newtheorem{remark}{Remark}

\newtheorem{cor}[theo]{Corollary}
\newtheorem{defi}[theo]{Definition}

\newtheorem{prop}[theo]{Proposition}
\newtheorem{lem}[theo]{Lemma}

\newtheorem{exo}{Exercise}

\usepackage{authblk}

 \usepackage{epigraph}

\begin{document}
	\title{ Hyperelliptic tangential covers and even elliptic finite-gap potentials, back and forth. }

	\date{\today}

	\author{{\Large Armando Treibich}\footnote{E-mail: treibicharmando@gmail.com. 
	ORCID: 0000-0002-0105-7717.}}
	\affil{DMA CURE - Udelar, Uruguay\vspace{1cm}
	
	\vspace{1cm}

\rightline{\normalsize\textit{Dedicated to the memory of Jean Louis Verdier and Igor M. Krichever.}}
	}
	\date{}
	
	\maketitle

	\renewcommand{\abstractname}{Abstract}
	\begin{abstract} 
	Let $(X,\omega_0):=(\mathbb{C}/\Lambda,0)$ denote the elliptic curve associated to the lattice $\Lambda$, $X_2:=\{\omega_0,\cdots, \omega_3\}$ its set of half-periods and $\wp:X \to \mathbb{P}^1$ the usual Weierstrass $\wp$ function, with a double pole at the origin $\omega_0$. Let $\pi:(\Gamma,p) \to (X,\omega_0)$ be a degree-$n$ ramified cover, marked at a smooth point $p$ over $\omega_0$, and consider the (rational) Abel map $Ab_p:\Gamma \to Jac\Gamma$ and the dual map $\pi^\vee:= Ab_p\circ\pi^* : X \to Jac\Gamma$ into the Jacobian variety of $\Gamma$. We then call $\pi$ a \textit{hyperelliptic tangential} (ht) cover if, and only if; $\Gamma$ is a hyperelliptic curve, $p\in \Gamma$ a Weierstrass point and the images of $\Gamma$ and $X$ in $Jac\Gamma$ are tangent at its origin. To any such ht cover $\pi$ we attach an integer vector $\mu \in \mathbb{N}^4$, so-called the \textit{type}, satisfying $\mu_0+1\equiv \mu_1 \equiv \mu_2 \equiv \mu_3 \equiv n \,\textrm{mod}.2\;\textrm{and} \,2n+1\,\;\textrm{-}\,\sum_i \mu_i^2=4d$, for some $d\in \mathbb{N}$. Whenever $\Gamma$ is smooth, the type $\mu$ gives the number of Weierstrass points of $\Gamma$ (different from $p$) over each half-period $\omega_i$ ($i=0,\dots,3$). We let $\mathbb{T}_0:= \{\mu\in \mathbb{N}^4, \mu_0+1\equiv \mu_1 \equiv \mu_2 \equiv \mu_3 \equiv n \,\textrm{mod}.2\}$, $\mathcal{S}\mathcal{C}_X(\mu,d)$ and  $\mathcal{S}\mathcal{D}_X(\mu,d)$ denote, respectively, the set of \textit{types}, the set of degree-$n$ $\it{ht}$ covers of type $\mu$, and the set of spectral data $(\pi,\xi)$, with $\pi \in \mathcal{S}\mathcal{C}_X(\mu,d)$ and $\xi$ a theta-characteristic. For any such $(\pi,\xi)$ the corresponding even, $\Lambda$-periodic, finite-gap potential, decomposes as 
	
		$$u_\xi(x)=\sum_0^3\alpha_i(\alpha_i+1)\wp(x\,\textrm{-}\,\omega_i) +2\sum_{j=1}^m \left(\wp(x\,\textrm{-}\,\rho_j)+\wp(x+\rho_j)\right),$$
		
for some $(\alpha,m)\in \mathbb{N}^4\times \mathbb{N}$ satisfying $2n=\sum_i\alpha_i(\alpha_i+1)+4m$, and any solution $\{\rho_j\} \in (X \setminus X_2)^{(m)}$ of the so-called (D-G) square system of equations.
  
	We let $\mathcal{P}ot_X(\alpha,m)$ denote the set of such potentials.\\
\indent Conversely, given any potential $u(x)\in \mathcal{P}ot_X(\alpha,m)$, its spectral data belongs to some $\mathcal{S}\mathcal{D}_X(\mu,d)$ as above. We thus have a bijection 
	$$(\pi,\xi) \in \bigcup_{(\mu,d)}\mathcal{S}\mathcal{D}_X(\mu,d) \mapsto u_\xi \in \bigcup _{(\alpha,m)}\mathcal{P}ot_X(\alpha,m)\,.$$
	
The problem at stake is to find out its inverse map and in particular the preimage of each $m$-th family $\mathcal{P}ot_X(m):=\bigcup _{\alpha\in \mathbb{N}^4}\mathcal{P}ot_X(\alpha,m).$
	The latter problem has been thoroughly studied for $\mathcal{P}ot_X(0)$ and $\mathcal{P}ot_X(1)$. In this article we go one step further, by constructing all spectral data in the pre-image of each $\mathcal{P}ot_X(\alpha,2)$.\\
	\indent It turns out that any $\it{ht}$ cover can be naturally equipped with eight fundamental theta-characteristics $\{\xi_{j,k}, (j,k) \in  \mathbb{Z}_2 \times \mathbb{Z}_4\}$. Better still, there are eight maps $\{\mathcal{C}^{j,k}:\mathbb{T}_0 \to \mathbb{N}^4\}$ such that for $\alpha=\mathcal{C}^{j,k}(\mu)$, any spectral data $(\pi,\xi) \in \mathcal{S}\mathcal{D}_X(\mu,d)$ with $\xi =\xi_{j,k}$, gives rise to a potential in $\mathcal{P}ot(\alpha,d)$. Hence, $\#\mathcal{S}\mathcal{C}_X(\mu,d) \leq \#\mathcal{P}ot_X(\alpha,d)$. We then focus on the case $d=2$ and obtain the bound $\#\mathcal{P}ot_X(\alpha,2)\leq 27$ for any $\alpha \in \mathbb{N}^4$. Prior to that, we construct for any given $\mu \in \mathbb{T}_0$ a weak del Pezzo surface, with a suitable family of rational curves from which we can reconstruct any $\it{ht}$ cover in $\mathcal{S}\mathcal{C}_X(\mu,2)$, for any $X$. We then deduce $\#\mathcal{S}\mathcal{C}_X(\mu,2)$ for a generic $X$. For example, let $M:=Max\{\alpha_i\}\,,\;S=\sum_i \alpha_i$, $m:=min\{\alpha_i\}$ and denote $g_\alpha:=\frac{1}{2}Max\left(2M,S+1\,\textrm{-}\,(1+(\,\textrm{-}\,1)^S)(m +\frac{1}{2})\right)$. Then $\#\mathcal{S}\mathcal{C}_X(\mu,2)\leq 27$, with equality if and only if $|2M\,\textrm{-}\,S+(1+(\,\textrm{-}\,1)^S)m| \geq 4$, in which case $\#\mathcal{P}ot_X(\alpha,2)=27$ as well. In particular $\{(\pi,\xi_{j,k}), \pi \in \mathcal{S}\mathcal{C}_X(\mu,2)\}$ is the set of spectral data giving rise to $\mathcal{P}ot_X(\alpha,2)$. They all have same arithmetic genus, namely: $g_\alpha$.\\
	 If $|2M\,\textrm{-}\,S+(1+(\,\textrm{-}\,1)^S)m| \leq 3$ instead, $\#\mathcal{S}\mathcal{C}_X(\mu,2)<27$ but we find $27\,\textrm{-}\,\#\mathcal{S}\mathcal{C}_X(\mu,2)$ other spectral data giving rise to elements in $\mathcal{P}ot_X(\alpha,2)$. Therefore, we still have $\#\mathcal{P}ot_X(\alpha,2)=27$ but the arithmetic genus ranges between $g_\alpha$ and $g_\alpha+2$. \\
	\indent Summing up the results obtained for any $d\leq 2$, we conclude with a natural conjecture, leading to a recursive formula in $d\in \mathbb{N}$, for the cardinals of both, $\mathcal{P}ot_X(\alpha,d)$ and $\mathcal{S}\mathcal{C}_X(\mu,d)$.

	\end{abstract}
	
	\section{Introduction}
	In this survey we consider a (third) family of (even) finite-gap elliptic potentials, and find their spectral data as well. We continue earlier work done with J.-L.Verdier (cf. \cite{TV1}, \cite{TV2}, \cite{TV3} \& \cite{T2}) and prove a refined version of the last conjecture stated in \cite{TV1}. Recall that one of many equivalent definitions states that an elliptic finite-gap potential (also called hyper-elliptic potential) is the initial condition function, $v(x)=u(x,0)$, of an elliptic KdV soliton $u(x,t)$. The latter, either $u(x,0)$ or $u(x,t)$, can be recovered from its spectral data, namely, a $\it{ht}$ cover over $(X,\omega_0):=(\mathbb{C}/\Lambda, 0)$, equipped with a point of its compactified Jacobian. In particular, the even hyper-elliptic potentials correspond to theta-characteristics on the corresponding spectral curves. Their study goes back, in modern terminology, to the seminal work of H.Airault, H.McKean and J.Moser \cite{AKM}, and S.Novikov and co-workers (cf. \cite{N},\cite{D1},\cite{DN},\cite{DMN},\cite{S1}, \cite{S2}, \cite{K1},\cite{K2},\cite{K3}, \cite{D2}, or the encyclopedia \cite{BeBoEM}, \cite{SW} and references therein). However, the first examples were found by Ince, decades before, namely: the family $\{g(g+1)\wp(x\,\textrm{-}\,ic), g > 0\}$, for a rectangular lattice $\Lambda=\mathbb{Z}\oplus i2c\mathbb{Z}$ (ref. \cite{I}). Later on, this family was completed to $\left\{\sum_{i=0}^3 \alpha_i(\alpha_i+1)\wp(x\,\textrm{-}\,\omega_i), \alpha\in \mathbb{N}^4\setminus\{\vec{0}\}\right\}$, with respect to any lattice $\Lambda \subset \mathbb{C}$, and their spectral data were obtained as well (cf. \cite{TV2}, see also \cite{GW}). A second family followed, namely, the set of functions
	 $$\sum_{i=0}^3 \alpha_i(\alpha_i+1)\wp(x\,\textrm{-}\,\omega_i)+ 2\wp(x\,\textrm{-}\,\rho)+2\wp(x+\rho),$$
	 
	 for any $ \alpha\in \mathbb{N}^4$ and $\rho \notin X_2$ such that $\sum_{i=0}^3 (2\alpha_i+1)\wp'(\rho\,\textrm{-}\,\omega_i)=0$. Their complete spectral data were also found, in terms of $\alpha$ (cf. \cite{T2}).\\
	 \indent More generally, it is a well known fact that any (even) hyper-elliptic potential is an (even) Picard potential and decomposes, up to a constant, as 
	 $$\sum_{i=0}^3 \alpha_i(\alpha_i+1)\wp(x\,\textrm{-}\,\omega_i)+ \sum_{l=1}^h \beta_l(\beta_l+1)\big(\wp(x\,\textrm{-}\,\rho_l)+\wp(x+\rho_l)\big),$$
	 
	 for some $\alpha\in \mathbb{N}^4$, $l\geq 0$, $\beta_l\geq 1$ for any $l=1,\cdots,h$ and $h$ distinct points $\{\rho_l\}$ in $X\setminus X_2$, satisfying a so-called Duistermaat-Gr$\ddot{\textrm{u}}$nbaum set of $\sum_l \beta_l$ equations (cf. \cite{GW}). For fixed $n:=\frac{1}{2}\sum_{i=0}^3 \alpha_i(\alpha_i+1)+ \sum_{l=1}^m \beta_l(\beta_l+1)$ the set of solutions, for all possible $(\alpha,h,\beta)$, is finite (although may be empty), and corresponds  to the set of theta-characteristics of all degree-$n$ $\it{ht}$ covers (cf. \cite{TV1}). If $\sum_l \beta_l > h >0$, the latter system is overdetermined and we conjecture it has no solution for a generic lattice $\Lambda \subset \mathbb{C}$. It seems therefore natural to restrict to the remaining cases, i.e.: $\beta_l=1$ for all $l$.\\
	  \indent Given $(n,\alpha,h)$ such that $2n=\sum_i\alpha_i(\alpha_i+1)+4h$, bounds for the spectral invariants, such as the type or arithmetic geni of the corresponding $\it{ht}$ covers can be obtained. However, despite its apparent simplicity, proving that for any fixed lattice $\Lambda$, $h>0$ and $\alpha\in \mathbb{N}^4$, the set of hyper-elliptic potentials is not empty, and (bounding or better still) calculating its cardinality, is yet out of reach, although we conjecture its value (see §$5$, Remark $8$). The latter problem is equivalent to finding all degree-$n$ $\it{ht}$ covers over $X$, and correspond to certain rational irreducible curves in a rational surface denoted $\widetilde{S}_X$, canonically associated to the elliptic curve $X$ . The first two families, i.e.: arbitrary $\alpha\in \mathbb{N}^4$ and $h=0$ or  $h=1=\beta_1$, were reduced to the construction of infinitely many exceptional curves and associated pencils of curves in $\widetilde{S}_X$. In the present paper we continue the study of rational curves in $\widetilde{S}_X$ and apply it to the case $h=2=\beta_1+\beta_2$, obtaining a third family of hyper-elliptic potentials, including their spectral data. To put our results in perspective and sketch the structure of this article, we start with a brief historical review of the subject and introduce the basic definitions.\\
	 \indent For almost half a century the Theory of finite-gap integration and its algebro-geometric counterpart, were developed and applied to the construction of exact solutions of many non-linear partial differential equations, as well to the study of (algebraically completely) integrable systems. It is hardly possible to present a fair list of publications, but one can trace back its origins (one of them at least), to the russian school led initially by S. Novikov and collaborators. Starting from the Cauchy problem for the Korteweg-deVries (KdV) equation and its connection with finite-gap potentials, it was shown that their spectral data have a deep algebro-geometric background, e.g.: they could be reinterpreted as (and recovered from) suitable linear orbits in Jacobians of hyperelliptic curves, marked at a Weierstrass point. This led in particular to striking exact formulas in terms of Riemann's theta functions, e.g.: the Its-Matveev formula for the KdV equation. Soon after I. Krichever developed a most profound and yet concrete tool, through the (re-)creation of Baker-Akhiezer's (BA) function. It enlarged the whole picture, on one side to the Kadomtsev-Petviashvili (KP) equation, a three-variable generalization of KdV, and to Jacobians of any projective curve on the other side. From then on, a rapidly increasing list of connections flourished in both directions, between Soliton Theory and Algebraic Geometry (e.g.: the characterization of Jacobian and Prym varieties in terms of geometric properties of their theta divisors, and effective constructions of particular solutions, like rational, periodic or doubly periodic ones with respect to the first flow, to mention a few ones).\\
	 \indent Recall that a Schr\"{o}dinger operator $\partial_x^2\,\textrm{-}\,v(x)$ is called finite-gap, if and only if there exists an odd-degree differential operator, say $P$, which commutes with $L$. In the latter case, a century old result of Burchnall and Chaundy states that they are polynomially related. More precisely, there exists a polynomial $\Pi_0^{2g} (x\,\textrm{-}\,\alpha_k) $, of degree $2g+1=degP$, such that $P^2= \Pi_0^{2g} (L\,\textrm{-}\,\alpha_k) $. Let $\Gamma$ denote the one-point compactification of the complex plane curve $\{y^2=\Pi_0^{2g} (x\,\textrm{-}\,\alpha_k)\}$, $p\in \Gamma$ the corresponding smooth Weierstrass point (at "$\infty$") and $\theta_\Gamma$ its Riemann theta function (or eventually the $\tau$-function if $\Gamma$ has singularities, i.e.: if $\Pi_0^{2g} (x\,\textrm{-}\,\alpha_k)$ has multiple roots). The finite-gap potential $v(x)$ can then be recovered from suitable data $(\Gamma,p,\lambda, \xi)$, where $\lambda$ is a local coordinate at the Weierstrass point $p\in \Gamma$, and $\xi$ a point of its Jacobian. Consider the marked Abel embedding $A_p:p'\in\Gamma \mapsto \mathcal{O}_\Gamma(p'\,\textrm{-}\,p) \in Jac\Gamma$ into its Jacobian variety, pick $Z \in \mathbb{C}^g$ projecting onto $\xi \in Jac\Gamma$ and let $U\in \mathbb{C}^g$ denote the first derivative $U:=\frac{\partial}{\partial \lambda}A_p(\lambda)_{|\lambda =0} $. The latter generates the tangent to $A_p(\Gamma)$ at the origin $\mathcal{O}_\Gamma$. The Its-Matveev formula (cf. \cite{IM}) tells us that $$v(x)=\textrm{-}\,2\partial_x^2log\theta_\Gamma(xU+Z).$$
	 In case $\Gamma$ is singular one should replace $\theta_\Gamma$ by the corresponding $\tau$-function, and $Jac\Gamma$ by the compactified Jacobian $W(\Gamma)$ (i.e.: the moduli space of rank $1$, degree $g\,\textrm{-}\,1$, torsion-free sheaves on $\Gamma$). We thus obtain a family of isospectral finite-gap potentials, parameterized by $Jac\Gamma$. We also remark that any such potential $v(x)$ is $\Lambda$- periodic, for some lattice $\Lambda \subset \mathbb{C}$, if and only if there exists a copy of the elliptic curve $(X,\omega_0):=(\mathbb{C}/\Lambda,0)$ inside $Jac\Gamma$, which is tangent to $A_p(\Gamma)$ at $A_p(p)$. Dualizing the above group homomorphism $( X ,\omega_0)\to  (Jac\Gamma, \mathcal{O}_\Gamma)$, and composing with the Abel map $A_p$, one obtains a $\it{ht}$ cover $\pi:( \Gamma,p) \to (X,\omega_0 )$. Conversely, given such a marked cover $\pi$, the image of $X$ by its dual morphism $\pi^\vee:X \to Jac\Gamma$ satisfies the above tangency condition if and only if $\pi$ is a $\it{ht}$ cover. The first basic invariant of $\pi$ is its degree, say $n:=deg\pi$. Let $X^{(n)}$ denote the $n$th symmetric power of $X$, $\Delta_n \subset X^{(n)}$ its generalized diagonal and $\mathbb{L}_n \subset X^{(n)}\setminus \Delta_n$ the common zero-locus of the $n$ equations $\sum_{j\neq i}\wp'(x_j\,\textrm{-}\,x_i)=0$, for $i=1,\cdots,n$ (cf. \cite{AKM}). For any generic $\xi \in Jac\Gamma$ the corresponding hyper-elliptic potential has $n$ double poles and decomposes as $v(x)= 2\sum_{i=1}^n\wp(x\,\textrm{-}\,x_i)$, for a unique $\{x_i\}\in \mathbb{L}_n$. Conversely, given any point in the closure of $\mathbb{L}_n$, say $\{x_i\}\in \overline{\mathbb{L}}_n\subset X^{(n)}$, the function $v(x)= 2\sum_{i=1}^n \wp(x\,\textrm{-}\,x_i)$ is a hyper-elliptic potential. We will say that $v(x)$ has degree $n$. \\
	 \indent The first concrete examples of such potentials had already been constructed by E. Ince (ref. \cite{I}), but only for a rectangular lattice $\Lambda_c = \mathbb{Z}\oplus i2c\mathbb{Z}$ ($c\in \mathbb{R}^\ast$). He proved in fact that for any $g\in \mathbb{N}^\ast$, $L:=\partial_x^2\,\textrm{-}\,g(g+1)\wp(x+ic)$ is a self-adjoint operator with finite spectrum. Its interpretation as finite-gap potential, and calculation of corresponding spectral data, were given much later (e.g.: \cite{BeBoEM}). \\
	 \indent We recall hereafter the basic invariants and properties of any $\it{ht}$ cover. \\
	 \indent Let $\pi_X:S_X \to X$ denote the unique ruled surface possessing a unique section $C_0 \subset S_X$ of zero self-intersection, $[\,\textrm{-}\,1]:(X,\omega_0) \to (X,\omega_0)$ the inverse homomorphism, $\tau:S_X \to S_X$ its lift to $S_X$ and $e:\S_X^\perp \to S_X$ the blowing-up of the eight fixed points of $\tau$. The involution $\tau$ lifts to $ \tau^\perp: S_X^\perp \to S_X^\perp$ and does not have isolated fixed points. Taking the quotient by its action we finally get a degree-$2$ projection, $\varphi: S_X^\perp \to \widetilde{S}_X$, onto a smooth rational anticanonical surface. The pair of morphisms $(e,\varphi)$ gives a common frame for the study of all $\it{ht}$ covers over $X$. Let us state their fundamental properties.
	 \begin{enumerate}
	 	\item The subset of half-periods $X_2:=\{\omega_0,\cdots,\omega_3\}\subset X$ acts canonically on $\widetilde{S}_X$ and induces a free action on its set of exceptional curves. Modulo this action, any such curve corresponds to a unique \textit{type} $\mu\in\mathbb{T}_0:= \big\{\mu \in \mathbb{N}^4,\mu_0+1 \equiv \mu_1\equiv \mu_2\equiv \mu_3\textrm{  mod}.2\big\} $ (cf.\cite{TV1}-§$5.5.1)$). We denote hereafter $\mu^{(2)}:=\sum_i \mu_i^2$ and $\mu^{(1)}:=\sum_i \mu_i$ to simplify the presentation.
	 	
	 	\item Any degree-$n$ $\it{ht}$ cover $\pi: (\Gamma,p) \to (X,\omega_0)$ factors via $\pi_X\circ e:S_X^\perp \to X$ (cf.\cite{TV1}-§$4.3$), and its projection in $\widetilde{S}_X$, say $\widetilde{\Gamma}$, is a rational irreducible curve. More precisely, there exists a unique exceptional curve $\widetilde{\Gamma}_\mu$ ($\mu \in \mathbb{T}_0$), such that $\widetilde{\Gamma}$ is numerically equivalent to $\widetilde{\Gamma}_\mu+d(\,\textrm{-}\widetilde{K}_X)$, where $\widetilde{K}_X$ denotes the canonical divisor, $\mu_0+1\equiv n \textrm{  mod.}2$ and $d:=\frac{1}{4}\big(2n+1\,\textrm{-}\,\mu^{(2)}\big)$ (cf.\cite{TV1}-§$4.5.1)$). Its arithmetic genus is bounded by $\frac{1}{2}(\mu^{(1)}\,\textrm{-}\,1)$. We call $\mu$ the $\it{type}$ of $\pi$.

	 	\item The linear system $|\widetilde{\Gamma}_\mu+d(\textrm{-}\widetilde{K}_X)|$ has dimension and arithmetic genus equal to $d$, and carries a finite number of rational irreducible curves, maybe none (\cite{TV1}-§$4.9$). It reduces to a single element when $d=0$, namely $\widetilde{\Gamma}_\mu$.  We let $\mathcal{S}\mathcal{V}_X(\mu,d)$ denote the (finite) so-called Severi set of rational irreducible divisors in $|\widetilde{\Gamma}_\mu+d(\textrm{-}\widetilde{K}_X)|$.
	 	
	 	\item Conversely, given such a rational irreducible curve $\widetilde{\Gamma}$, there exists a unique degree-$n$ $\it{ht}$ cover of type $\mu$ and arithmetic genus $g:=\frac{1}{2}(\mu^{(1)}\,\textrm{-}\,1)$ projecting onto $\widetilde{\Gamma} \subset \widetilde{S}_X$. Its compactified Jacobian is isomorphic to an irreducible component of $\overline{\mathbb{L}}_n \subset X^{(n)}$ (cf.\cite{TV4}-§$8.10$ \& $8.11$, see also \cite{T1}-§$1.12$). It follows that $\overline{\mathbb{L}}_n$ decomposes as a disjoint union of irreducible components, each one isomorphic to the Jacobian of such a spectral curve. In particular, there is a finite number of degree-$n$ even hyper-elliptic potentials.
	 	
	 	\item  Consider any triplet $(n,\mu,d)$ as above and assume $\mathcal{S}\mathcal{V}_X(\mu,d)$ is not empty. Pick any $\widetilde{\Gamma} \in \mathcal{S}\mathcal{V}_X(\mu,d)$ and let $\pi:(\Gamma,p) \to (X,\omega_0)$ denote the corresponding degree-$n$ $\it{ht}$ cover of type $\mu$ and arithmetic genus $g:=\frac{1}{2}(\mu^{(1)}\,\textrm{-}\,1)$. Given $\xi \in W(\Gamma)$, we let $Orb_\xi : X \to  W(\Gamma),\quad \omega \mapsto \xi(\pi^\ast(\omega)\,\textrm{-}\,\pi^\ast(\omega_0))$ denote the orbital morphism. Recall also the theta divisor $\Theta_{\Gamma} := \{\xi \in W(\Gamma), h^0(\Gamma, \xi) > 0\}$ and let $\sum_{i=1}^n \alpha_i:= Orb_\xi^\ast(\Theta_{\Gamma} )$. Then, the hyper-elliptic potential associated to $\xi$ is equal to $v_\xi (x):=2\sum_{i=1}^n \wp(x\,\textrm{-}\,\alpha_i)$ (\cite{TV4}-App.A, see also \cite{T2}). In order to calculate $\sum_{i=1}^n \alpha_i$ we will need the following formula (cf. \cite{TV2}-§$3.2$):\\
	 	let $m:=Max\{k \in \mathbb{N}, h^0(\Gamma,\xi(\textrm{-}kp) > 0\}$, then $\Theta_\Gamma$ and $Orb_\xi(X)$ intersect at $\xi$ with multiplicity $\frac{1}{2}(m+1)(m+2)$. 
	 	
	 	\item We apply the above formula to the following theta-characteristics ($k\in \mathbb{Z}_4)$:
	 	$$\xi_{0,k}:=\mathcal{O}_\Gamma\big((g\,\textrm{-}1\,\textrm{-}\,2n)p+\pi^\ast(\omega_0+\omega_k)\big) \quad\textrm{and}\quad \xi_{1,k}:=\mathcal{O}_\Gamma\big((g\,\textrm{-}1\,\textrm{-}\,n)p+\pi^\ast(\omega_k)\big)\,.$$
	 	
	 	The even hyper-elliptic potential associated to each $\xi_{j,k}$, say $v_{\xi_{j,k}}(x)$, can be (partially) calculated. It decomposes as:
	 		$$ v_{\xi_{j,k}}(x)=\sum_{i=0}^3 \alpha_i(\alpha_i+1)\wp(x\,\textrm{-}\,\omega_i)+\sum_{l=1}^d 2\big(\wp(x\,\textrm{-}\,\rho_l)+\wp(x+\rho_l)\big)\;,$$
	 		
	 	for some $\{\rho_l\}\in X^{(d)}$ disjoint from $X_2$. The vector $\alpha\in \mathbb{N}^4$ is explicitely calculated in terms of $\big(\mu,(j,k)\big)$ in cf. \cite{TV2}-§3.4, and satisfies $\sum_i \alpha_i(\alpha_i+1)=\mu^{(2)} \,\textrm{-}\,1$. For any other theta-characteristic $\xi$ we know at least that,
	 	$$ v_\xi(x)=\sum_{i=0}^3 a_i(a_i+1)\wp(x\,\textrm{-}\,\omega_i)+\sum_{l=1}^h \big(2\wp(x\,\textrm{-}\,\rho_l)+2\wp(x+\rho_l)\big)$$
	 	
	 	with $\sum_i a_i(a_i+1) < \mu^{(2)} \,\textrm{-}\,1$ and $\sum_i a_i(a_i+1) + 4h =2n=\sum_i \alpha_i(\alpha_i+1)+4d $ (hence $h>d$).
	 	Conversely, given any $\alpha\in \mathbb{N}^4$, there exists  a unique $\mu\in  \mathbb{T}_0$ such that $\sum_i \alpha_i(\alpha_i+1)=\mu^{(2)} \,\textrm{-}\,1$, explicitely calculated in terms of $\alpha$, and $(j,k) \in \mathbb{Z}_2\times \mathbb{Z}_4$ (cf.\cite{TV2}-§$4.4$) with the following property. Let $d\in \mathbb{N}$ and pick any $\{\rho_l\} \in X^{(d)} \setminus \Delta_d$, disjoint from $X_2$, such that $\sum_i\alpha_i(\alpha_i+1)\wp(x\,\textrm{-}\,\omega_i)+\sum_{l=1}^d \big(2\wp(x\,\textrm{-}\,\rho_l)+2\wp(x+\rho_l)\big)$ is a hyper-elliptic potential. Then, its spectral curve is a degree $n:=\frac{1}{2}\sum_i \alpha_i(\alpha_i+1) +2d$ $\it{ht}$ cover of type $\nu$ such that $\nu\,\textrm{-}\,\mu \in 2\mathbb{N}^4$. It can also be proved that $\nu=\mu$ if and only if the corresponding theta-characteristic is equal to $\xi_{j,k}$ (see for example the proof of §$5.12$ or the preceding lemmas). We deduce the following bounds for its arithmetic genus $g$. \\
	 	
	 	\indent Let us denote $M:=Max\{\alpha_i\}$, $m:=min\{\alpha_i\}$, $S:=\sum_i \alpha_i$ and  
	 	$$g_{\alpha}:=\frac{1}{2}(\mu^{(1)} \,\textrm{-}\,1)=\frac{1}{2}Max\{2M, S + 1\,\textrm{-}\,(1+(\textrm{-}1)^S )(m+\frac{1}{2})\}\;.$$
	 		 
	 	Then $g$ ranges from $ g_{\alpha}$ to at most $g_{\alpha} +2\sqrt{d}$. \\
	 	\indent Last but not least, they all have the same arithmetic genus, i.e.: $g=g_{\alpha} $, in case $min\{\mu_i\}\geq d$, which amounts to impose $|2M\,\textrm{-}\,S+(1+(\textrm{-}1)^S )m|\geq 2d$.
	 		 
	 	\end{enumerate}
 	
 	Let me also make a disgression, following the Referee's request, to mention previous results which generalize some of the latter results to the discrete and matrix cases:
 	\begin{enumerate}
 		\item in IMRN, Vol. 2003, Issue 6, 2003, pp. 313–360, doi.org/10.1155/S1073792803204104, we find new classes of difference elliptic finite-gap operators which are difference analogs of the classically known elliptic finite-gap Schr\"odinger operators;
 		
 		\item in Funct. Anal. \& Its Appl 50, 308–318 (2016), https://doi.org/10.1007/s10688-016-0161-0 we study doubly periodic (in $x$) solutions of the matrix KdV equation $$U_t=\frac{1}{4}(3UU_x +3U_xU + U_{xxx})$$ and its KP generalization. We find in particular, arbitrarily high genus curves giving rise to such solutions and present simple polynomial equations defining them.
 	\end{enumerate}
 	The above results go back and forth from theta-characteristics of $\it{ht}$ covers to even hyper-elliptic potentials, and sketch the global picture. However, the following issues remain unanswered:
 	\begin{enumerate}
 		\item given an elliptic curve $X$ and $(\mu,d)\in \mathbb{T}_0\times\mathbb{N}$, is $\mathcal{S}\mathcal{V}_X(\mu,d)$ empty or not?;
 		\item for a generic elliptic curve $X$, calculate $\#\mathcal{S}\mathcal{V}_X(\mu,d)$.
 	\end{enumerate} 
	 	\indent Roughly speaking, the latter issues are equivalent to the following unanswered questions. \\
	 	\indent Given an elliptic curve $X$ and $(\alpha,d)\in \mathbb{N}^4\times\mathbb{N}$, are there hyper-elliptic potentials over $X$ decomposing as $$ \sum_{i=0}^3 \alpha_i(\alpha_i+1)\wp(x\,\textrm{-}\,\omega_i)+\sum_{l=1}^d \big(2\wp(x\,\textrm{-}\,\rho_l)+2\wp(x+\rho_l)\big)\,\, ,$$ for some $\{\rho_l\}\in X^{(d)}\setminus \Delta_d$ disjoint from $X_2$, and how many? \\
	 	 \\
	 	\indent In previous work we calculated their number when $d \leq 1$, and found the spectral data for each one of them.
	 	
	 	\indent For any $d> 0$ we do know that any such potential is hyper-elliptic if and only if $\{\rho_l\}$ satisfies the Duistermaat-Gr\"{u}nbaum (D-G) set of equations (cf. \cite{GW}), i.e., $ \forall l=1,\cdots,d \;$:
	 	
	 	$$ 8\sum_{k\neq l}\left(\wp'(\rho_l\,\textrm{-}\,\rho_k)+\wp'(\rho_l+\rho_k)\right)+\sum_{i=0}^3(2\alpha_i+1)^2\wp'(\rho_l\,\textrm{-}\,\omega_i)=0\, .$$ We also have partial results about their spectral data, assuming $\mathcal{P}ot_X(\alpha,d)\neq \emptyset$. \\
	 	\indent  In this article we work out the next cases, i.e.: $\mathcal{S}\mathcal{V}_X(\mu,2)$ and $\mathcal{P}ot_X(\alpha,2)$.\\
	 	 \indent When $d=2$ the (D-G) system reduces, via the classical addition formulas for $\wp$ and $\wp'$, to the common zeroes of two degree-$9$ homogeneous polynomials in the variables $[x_1:x_2:e_1:e_2]:=[\wp(\rho_1):\wp(\rho_2):\wp(\omega_1):\wp(\omega_2)] \in \mathbb{P}^3(\mathbb{C})$. It easily follows that for fixed $[e_1:e_2]$, there are at most $27$ even hyper-elliptic potentials as above. According to the preceding discussion, their spectral curves are degree-$n$ $\it{ht}$ covers, corresponding to an element of $\mathcal{S}\mathcal{V}_X(\nu,h)$ such that $h=\frac{1}{4}(2n+1\,\textrm{-}\,\nu^{(2)})\leq 2$.
	 	
	 	\indent  At last we calculate $\#\mathcal{S}\mathcal{V}_X(\mu,2)$ in terms of $\mu$ for a generic elliptic curve $(X,\omega_0)$, from which we deduce that $\#\mathcal{P}ot_X(\alpha,2)=27$ for any $\alpha \in \mathbb{N}^4$, and construct their spectral data.\\
	 	 
	 	 \indent The structure of the paper is as follows.\\
	 	 
	 	 \indent Section §$2$ recalls the construction of the family $\big \{\widetilde{S}_X \big \}$, of rational anticanonical surfaces we will deal with. Each one of them is canonically attached to a given elliptic curve $(X,\omega_0)=(\mathbb{C}/\Lambda,0)$ and, up to automorphism, the exceptional curves are parameterized by $\mathbb{T}_0$. We also present, for any $(\mu,d) \in \mathbb{T}_0 \times \mathbb{N}$, the  linear system $|\widetilde{\gamma}_X(\mu,d)|$ of $\widetilde{S}_X $, whose Severi set of rational irreducible curves, denoted $\mathcal{S}\mathcal{V}_X(\mu,d)$, gives rise to all $\it{ht}$ covers, of degree $n:=\frac{1}{2}(\mu^2\,\textrm{-}\,1)+2d$ over $X$ and type $\mu$.\\
	 	 
	 	  \indent From then on we focus on the case $d=2$. We construct a pair of degree-$2$ projections of $\widetilde{S}_X$ onto $\mathbb{P}^2(\mathbb{C})$, and determine their main properties. They give two complementary approaches to the study of $\mathcal{S}\mathcal{V}_X(\mu,2)$, reducing it to tangency properties between lines and conics with cubics.\\
	 	   
	 	 \indent We start Section §$3$ fixing a type $\mu \in \mathbb{T}_0$ and defining a birational morphism onto a degree-$2$ weak del Pezzo surface, $\widetilde{S}_X \to \underline{S}_\mu$, which according to Appendix \textbf{A} is independent of $X$. As such it carries a canonical so-called Geiser degree-$2$ projection $g_\mu :\underline{S}_\mu \to \mathbb{P}^2(\mathbb{C})$, with discriminant a reduced quartic $\mathcal{Q}$. Any $\widetilde{\Gamma} \in |\widetilde{\gamma}_X(\mu,2)|$ projects with degree $2$ onto a projective conic, and belongs to $\mathcal{S}\mathcal{V}_X(\mu,2)$ if and only if the conic has specific tangency properties with $\mathcal{Q}$. We thus obtain a common frame for studying $\mathcal{S}\mathcal{V}_X(\mu,2)$, which reduces our problem to checking the tangency properties between an arbitrary conic and $\mathcal{Q}$. It follows that for a generic elliptic curve, any $\widetilde{\Gamma} \in \mathcal{S}\mathcal{V}_X(\mu,2)$ has two nodes (conclusion done in Appendix \textbf{B}).\\
	 	 
	 	 \indent Section §$4$ works out the (rational) degree-$2$ projection $\widetilde{S}_X \to \mathbb{P}^2(\mathbb{C})$ defined by the linear system $|\widetilde{\gamma}_X(\mu,2)|$, with discriminant a sextic $\mathcal{D}$. A divisor $\widetilde{\Gamma} \in |\widetilde{\gamma}_X(\mu,2)|$ belongs to $\mathcal{S}\mathcal{V}_X(\mu,2)$ if and only if it projects onto a line multitangent to $\mathcal{D}$. It follows that $\mathcal{S}\mathcal{V}_X(\mu,2)$ corresponds to the singular points of its dual plane curve $\mathcal{D}^\vee$. Applying Pl\"{ u}cker's formulas we obtain bounds for $\#\mathcal{S}\mathcal{V}_X(\mu,2)$, and its exact value when any $\widetilde{\Gamma} \in \mathcal{S}\mathcal{V}_X(\mu,2)$ has two nodes. Putting together both approaches we end up calculating $\#\mathcal{S}\mathcal{V}_X(\mu,2)$ for a generic elliptic curve $X$: there are $27\,\textrm{-}\,14I_0(\mu)+2I_0(\mu)^2\,\textrm{-}\,3I_1(\mu)$ such covers, where $I_0(\mu):=\#\{i,\mu_i=0\}$ and $I_1(\mu):=\#\{i,\mu_i=1\}$.\\
	 	 
	 	 \indent In the final section §$5$, in order to prove that $\#\mathcal{P}ot_X(\alpha,2) = 27$ for a generic elliptic curve, we start deducing the bound $\#\mathcal{P}ot_X(\alpha,2) \leq 27$ from the (D-G) set of equations. We then recall that if $\alpha= \mathcal{C}^{j,k}(\mu)$, each spectral data in $\mathcal{S}\mathcal{V}_X(\mu,2)$, equipped with $\xi_{j,k}$, gives rise to a potential in $\mathcal{P}ot_X(\alpha,2)$. Hence $\#\mathcal{S}\mathcal{V}_X(\mu,2)\leq \#\mathcal{P}ot_X(\alpha,2)$, implying the equality $\#\mathcal{P}ot_X(\alpha,2)=27$ whenever $\left(I_0(\mu),I_1(\mu)\right)=(0,0)$. When $\left(I_0(\mu),I_1(\mu)\right)\neq (0,0)$, we find another $27\,\textrm{-}\,\#\mathcal{S}\mathcal{V}_X(\mu,2)$ potentials in $\mathcal{P}ot_X(\alpha,2)$, deducing that for any $\alpha\in \mathbb{N}^4$ and generic elliptic curve $X$, its cardinal is equal to $\#\mathcal{P}ot_X(\alpha,2)=27$. Moreover, we calculate in both cases (at least bounds for) their spectral data. Summing up what is known for the three families $\{\mathcal{P}ot_X(d), d\leq 2\}$, we state the following conjecture: $$\#\mathcal{P}ot_X(\alpha,d)=\#\mathcal{P}ot_X(\overrightarrow{0},d)\quad\textrm{for any} \quad(\alpha,d) \in \mathbb{N}^4 \times \mathbb{N},$$ implying recursive formulas with respect to $d$ for both,  $\#\mathcal{P}ot_X(\alpha,d)$ and $\#\mathcal{S}\mathcal{V}_X(\mu,d)$.
	 	 
  \section{The family $\big\{\widetilde{S}_X,X\in\mathfrak{X}\big\}$ of rational surfaces} 
  
  {\bf 2.1}
  Let $\mathbb{P}^1$ denote the projective line over $\mathbb{C}$ (or any algebraically closed field of characteristic $0$), marked at three distinct points $\{p_1,p_2,p_3\}$, and $ \mathfrak{X}:=\mathbb{P}^1\setminus \{p_1,p_2,p_3\}$. To any  $p_0 \in \mathfrak{X}$ we associate the unique double cover $\wp:X \to \mathbb{P}^1$ with discriminant $\{p_i\}$, and denote $\omega_i:=\wp^{\textrm{-1}}(p_i)$ for any $i=0,\cdots,3$. The marked curve $(X,\omega_0)$ is an elliptic curve, and $X_2:=\{\omega_i\}$ its set of half-periods. There exits a unique indecomposable, degree $0$, rank-$2$ vector bundle over $X$. The corresponding ruled surface, say $\pi_X: S_X \to X$, can be characterized among all ruled surfaces over $X$, by the existence of a unique section $C_0 \subset S_X$ of zero self-intersection.
  
  The inverse homomorphism $[\,\textrm{-}\,1] : (X,\omega_0)\rightarrow (X,\omega_0)$ can be lifted to an involution of $S_X$, denoted $\tau$, leaving $C_0$ invariant and having two fixed points at each $\tau$-invariant fiber $S_i:=\pi_X^{\,\textrm{-}\,1}(\omega_i)$: one at $C_0$, denoted $s_i$, and the second one denoted $r_i$. We consider the blowing-up of $\big\{s_i,r_i\big\} \subset S_X$, denoted $e:S_X^\perp \to S_X$. The eight exceptional curves $\big\{s_i^\perp :=e^\ast(s_i),r_i^\perp:=e^\ast(r_i)\big\}$ are fixed by $\tau^\perp : S_X^\perp \to S_X^\perp$, the canonical lift of $\tau$ to $S_X^\perp$. There are no other fixed points, so the quotient surface $\widetilde{S}_X:= S_X^\perp/{\tau^\perp}$ is smooth. We recall hereafter the main relations between the latter surfaces and canonical morphisms.
  
  \begin{prop}\hspace*{2mm}
  	
  	Let $\varphi: S^\perp_X \to \widetilde{S}_X$ denote the canonical degree-$2$ projection and, for any $i=0,\cdots,3$, $S_i^\perp, \widetilde{S}_i,\widetilde{s}_i$ and $\widetilde{r}_i$ the strict transform of $S_i$ by the blowing-up $e$, and the images of $S_i^\perp,s_i^\perp$ and $r_i^\perp$ respectively. Then$:$
  	
  	\begin{enumerate}
  		\item for any $n>0$, $nC_0$ is the only element of the linear system $|nC_0|;$
  		\item the curves $\big\{s_i^\perp,r_i^\perp\big\}$ are fixed by $\tau^\perp$, and $\varphi: S_X^\perp \to \widetilde{S}_X$ is a double cover ramified along them$;$
  		
  		\item the composed morphism $ \pi_X \circ e:S_X^\perp \to X$ can be pushed down to $ \widetilde{S}_X \to \mathbb{P}^1$. In particular $\widetilde{S}_X$ is a smooth rational surface.

  	\end{enumerate}
  	
  \end{prop}
  \textbf{Proof.}  
  \begin{enumerate}
  	\item This result (proved in \cite{TV1}, §$3.2.1$, also) follows from the following irreducibility criterion for the ruled surface $S_X$: any divisor $\Gamma \in |nC_0+S_0|$ is equal to $aC_0 +D$, for some $ 0\leq a \leq n$ and an irreducible curve $D\in |(n\,\textrm{-}\,a)C_0+S_0|$ (cf. \cite{T1}, §$3.4$). Thus, given any $C \in |nC_0|$, there exists $ 0\leq a \leq n$ such that $C+ S_0\in |nC_0+S_0|$ is equal to $aC_0 +D$, with $D$ an irreducible curve. Hence ($D=S_0$ and) $C=nC_0$.
  	
  	\item The eight points $\{s_i,r_i\}$ being the unique fixed points of the involution $\tau$, it immediately follows that its lift $\tau^\perp$ fixes the corresponding exceptional curves, and $\varphi :S_X^\perp \to \widetilde{S}_X$ is ramified along them. Moreover, since $\tau^\perp$ has no other fixed points, the quotient surface $\widetilde{S}_X$ is smooth.
  	
  	\item  The inverse homomorphism $[\,\textrm{-}\,1]:(X,\omega_0) \to (X,\omega_0)$ also lifts to $\tau^\perp$. Taking the quotient by the latter involutions, we deduce a morphism $\widetilde{S}_X \to \mathbb{P}^1$, with generic fiber isomorphic to $\mathbb{P}^1$. In particular $\widetilde{S}_X$ is a rational surface.  $\blacksquare$
  	
  \end{enumerate}
  
  The symbol $\equiv$ will design hereafter linear equivalence between divisors of a surface. Given $C\subset S_X$ we will denote $C^\perp$ its strict transform in $S_X^\perp$, and $\widetilde{C}$ its reduced image in $\widetilde{S}_X$. For any $\nu \in \mathbb{N}^4$, we will denote $\nu^{(1)}:=\sum _i \nu_i $,  $\nu^{(2)}:=\sum _i \nu_i^2 $ and $\nu.r^\perp:= \sum_i \nu_i r_i^\perp$. The following Lemma and Proposition recall the main properties of each surface $\widetilde{S}_X$ and its set of exceptional curves (cf. \cite{TV1} §$4$ \& §$5$).
  
  \begin{lem} \hspace*{2mm} 
  	
  	Let $K_X, K_X^\perp$ and $\widetilde{K}_X$ denote the canonical divisors of $S_X,S_X^\perp$ and $\widetilde{S}_X$ respectively, $l\subset S_X$ an arbitrary fiber of $\pi_X$ and $l^\perp $ its strict transform by the blowing-up $e:S_X^\perp \to S_X$. Then$:$
  	\begin{enumerate}
  		\item the inverse images $e^\ast(K_X)$ and $\varphi^\ast(\widetilde{K}_X)$ are both linearly equivalent to\\ $K_X^\perp\,\textrm{-}\,\sum_0^3(s_i^\perp+r_i^\perp)$, while $K_X\equiv \textrm{-}\,2C_0$ and $\widetilde{K}_X\equiv \textrm{-}\,2\widetilde{C}_0\,\textrm{-}\,\sum_0^3\widetilde{s}_i;$
  		\item $\widetilde{S}_X$ is a rational anticanonical surface $(i.e.:$ $\widetilde{S}_X$ is rational and $\textrm{-}\,\widetilde{K}_X$ is effective$);$
  		\item for any $i=0,\cdots,3$ we have $\widetilde{l}:=\varphi(l^\perp) \equiv 2\widetilde{S}_i+\widetilde{r_i}+\widetilde{s}_i ;$
  		\item $Pic(\widetilde{S}_X)$ is a free group of rank $10$ and $\varphi^\ast: Pic(\widetilde{S}_X) \to Pic(S_X^\perp)$ is injective$;$
  		\item the ten divisors $\big\{\widetilde{C}_0,\widetilde{l},\widetilde{S}_i,\widetilde{r_i},i=0,\cdots,3\big\}$ make a $\mathbb{Z}$-base of $Pic(\widetilde{S}_X)$. Their intersection multiplicities are given hereafter$:$ 
  	\end{enumerate}
  	 $\widetilde{C}_0^2=\textrm{-}\,2=\widetilde{r_i}^2\,;\;\widetilde{C_0}.\widetilde{l}=1; \;\widetilde{S}_i^2=\textrm{-}\,1\,;\;\widetilde{S}_i.\widetilde{r_j}=\delta_{i,j}\,;\;\;\widetilde{l}^{\,2}=\widetilde{S}_i.\widetilde{C}_0=\widetilde{r}_i.\widetilde{C}_0=0=\widetilde{S}_i.\widetilde{l}=\widetilde{r}_i.\widetilde{l}\,.$

  \end{lem}
  
  \begin{prop} \hspace*{2mm}
  	
  	Fix $k\in \mathbb{Z}_4$ and denote $\mathbb{T}_k:=\{\nu\in \mathbb{N}^4,\nu_k+1\equiv \nu_i\, \textrm{mod.}\,2, \forall i \neq k\}$. For any exceptional curve $\widetilde{\Gamma}\subset\widetilde{S}_X$, there exists $k\in \mathbb{Z}_4$ and a unique vector $\nu \in \mathbb{T}_k$ such that $$\varphi^\ast(\widetilde{\Gamma})\equiv e^\ast (n_\nu {C}_0+S_k)\,\textrm{-}\,s_k^\perp\,\textrm{-}\,\nu .r^\perp\quad \textrm{, with}\quad n_\nu:=\frac{1}{2}(\nu^{(2)}\,\textrm{-}\,1).$$
  	
  	Conversely, given any $\nu \in \mathbb{T}_k$, there exists a unique exceptional curve, denoted $\widetilde{\Gamma}_\nu$, such that
  	$$\varphi^\ast(\widetilde{\Gamma}_\nu)\equiv e^\ast (n_\nu {C}_0+S_k)\,\textrm{-}\,s_k^\perp\,\textrm{-}\,\nu .r^\perp\;.$$
  	In particular, for any $X\in \mathfrak{X}$, the set of exceptional curves in $\widetilde{S}_X$ is parameterized by $\mathbb{T}:=\bigcup_k \mathbb{T}_k=\big\{\nu \in \mathbb{N}^4, \nu^{(1)}\equiv 1 \;\textrm{mod}.2\big\}$.\\
  \end{prop}
  
  \begin{remark}\hspace*{2mm}
  	
  	\begin{enumerate}
  		\item 
  		The homomorphism $\varphi^\ast: Pic(\widetilde{S}_X) \to Pic(S_X^\perp)$ being injective, the above formula characterizes the equivalence class of $\widetilde{\Gamma}_\nu \subset \widetilde{S}_X$.
  		
  		\item We deduce the following formula, giving the intersection multiplicity between any two exceptional curves $\widetilde{\Gamma}_\nu$ et $\widetilde{\Gamma}_\sigma$ $( \nu \in \mathbb{T}_i$ and $ \sigma \in \mathbb{T}_k)\,:$
  		
  		$$ \widetilde{\Gamma}_\nu.\widetilde{\Gamma}_\sigma=\frac{1}{4}\big((\nu\,\textrm{-}\,\sigma)^{(2)}\,\textrm{-}\,2\big) \quad \textrm{if}\quad i\neq k \quad\textrm{, }\quad \widetilde{\Gamma}_\nu.\widetilde{\Gamma}_\sigma=\frac{1}{4}\big((\nu\,\textrm{-}\,\sigma)^{(2)}\,\textrm{-}\,4\big)\quad \textrm{if}\quad i=k.$$

  		\item Let $\mathcal{S}\mathcal{V}_X(\mu,d)$ denote the Severi Variety of  $|\widetilde{\gamma}_X(\mu,d)|:=|\widetilde{\Gamma}_\mu +d(\textrm{-}\,\widetilde{K}_X)|$, pick any $\widetilde{C} \in \mathcal{S}\mathcal{V}_X(\mu,d)$ and let $f: \mathbb{P}^1 \to \widetilde{C} $ denote its normalization morphism. Then, the fiber product $C:=\varphi^\ast(\widetilde{C}) \times_{\widetilde{C}} \mathbb{P}^1 $ is a hyperelliptic curve, the inverse image of the unique point in $\widetilde{C}\cap \widetilde{s}_0$, say $p\in C$, is a smooth Weierstrass point, and the canonical marked projection $(C,p) \to X$ is a degree-$n$ $\it{ht}$ cover of arithmetic genus $g$ satisfying $2n+1=\mu^{(2)}+4d$ and $2g+1=\mu^{(1)}$.
  		
  			\item  Conversely, let $\pi : (\Gamma,p) \to X$ be a degree-$n$ $\it{ht}$ cover. Then, there exists a canonical factorization $\iota^\perp : \Gamma \to S_X^\perp $, such that $\iota^\perp(\Gamma)$ is $\iota^\perp$-invariant and projects onto a rational curve $\widetilde{\Gamma} \in |\widetilde{\gamma}_X(\mu,d)|$, for some $(\mu,d) \in \mathbb{T}_0 \times \mathbb{N}$ satisfying $2n+1=\mu^{(2)}+4d$.
  			
  			\item We thus reduce the study of all degree-$n$ $\it{ht}$ covers, and corresponding hyper-elliptic potentials, to the above linear systems and corresponding Severi sets. We recall hereafter their main properties. 
  				
  				\end{enumerate}
  		\end{remark}
  	
  	\begin{lem}\hspace*{2mm}
  		
  		For any $i=0,\cdots,3$ and $X\in \mathfrak{X}$, there is a unique isomorphism between $\widetilde{s}_i \subset \widetilde{S}_X $ and $ \mathbb{C}\cup \{\infty\}$, identifying $\widetilde{s}_i\cap \widetilde{S}_i$ and $\widetilde{s}_i\cap \widetilde{C}_0$ with $0$ and $\infty$ respectively, and, for any $\nu \in \mathbb{T}_i$, the points $\widetilde{s}_i\cap \widetilde{\Gamma}_\nu$  and $n_\nu:=\frac{1}{2}(\nu^{(2)}\,\textrm{-}\,1)$.
  		
  	\end{lem}
  \textbf{Proof.}\\
   The surface $S_X$ can be obtained from $\mathbb{P}^1 \times X$ after two suitable elementary transformations (cf. \cite{T1} §$3.6$). Choosing $(T,z\,\textrm{-}\,\omega_i)$ as pair of local coordinates at $(0,\omega_i) \in \mathbb{P}^1 \times X$, we deduce a natural isomorphism between $\widetilde{s}_i$  and $\mathbb{C}\cup \{\infty\}$: it identifies the points  $\widetilde{s_i}\cap \widetilde{C_0}$ and $\widetilde{s_i}\cap \widetilde{S}_i$, with $0$ and $\infty$ respectively. Moreover, given any exceptional curve $\widetilde{\Gamma}_\nu \subset \widetilde{S}_X$ ($\nu \in \mathbb{T}_i$), the unique point in the intersection $\widetilde{\Gamma}_\nu \cap \widetilde{s}_i$ gets identified with $n_\nu :=\frac{1}{2}( \nu^{(2)} \,\textrm{-}\,1) \in \mathbb{C}$. Let indeed $\Gamma_\nu^\perp $ denote the strict transform of $\widetilde{\Gamma}_\nu$ in $S_X^\perp$. Its direct image $\Gamma_\nu:=e_\ast(\Gamma_\nu^\perp) \subset S_X$ corresponds to the zero-divisor (in $\mathbb{P}^1 \times X$) of a so-called \textit{ tangential polynomial} of degree $n_\nu$. The latter has a local equation at $(0,\omega_i)$, in terms of the local coordinates $(T,z\,\textrm{-}\,\omega_i)$, with slope $n_\nu$ at $s_i$. Hence the result. $\quad \quad \blacksquare$
  	
  	 \begin{prop}\hspace*{2mm}
  	 	
  	 	For any $(\mu,d)\in \mathbb{T}_0 \times \mathbb{N}$ let $n:=\frac{1}{2}(\mu^{(2)}\,\textrm{-}\,1+4d)$ and $b_{n} \in \widetilde{s}_0 \simeq \mathbb{C} \cup \{\infty\}$ the point of $\widetilde{s}_0 $ identified with $n$. Then, $|\widetilde{\gamma}_X(\mu,d)|$ has the following properties $($cf. \cite{TV1}$):$
  	 	\begin{enumerate}
  	 		\item its generic element is a smooth irreducible curve of genus $d;$
  	 		\item $dim(|\widetilde{\gamma}_X(\mu,d)|)=d$ and $b_n $ is its unique base point$;$
  	 		\item $\mathcal{S}\mathcal{V}_X(\mu,d)\subset|\widetilde{\gamma}_X(\mu,d)|$ is a finite, eventually empty, subset$;$
  	 		
  	 		\item for any divisor $\widetilde{\Gamma} \in |\widetilde{\gamma}_X(\mu,d)|$ there are a unique $(\nu,e)\in \mathbb{T}_0 \times \mathbb{N}$ and an irreducible curve $\widetilde{C} \in |\widetilde{\gamma}_X(\nu,e)|$, such that$:$ 
  	 		$$\nu\,\textrm{-}\,\mu \in 2\mathbb{N}^4,\quad m:=\frac{1}{2}(\nu^{(2)}\,\textrm{-}\,1+4e)\leq n:=\frac{1}{2}(\mu^{(2)}\,\textrm{-}\,1+4d)$$ and
  	 		
  	 		$$\widetilde{\Gamma} = \widetilde{C}+ \frac{1}{2}(n\,\textrm{-}\,m)\big(2\widetilde{C}_0+\sum_i \widetilde{s}_i\big) +\sum_i \frac{1}{2}(\nu_i\,\textrm{-}\mu_i)\widetilde{r}_i \,.$$
  	 		In particular, $\widetilde{\Gamma}$ is reducible if and only if $(\nu,e) \neq (\mu, d)$.
  	 	\end{enumerate}

  	 	\end{prop}
   	
  	 \indent We fix hereafter $\mu \in \mathbb{T}_0$ and study the $\it{ht}$ covers of type $\mu$ and degree $n:=\frac{1}{2}(\mu^{(2)}\,\textrm{-}\,1+8)$, by characterizing the singularities of each curve in $\mathcal{S}\mathcal{V}_X(\mu,2)$. It turns out in this case, that the theta-characteristics $\{ \xi_{j,k} , (j,k)\in \mathbb{Z}_2 \times \mathbb{Z}_4\}$ give rise to even hyper-elliptic potentials in $\mathcal{P}ot_X(\alpha,2)$, i.e.: decomposing as 
  			$$\sum_i\alpha_i(\alpha_i+1)\wp(x\,\textrm{-}\,\omega_i)+2\wp(x\,\textrm{-}\,\rho_1)+2\wp(x+\rho_1)+2\wp(x\,\textrm{-}\,\rho_2)+2\wp(x+\rho_2) \, ,$$
  			for some $ \{\rho_1,\rho_2\}\in X^{(2)} \setminus\Delta_2$ disjoint from $X_2$. The integer vector $\alpha\in \mathbb{N}^4$ is explicitely calculated in terms of $\big(\mu,(j,k)\big)$.
  			 For example, in case $(j,k)=(1,0)$ we deduce each $\alpha_i$ from the equality: $4\alpha_i(\alpha_i+1)= (\mu^{(1)}\,\textrm{-}\,2\mu_i\,\textrm{-}\,1)(\mu^{(1)}\,\textrm{-}\,2\mu_i+1)$ (cf. \cite{TV2}). The cardinal of $\mathcal{P}ot_X(\alpha,2)$ is greater than $\#\mathcal{S}\mathcal{V}_X(\mu,2)$, and can be easily bounded from above by $27$ (see §$5.3$). We prove in Appendix \textbf{B} that for generic $X \in \mathfrak{X}$, any element in the Severi set $\mathcal{S}\mathcal{V}_X(\mu,2)$ has two simple nodes. The latter technical result enables us to find $\#\mathcal{S}\mathcal{V}_X(\mu,2)$ for any $\mu$ (see §$4$). At last, we deduce $27$ (the maximum possible number of) such potentials in §$5$, and obtain their main spectral invariants.

  \section{A global model for $\cup_\mathfrak{X}\mathcal{S}\mathcal{V}_X(\mu,2)$ via a unique surface.}
  
  Given any $\mu \in \mathbb{T}_0$ and $X \in \mathfrak{X}$, we will define hereafter a birational morphism of $\widetilde{S}_X$ onto a so-called weak del Pezzo surface, naturally endowed with a degree-$2$ projection onto the projective plane. The latter turns out to be independent of $X \in \mathfrak{X}$ (see Appendix\textbf{A}) and gives a common frame for studying the family of Severi sets $\cup_\mathfrak{X}\mathcal{S}\mathcal{V}_X(\mu,2)$.
  \begin{defi}
  Let $\psi_\mu : \widetilde{S}_X \to \underline{S}_\mu$ denote the contraction of the divisor $\widetilde{\Gamma}_\mu+\widetilde{s}_0$, $\underline{K}_\mu$ its canonical divisor and $\underline{L}_\mu:=\textrm{-}\underline{K}_\mu$. Given any curve $\widetilde{C}\subset \widetilde{S}_X$, we will denote hereafter $\underline{C}:=\psi_\mu (\widetilde{C})$.
   \end{defi}
 The following properties follow directly from those satisfied by $\widetilde{S}_X$. We prove in Appendix \textbf{A} a more fundamental one, namely, that $\underline{C}_0$ and $\underline{S}_\mu$ are independent of $X \in \mathfrak{X}$. 
 
 \begin{prop}\hspace*{2mm}
 	
 	\begin{enumerate}
 		\item $\underline{S}_\mu$ is rational and $\underline{L}_\mu=\textrm{- }\psi_\mu (\widetilde{K}_X)=2\underline{C}_0+\sum_1^3 \underline{s}_j;$
 		
 	\item for any $j=1,2,3$ the linear system $|\underline{L}_\mu\,\textrm{-}\,\underline{s}_j|$ is a pencil of generically rational irreducible curves of zero self-intersection $($see also §$3.5.2))$;
 		
 \item $\underline{L}_\mu$ has self-intersection $2$, $|\underline{L}_\mu|$ has no base point and $dim|\underline{L}_\mu|=2;$
 
 \item $\underline{L}_\mu$ defines a degree-$2$ morphism, denoted $g_\mu: \underline{S}_\mu \to |\underline{L}_\mu|^\vee \simeq \mathbb{P}^2(\mathbb{C})$.\end{enumerate} In other words, $\underline{S}_\mu$ is a so-called $($weak$)$ del Pezzo surface of degree $2$. The corresponding so-called Geiser involution, will be denoted $\iota_\mu$ $($cf. \cite{Do}, § $8.7.1)$.
  \end{prop}

  \indent  For any  $\underline{p} \in \underline{S}_\mu$, let $|\underline{L}_\mu|_{\underline{p}}$ denote the pencil in $|\underline{L}_\mu|$, made of all divisors going through $\underline{p}$. Let also $\underline{p}_X \in \underline{C}_0$ denote the point corresponding to $X\in \mathfrak{X}$, under the natural identifications $\mathfrak{X}\simeq \widetilde{C}_0 \setminus \bigcup_1^3 \widetilde{s}_j \simeq\underline{C}_0 \setminus \bigcup_1^3 \underline{s}_j$. 
  The following lemmas collect their basic properties, some of which depend upon $I_\mu(0):=\#\big\{i,\,\mu_i=0\big\}$. 
  
  \begin{lem}\hspace*{2mm}
  	
  	\begin{enumerate}
  		
  		\item The projection $\psi_\mu : \widetilde{S}_X \to \underline{S}_\mu$ identifies $|\widetilde{\gamma}_X(\mu,1) |$ with $|\underline{L}_\mu|_{\underline{p}_X}$.
  		
  		\item The pencil $ |\underline{L}_\mu|_{\underline{p}_X}$ has $1+I_\mu(0)$ reducible elements, namely, \\
  		$ 2\underline{C}_0+ \sum_1^3 \underline{s}_j$ and, for each index $k$ such that $\mu_k=0$, $\underline{\Gamma}_\nu+\underline{r}_k$, with $\nu_i=\mu_i + 2\delta_{i,k}$.
  		
  		\item The generic element $ \underline{\Gamma} \in |\underline{L}_\mu|$ is a smooth irreducible curve of genus $1$, double cover of the line $g_\mu(\underline{\Gamma} ) \subset\mathbb{P}^2(\mathbb{C})$.
  		
  	\end{enumerate}
  \end{lem}
  
  \textbf{Proof.}
  \begin{enumerate}
  	\item The pencil $|\widetilde{\gamma}_X (\mu,1) |$ gets embedded by $\psi_\ast$ into  $|\underline{L}_\mu|$, and all its elements intersect $\underline{C}_0$ at $\underline{p}_X $. Hence the claim.
  	
  	\item According to the preceding result, any reducible element of $|\underline{L}_\mu|_{\underline{p}_X }$ comes from $|\widetilde{\gamma}_X(\mu,1)|$. The latter are known (cf. \cite{TV1} §5.2.4)).
  	
  	\item Recall that for any $\underline{p}_X\in \underline{C}_0 \setminus\bigcup_1^3 \underline{s}_j\simeq \mathfrak{X}$, the pencil $|\underline{L}_\mu|_{\underline{p}_X}$ is equal to ${\psi_\mu}_\ast(|\widetilde{\gamma}_X(\mu,1)|)$, and its generic element is a smooth irreducible curve of genus $1$. Let $\underline{C} \in |\underline{L}_\mu|_{\underline{p}_X}$ and take any line $H \subset \mathbb{P}^2(\mathbb{C})$. We have $$2=\underline{C}.\underline{L}_\mu=\underline{C}.g_\mu^\ast(H)={g_\mu}_\ast(\underline{C}).H\;.$$ We claim that $\underline{C}$ is a double cover of its projection, i.e.: ${g_\mu}_\ast(\underline{C})=2g_\mu(\underline{C})$ (and $g_\mu(\underline{C})$ is a line). Otherwise, ${g_\mu}_\ast(\underline{C})=g_\mu(\underline{C})$, implying that its inverse image decomposes as:
  	$$g_\mu^\ast\big(g_\mu(\underline{C})\big)=\underline{C}+\iota_\mu(\underline{C})\;.$$
  	But the reducible elements in $|\underline{L}_\mu|_{\underline{p}_X}$ are $2\underline{C}_0+\sum_{1}^3 \underline{s}_j $ and, for each index $k\in \mathbb{Z}_4$ such that $\mu_k=0$,  $\underline{\Gamma}_{\nu}+ \underline{r}_k$ (where $\nu_i=\mu_i+2\delta_{i,k}$ for all $i\in \mathbb{Z}_4$). Contradiction!	$\blacksquare$
  \end{enumerate}
  
  \begin{lem}\hspace*{2mm}
  	
  	\begin{enumerate}
  		\item The set $\big\{\underline{\Gamma}_\nu, (\nu\,\textrm{-}\,\mu)^{(2)}=2 \big\}\cup\{\underline{C}_0\}$ contains all the exceptional curves in $\underline{S}_\mu$. There are, therefore, $25$, $19$, $14$ or $10$, if $I_\mu(0)$ is equal to $0$, $1$, $2$ or $3$.
  		
  		\item The three curves $\{\underline{s}_1,\underline{s}_2,\underline{s}_3\}$ are $\iota_\mu$-invariant and project onto three distinct points of the line $g_\mu(\underline{C}_0)$.  
  		
  		\item For any $\nu \in \mathbb{T}_j$, $g_\mu(\underline{\Gamma}_\nu)$ is a line going through $g_\mu(\underline{s}_j)$.
  		
  		\item The involution $\iota_\mu$ fixes $\underline{C}_0$, and permutes the other exceptional curves. More precisely$:$ $$\iota_\mu(\underline{\Gamma}_\nu)=\underline{\Gamma}_{\iota_\mu(\nu)} \quad\textrm{where}\quad \forall i, \iota_\mu(\nu)_i=|2\mu_i\,\textrm{-}\,\nu_i|\;,$$ 
  		and $$\underline{\Gamma}_\nu+ \underline{\Gamma}_{\iota_\mu(\nu)}+ \underline{s}_j+ \sum_{i,\mu_i=0} \underline{r}_i \equiv \underline{L}_\mu \quad\textrm{if}\;\;\nu \in \mathbb{T}_j.$$
  		
  		\item We have $\iota_\mu(\underline{\Gamma}_\nu)=  \underline{\Gamma}_\nu$ $($i.e.$:$ $\iota_\mu(\nu)=\nu)$ if, and only if, there exists $j\neq k$ such that $\mu_j=0=\mu_k$ and $\nu=\mu + \epsilon$, where $\epsilon_i=1$ if $i\in \{j,k\}$ and $\epsilon_i=0$ otherwise. In the latter case $\iota_\mu$ fixes $\underline{\Gamma}_\nu$ and
  		$$2\underline{\Gamma}_\nu+\underline{r_j}+\underline{r_k}\equiv 2\underline{C_0}+\underline{s_j}+\underline{s_k}\in |\underline{L}_\mu \,\textrm{-}\,\underline{s}_l|.$$ 
  		For example, if $\mu_1=0=\mu_2$ we have $\nu= \mu +(0,1,1,0)$.
  		
  		\item The involution $\iota_\mu$ acts non-trivially on each $(\textrm{-}2)$-curve  $\{\underline{s}_j\}\cup\{\underline{r}_i,\mu_i=0\}$, $($i.e.$:$ it only fixes two points on each one of them$)$. In particular, on each $\underline{s}_j$, identified as precedently with $\mathbb{C}\cup\{\infty\}$, it coincides with the simmetry with center $\frac{1}{2}(\mu^{(2)}+1))$, $$t\in \mathbb{C} \mapsto \mu^{(2)}+1\,\textrm{-}\,t \in \mathbb{C},\quad\textrm{and}\quad \infty \mapsto \infty\,.$$

  	\end{enumerate}
  \end{lem}
  
  \textbf{Proof.}
  \begin{enumerate}
  	
  	\item An exceptional curve of $\underline{S}_\mu$ is, either the image of an exceptional curve $\widetilde{\Gamma}_\nu \subset\widetilde{S}_X$ disjoint from $\widetilde{\Gamma}_\mu$ and $\widetilde{s}_0$, which is equivalent to asking $(\nu\,\textrm{-}\,\mu)^{(2)}=2$, or the image of a $(\textrm{-}2)$-curve  intersecting $\widetilde{s}_0$, only satisfied by $\widetilde{C}_0$. As for their number, if $I_\mu(0)=0$ there are $24$ solutions to the equation $(\nu\,\textrm{-}\,\mu)^{(2)}=2$. The other cases also follow directly.
  	
  	\item The equalities $1=\underline{C}_0.\underline{L}_\mu=\underline{C}_0.g_\mu^\ast(H)={g_\mu}_\ast(\underline{C}_0).H$ imply that ${g_\mu}_\ast(\underline{C}_0)$ is a line equal to ${g_\mu}(\underline{C}_0)$. Analogously, the equalities $0=\underline{s}_j.\underline{L}_\mu=\underline{s}_j.g_\mu^\ast(H)={g_\mu}_\ast(\underline{s}_j).H$ imply that ${g_\mu}(\underline{s}_j)$ is a point, for any $j=1,2,3$. The latter are distinct points of the line ${g_\mu}(\underline{C}_0)$. Indeed, $\underline{C}_0$ intersects each $\underline{s}_j$, so ${g_\mu}(\underline{C}_0)$ goes through ${g_\mu}(\underline{s}_j)$.\\
  	\indent On the other hand, had we, for example, ${g_\mu}(\underline{s}_1)={g_\mu}(\underline{s}_2)$, the subspace\\ $\underline{s}_1+\underline{s}_2+|\underline{L}_\mu\,\textrm{-}\,\underline{s}_1\,\textrm{-}\,\underline{s}_2|\subset |\underline{L}_\mu|$ would be a pencil. But $2\underline{C}_0+\underline{s}_3$ is the unique element of $|\underline{L}_\mu\,\textrm{-}\,\underline{s}_1\,\textrm{-}\,\underline{s}_2|$. Contradiction!
  	
  	\item Every exceptional curve $\underline{\Gamma}_\nu$ satisfies $\underline{\Gamma}_\nu.\underline{K}_\mu= \,\textrm{-}\,1$, hence
  	$$1=\underline{\Gamma}_\nu.\underline{L}_\mu=\underline{\Gamma}_\nu.g_\mu^\ast(H)={g_\mu}_\ast(\underline{\Gamma}_\nu).H\;.$$
  	Therefore, ${g_\mu}_\ast(\underline{\Gamma}_\nu)={g_\mu}(\underline{\Gamma}_\nu)$, implying it is a line. Moreover, if $\nu \in \mathbb{T}_j$, it must go through $g_\mu(\underline{s}_j)$ because $\underline{\Gamma}_\nu$ intersects $\underline{s}_j$.
  	
  	\item According to §$3.4.2)$ and §$3.4.3)$, for any $j=1,2,3$ the involution $\iota_\mu$ permutes the exceptional curves intersecting $\underline{s}_j$. On the other hand, $\underline{C}_0 $ being the unique exceptional curve intersecting each $\underline{s}_j$, it must be (at least) invariant by $\iota_\mu$. Hence, $\iota_\mu$ fixes the three points $\underline{C}_0 \cap \bigcup_j \underline{s}_j$, as well as $\underline{C}_0$ because it is isomorphic to $\mathbb{P}^1$. Let $\nu \in \mathbb{T}_j$ satisfy  $ \;\big(\nu\,\textrm{-}\,\mu\big)^{(2)}=2$, so that $\underline{\Gamma}_\nu \subset \underline{S}_\mu$ is an exceptional curve, and let us denote $\iota_\mu(\nu)\in \mathbb{N}^4$ the vector defined by $\iota_\mu(\nu)_i=|2\mu_i\,\textrm{-}\,\nu_i|$, for any $i=0,\cdots,3$. We can also check that $|\nu_i\,\textrm{-}\,\mu_i|=|\iota_\mu(\nu)_i\,\textrm{-}\,\mu_i|$, hence $\nu_i\equiv \iota_\mu(\nu)_i \,\textrm{mod.}2$, for any $i$, as well as  $(\iota_\mu(\nu)\,\textrm{-}\,\mu)^{(2)}=2$ and $\iota_\mu(\nu) \in \mathbb{T}_j$. In other words, $\underline{\Gamma}_{\iota_\mu(\nu)}$ is an exceptional curve of $\underline{S}_\mu$ intersecting $\underline{s}_j$. The following formula follows: $\nu_i +\iota_\mu(\nu)_i=2$ if $\mu_i=0$ and $\nu_i=1$ (while $\nu_i +\iota_\mu(\nu)_i=2\mu_i$ otherwise). We have two possibilities:
  	\begin{enumerate}
  		\item either $\nu_i +\iota_\mu(\nu)_i=2\mu_i$ for each $i$, implying $\iota_\mu(\nu)\,\textrm{-}\,\nu=2(\mu\,\textrm{-}\,\nu)$ and $\iota_\mu(\nu)\,\textrm{-}\,\mu=\mu\,\textrm{-}\,\nu$, hence $\underline{\Gamma}_{\iota_\mu(\nu)}.\underline{\Gamma}_\nu=1$ because $(\iota_\mu(\nu)\,\textrm{-}\,\nu)^{(2)}=8$ (cf. §$2$ Remark  $1.2)$). It then follows that $g_\mu(\underline{\Gamma}_\nu)=g_\mu(\underline{\Gamma}_{\iota_\mu(\nu)})$, because they both go through the point $g_\mu(\underline{s}_j)$ (cf. §$3.4.3)$) and the projection of $\underline{\Gamma}_{\iota_\mu(\nu)}\cap \underline{\Gamma}_\nu\,;$
  		
  		\item or there exists $i_0$ such that $\mu_{i_0}=0$ and $\nu_{i_0}=1$, in which case $\underline{\Gamma}_{\iota_\mu(\nu)}$ and $\underline{\Gamma}_\nu$ intersect the fibers $ \underline{s}_j$ and $\underline{r}_{i_0}$. Their projections being lines with two common points, they must coincide.
  	\end{enumerate}
  	Since $\underline{\Gamma}_{\iota_\mu(\nu)}$ and $\underline{\Gamma}_\nu$ project onto the same line, they are exchanged by the Geiser involution. On the other hand, $\underline{\Gamma}_\nu$ intersects $\underline{s}_j$, as well as $\underline{r}_i$ if $\mu_i=0$. Therefore, $\underline{\Gamma}_\nu+ \underline{\Gamma}_{\iota_\mu(\nu)}+\underline{s}_j+ \sum_{i,\mu_i=0} \underline{r}_i$ is a component of $g_\mu^\ast(g_\mu(\underline{\Gamma}_\nu))$. We can check at last (e.g.: by means of §$2.3$) that they are both numerically equivalent to $\underline{L}_\mu$. Thus, they coincide.
  	
  	\item An exceptional curve $\underline{\Gamma}_\nu $ is invariant by $\iota_\mu$, i.e.: $\iota_\mu(\underline{\Gamma}_\nu ) =\underline{\Gamma}_{\iota_\mu(\nu)}=\underline{\Gamma}_\nu $, if and only if $\iota_\mu(\nu)=\nu$. Equivalently, according to formula §$3.4.4)$, there exist $j\neq k$ such that $\mu_j=\mu_k=0$, but also $\nu_j=\nu_k=1$. In particular, the three curves $\underline{r}_j,\underline{r}_k,\underline{s}_l$ project onto points because $\underline{r}_j.\underline{L}_\mu = \underline{r}_k.\underline{L}_\mu= \underline{s}_l.\underline{L}_\mu=0$. Moreover, $\underline{\Gamma}_\nu.\underline{r}_j=\underline{\Gamma}_\nu.\underline{r}_k=\underline{\Gamma}_\nu.\underline{s}_l=1$, hence $\underline{\Gamma}_\nu$ projects with degree $1$ onto a line going through the three points $\{g_\mu(\underline{r}_j),g_\mu(\underline{r}_k),g_\mu(\underline{s}_l)\}$. Let us consider, on the other hand, an exceptional curve $\underline{C}$ intersecting two out of the three curves $\{\underline{r}_j,\underline{r}_k,\underline{s}_l\}$. Had both (-$2$)-curves the same projection in $g_\mu(\underline{C})$, the latter would have a node there. Contradiction!\\
  	It follows that the three points $\underline{\Gamma}_\nu\cap \big(\underline{r}_j \cup \underline{r}_k \cup \underline{s}_l)$ of $\underline{\Gamma}_\nu$ are fixed by the involution, hence $\underline{\Gamma}_\nu$ as well, the latter being isomorphic to $\mathbb{P}^1$.
  	
  	At last, knowing that $g_\mu(\underline{\Gamma}_\nu)$ contains $\{g_\mu(\underline{r}_j),g_\mu(\underline{r}_k),g_\mu(\underline{s}_l)\}$, we deduce immediately the following formula:
  	$$g_\mu^\ast(g_\mu(\underline{\Gamma}_\nu))=2\underline{\Gamma}_\nu + \underline{r}_j+\underline{r}_k+\underline{s}_l\equiv 2\underline{C}_0+\sum_{h=1}^3 \underline{s}_h.$$ 
  	
  	\item For any $j\neq 0$, $\underline{s}_j$ is the unique $(\textrm{-}2)$-curve intersecting each exceptional curve $\underline{\Gamma}_\nu$, with $\nu \in \mathbb{T}_j$. The involution $\iota_\mu$ permutes these curves (cf. §$3.4.4)$), hence leaves stable $\underline{s}_j$ and exchanges the points $\underline{\Gamma}_\nu\cap \underline{s}_j$ and $\underline{\Gamma}_{\iota_\mu(\nu)}\cap\underline{s}_j$. The latter correspond to the values $n_\mu$ and $n_{\iota_\mu (\nu)}$ (cf. §$2.4$), which satisfy the equalities $n_\nu+n_{\iota_\mu (\nu)}= 2n_\mu+2=\mu^{(2)}+1$. Therefore, they are symmetric with respect to $\frac{1}{2}(\mu^{(2)}+1)$.\\
  	Assume at last there exists $i$ such that $\mu_i=0$, pick any $i' $ for which $\mu_{i'}\neq 0$ and let $\nu \in \mathbb{N}^4$ denote the vector defined as $\nu_k=\mu_k + 1$ if $k \in \{i,i'\}$, and $\nu_k=\mu_k$ otherwise. For example, in case $\mu=(\mu_0,0,\mu_2,\mu_3)$, $\mu_0$ is odd, hence $ \neq 0$ and $\nu=(\mu_0+1,1,\mu_2,\mu_3)$. By definition, $\iota_\mu(\nu)=(\mu_0\,\textrm{-}\,1,1,\mu_2,\mu_3)$. It follows that the exceptional curves $\underline{\Gamma}_\nu$ and $\underline{\Gamma}_{\iota_\mu(\nu)}$ are disjoint, exchanged by $\iota_\mu$ and $\underline{\Gamma}_\nu.\underline{r}_i=1=\underline{\Gamma}_{\iota_\mu(\nu)}.\underline{r}_i$. Thus, $\iota_\mu$ exchanges their points of intersection with $\underline{r}_i$, proving that it acts non-trivially on $\underline{r}_i$.  							$\blacksquare$
  \end{enumerate}

  \begin{prop}\hspace*{2mm}
  	
  	\begin{enumerate}
  		
  		\item Let $\underline{C} \subset \underline{S}_\mu$ be an irreducible curve. Then $\underline{C}.\underline{L}_\mu=0$ if and only if it has self-intersection $\underline{C}.\underline{C}=\textrm{-}\,2$, in which case $g_\mu(\underline{C})$ is a point. The latter is equivalent to asking $\underline{C}\in  \big\{\underline{s}_j,j\neq 0\big\}\bigcup \big\{\underline{r}_i,\mu_i= 0\big\} $.

  		\item  For any curve $\underline{C}\in  \big\{\underline{s}_j,j\neq 0\big\}\bigcup \big\{\underline{r}_i,\mu_i= 0\big\} $, $dim|\underline{L}_\mu \,\textrm{-}\,\, \underline{C}|=1$. Moreover, its generic element is irreducible, isomorphic to $\mathbb{P}^1$ and projects with degree $2$ onto a line going through the point $g_\mu(\underline{C})$. In other words, $g_\mu(|\underline{L}_\mu \,\textrm{-}\,\, \underline{C}|)$ is the pencil of lines going through that point.
  		
  		\item There are $3+I_\mu(0)$ positive dimensional fibers of $g_\mu: \underline{S}_\mu \to \mathbb{P}^2(\mathbb{C})$. They all belong to $\big\{\underline{s}_j,j\neq 0\big\}\bigcup \big\{\underline{r}_i,\mu_i= 0\big\} $. 
  		
  		\item The corresponding pencils $\underline{C}+|\underline{L}_\mu \,\textrm{-}\,\underline{C}|$ contain all reducible curves in $|\underline{L}_\mu| $.
  		
  	\end{enumerate}
  \end{prop}
  
  \textbf{Proof.}
  \begin{enumerate}
  	
  	\item Assume $\underline{C}\subset \underline{S}_\mu$ is irreducible and satisfies $\underline{C}.\underline{L}_\mu=0$. We know that $g_\mu(\underline{C})$ is a point because ${g_\mu}_\ast(\underline{C}).H=\underline{C}.g_\mu^\ast(H)=\underline{C}.\underline{L}_\mu=0$, hence is contained in a positive dimensional fiber of $g_\mu$. In particular $dim|\underline{C}|=0$ and $p_a(\underline{C})\geq 0$, implying $\underline{C} .\underline{C}$ is even, either -$2$ or $0$. We immediately deduce  $\underline{C} .\underline{C}=\,\textrm{-}2$, by means of \cite{H} III.1.c), which means $\underline{C} \in \big\{\underline{s}_j,j\neq 0\big\}\bigcup \big\{\underline{r}_i,\mu_i= 0\big\}$. 
  	
  	\item Each element of $|\underline{L}_\mu\,\textrm{-}\,\underline{C}|$ intersects $\underline{C}$ and projects with degree $2$ over a line going through $g_\mu(\underline{C})$. According to §$3.3.3)$), $\underline{C}+|\underline{L}_\mu\,\textrm{-}\,\underline{C}| \subsetneq |\underline{L}_\mu|$, hence  $dim|\underline{L}_\mu\,\textrm{-}\,\underline{C}| < dim|\underline{L}_\mu|=2$. On the other side, the equality $(\underline{L}_\mu\,\textrm{-}\,\underline{C}).\underline{K}_\mu=\textrm{-}\,2$ alone implies $|\underline{L}_\mu\,\textrm{-}\,\underline{C}|$ is a pencil without base component or base point, according to \cite{H} III.1-a). Therefore, its generic element $\underline{D}\in  |\underline{L}_\mu\,\textrm{-}\,\underline{C}|$ is smooth, has zero self-intersection and projects with degree $2$ over a line going through the point $g_\mu(\underline{C})$. Should $D$ decompose, we would have $\underline{D}=\underline{D}'+\underline{D}''$, with $\underline{D}'$ and $\underline{D}''$ disjoint curves exchanged by the involution $\iota_\mu$, hence isomorphic to $\mathbb{P}^1$ and satisfying: 
  	
  		$$\underline{D}'.\underline{D}''=0,\; (\underline{D}'+\underline{D}'')^2=\underline{D}'^2+\underline{D}''^2=(\underline{L}_\mu\,\textrm{-}\,\underline{C})^2=0\;\;\textrm{and}\;\;(\underline{D}'+\underline{D}'').\underline{L}_\mu=2.$$
  		We would get, at last
  		 $$0=p_a(\underline{D}')+p_a(\underline{D}'')=2+ \frac{1}{2}\big(\underline{D}'^2+\underline{D}''^2\,\textrm{-}\,(\underline{D}'+\underline{D}'').\underline{L}_\mu\big)=2+\frac{1}{2}(0\,\textrm{-}\,2)
  		=1.$$
  		Contradiction!

  	\item For any $\underline{C}'\neq \underline{C}''  \in  \big\{\underline{s}_j,j\neq 0\big\}\bigcup \big\{\underline{r}_i,\mu_i= 0\big\} $, we have just proved $\underline{C}'+|\underline{L}_\mu\,\textrm{-}\,\underline{C}'|$ and $\underline{C}''+|\underline{L}_\mu\,\textrm{-}\,\underline{C}''|$ are distinct, hence $g_{\mu}(\underline{C}')\neq g_{\mu}(\underline{C}'')$. In particular $\underline{C}'$ and $\underline{C}''$ are contained in different fibers of $g_\mu$. Conversely, if $\underline{C}\subset \underline{S}_\mu$ is an irreducible component of a fiber, it must satisfy $\underline{C}. \underline{L}_\mu=0$, and therefore $\underline{C}. \underline{C}=\,\textrm{-}\,2$ according to §$3.5.1)$. In other words, the fibers of positive dimension are the $(\textrm{-}\,2$)-curves. There are $3+ I_\mu(0)$.
  	
  	\item Let $D=\underline{D}_1+\underline{D}_2 \in |\underline{L}_\mu|$ be reducible, with $\underline{D}_1$ irreducible and $\underline{D}_2$ a non-zero effective divisor. Then, $\underline{L}_\mu$ being \textit{nef}, we have $0\leq \underline{D}_1.\underline{L}_\mu  \leq (\underline{D}_1+\underline{D}_2).\underline{L}_\mu  =2$, implying:
  	\begin{enumerate}
  		\item if $\underline{D}_1.\underline{L}_\mu =0$, we have just proved that $\underline{D}_1 \in \big\{\underline{s}_j,j\neq 0\big\}\bigcup \big\{\underline{r}_i,\mu_i= 0\big\} $;
  		\item if $\underline{D}_1.\underline{L}_\mu =1$, we immediately check $\underline{D}_1\neq \underline{C}_0$, implying the existence of a unique index $j_0\neq 0$ such that $\underline{D}_1.\underline{s}_j=\delta_{j,j_0}$, and
  		$$ 0=\underline{L}_\mu.\underline{s}_{j_0}=(\underline{D}_1+\underline{D}_2).\underline{s}_{j_0}= 1+\underline{D}_2.\underline{s}_{j_0}.$$
  		But $\underline{D}_2.\underline{s}_{j_0}= \textrm{-}\,1$ can only happen if $\underline{s}_{j_0}$ is a component of $\underline{D}_2$.
  		
  		\item If $\underline{D}_1.\underline{L}_\mu =2$ we have $\underline{D}_2.\underline{L}_\mu =0$ and, $\underline{L}_\mu$ being \textit{nef}, each irreducible component $\underline{C}$ of $\underline{D}_2$ satifies $\underline{C}.\underline{L}_\mu=0$. At last, according to §$3.5.1)$, $\underline{C}$ must be one of the curves $\big\{\underline{s}_j,j\neq 0\big\}\bigcup \big\{\underline{r}_i,\mu_i= 0\big\}$. $\quad \blacksquare$\\
  	\end{enumerate}
  	
  \end{enumerate}

  For any $\underline{p}_0 \in \underline{S}_\mu$ outside $\big\{\underline{s}_j,j\neq 0\big\}\bigcup \big\{\underline{r}_i,\mu_i= 0\big\}$, its image $\iota_\mu(\underline{p}_0)$ is the other fixed point in $ |\underline{L}_\mu|_{\underline{p}_0}$. We describe hereafter the action of $\iota_\mu$ on the set of curves of negative self-intersection of $\underline{S}_\mu$. Recall (for any $j=1,2,3$) the natural isomorphism $\underline{s}_j \to \mathbb{C}\cup \{\infty\}$ (§$2.4$), identifying $\underline{s}_j \cap \underline{C}_0$ and $\underline{s}_j \cap \underline{\Gamma}_{\nu} $ with $\infty$ and $n_\nu:=\frac{1}{2}(\nu^{(2)}\,\textrm{-}\,1)$, respectively.

  \begin{defi}\hspace*{2mm}
  	
  	Let $\underline{R}$ denote the ramification divisor of the double cover $g_\mu: \underline{S}_\mu \to \mathbb{P}^2(\mathbb{C})$, $\mathfrak{D}: =g_\mu(\underline{R})$ its discriminant, $\underline{R}_r \subset  \underline{S}_\mu $ the reduced inverse image of $\mathfrak{D}$ and $\underline{R}_c:=\underline{R}\,\textrm{-}\, \underline{R}_r$.
  \end{defi}
  We prove hereafter that $\underline{R}$, $\underline{R}_r$ and $\underline{R}_c$ are reducible, the latter being the sum of positive dimensional fibers of $g_\mu$. We also compute their irreducible components.

  \begin{prop}\hspace*{2mm}
  	
  	\begin{enumerate}
  		\item $\underline{R}$ is numerically equivalent to $2\underline{L}_\mu\equiv 4\underline{C}_0+2\sum_{1}^3 \underline{s}_j$, while $\underline{R}_r$ is reducible, contains $\underline{C}_0$, and its image $\mathfrak{D}: =g_\mu(\underline{R})$ is a reduced quartic.
  		
  		\item $\underline{R}_c:=\underline{R}\,\textrm{-}\, \underline{R}_r$ is equal to $\underline{R}_c=\sum_{1}^3 \underline{s}_j+\sum_{i,\mu_i=0} \underline{r}_i$. In particular the curve $\underline{R}_r\,\textrm{-}\,\underline{C}_0$, henceforth denoted $\underline{\Omega}$, is numerically equivalent to
  		
  		 $$\underline{\Omega} \equiv 3\underline{C}_0+\sum_{1}^3 \underline{s}_j\,\textrm{-}\sum_{i,\mu_i=0} \underline{r}_i\;.$$
  		
  		\item $\underline{\Omega}$ is an irreducible curve if and only if $I_\mu(0) \leq 1$.
  		
  	\end{enumerate}
  \end{prop}
  
  \textbf{Proof.}
  \begin{enumerate}
  	\item Let$H$ denote hereafter the class of a projective line. By its very definition, $\underline{R}$ satisfies the numerical equivalence $g_\mu^\ast(\textrm{-}3H)+\underline{R}\equiv \underline{K}_\mu$, implying $\underline{R}\equiv 2\underline{L}_\mu$ because $g_\mu^\ast(H)\equiv \underline{L}_\mu= \,\textrm{-}\,\underline{K}_\mu$, and its image ${g_\mu}_\ast(\underline{R})\equiv 2{g_\mu}_\ast(\underline{L}_\mu)\equiv 4H$ is a plane quartic. In particular $|\underline{L}_\mu|$ intersects $\underline{R}$, hence $\underline{R}_r$, with multiplicity $4$. According to §$3.3.3)$, its generic element is a smooth curve of genus $1$, projecting with degree $2$ over a line. Therefore, it must have four (distinct) ramification points, namely: its intersection with $\underline{R}_r$, forcing $\underline{R}_r$ to be reducible. We also know that $\underline{C}_0$ is fixed by the involution $\iota_\mu$ (cf. §$3.4.4)$), therefore contained in $\underline{R}_r$.
  	
  	\item The latter result implies the equation $$\underline{\Omega} +\sum_{j=1}^3 \alpha_j \underline{s}_j+ \sum_{i,\mu_i=0} \beta_i \underline{r}_i \equiv 3\underline{C}_0+2\sum_{j=1}^3 \underline{s}_j\;,$$ with non-negative unknowns $\big\{\alpha_j,\;\beta_i\big\}$. Intersecting both sides with $\underline{s}_j$, and $\underline{C}_0$, we deduce:
  	$$\underline{\Omega}.\underline{s}_j\,\textrm{-}\,2\alpha_j=3\,\textrm{-}\,4=\,\textrm{-}\,1\;\textrm{,}\; \underline{\Omega}.\underline{C}_0+\sum_{j=1}^3 \alpha_j=\,\textrm{-}\,3+6=3\;,$$
  	
  	from which follows $\alpha_j=1$ for any $j=1,2,3$, because $\underline{\Omega}.\underline{s}_j \geq 0$ and $\underline{\Omega}.\underline{C}_0 \geq 0$.\\
  	It remains to prove that $\beta_{i_0}=1$, whenever $\mu_{i_0}=0$. In the latter case, we choose $i'_0$ such that $\mu_{i'_0}\neq 0$ and define $\nu \in \mathbb{N}^4$, with $\nu_k=\mu_k + 1$ if $k\in \{i_0,i'_0\}$, and $\nu_k=\mu_k $ otherwise. Recall that $\underline{\Gamma}_\nu$ and $\iota_\mu({\underline{\Gamma}_\nu})$ are disjoint (cf. the proof of §$3.4.6$)), hence disjoint with $\underline{R}_r=\underline{\Omega}+\underline{C}_0$. Furthermore, $\underline{\Gamma}_\nu$ satisfies 
  	$$\underline{\Gamma}_\nu.\sum_j \underline{s}_j=1=\underline{\Gamma}_\nu.\underline{r}_{i_0}\quad\textrm{but}\quad \underline{\Gamma}_\nu.\underline{r}_l=0 \;\;\textrm{if}\;\; \mu_l=0.$$
  	
  	Intersecting now with $\underline{\Gamma}_\nu$ gives $\beta_i=1$, hence
  	$\underline{\Omega} \equiv 3\underline{C}_0+\sum_{1}^3 \underline{s}_j\,\textrm{-}\sum_{i,\mu_i=0} \underline{r}_i\;$.
  	
  	\item We now know $\underline{\Omega}:=\underline{R}_r\,\textrm{-}\,\underline{C}_0$ projects with degree $1$ onto the reduced cubic $\mathcal{K}_0:=g_\mu\big(\underline{\Omega}\big)\equiv g_\mu\big(3\underline{C}_0\big)\equiv 3H$. The latter is therefore irreducible if it contains no line. Suppose on the contrary, there exists a line $D \subset \mathcal{K}_0$, distinct of $g_\mu(\underline{C}_0)$ because $\underline{R}_r$ is a reduced divisor, and assume $I_\mu(0) \leq 1$. Its reduced inverse image (by $g_\mu$), denoted $\underline{\delta}$, is $\iota_\mu$-invariant and, modulo positive dimensional fibers, $2\underline{\delta}\equiv \underline{L}_\mu$. In other words: $2\underline{\delta}+\sum_{j=1}^3\alpha_j \underline{s}_j \equiv \underline{L}_\mu$ if $I_\mu(0)=0$, or $2\underline{\delta}+\sum_{j=1}^3\alpha_j \underline{s}_j +\beta \underline{r}_k\equiv \underline{L}_\mu$ if $k$ is the unique index for which $\mu_k=0$. Intersecting with each $\underline{s}_j$ (plus $\underline{r}_k$ if $\mu_k=0$) we obtain $\alpha_j= \underline{\delta}.\underline{s}_j $ (and $\beta=\underline{\delta}.\underline{r}_k$), hence
  	
  	$$2=\underline{L}_\mu^2=(2\underline{\delta}+\sum_{j=1}^3\alpha_j \underline{s}_j +\beta \underline{r}_k)^2=4\underline{\delta}^2+2\big(\sum_{j=1}^3 \alpha_j^2 + \beta^2\big)\;.$$
  	
  	On the other hand, intersecting with $\underline{L}_\mu$ gives $\underline{\delta}.\underline{L}_\mu=1$. The curve $\underline{\delta}$ being isomorphic to $\mathbb{P}^1$, the latter implies, first  $\underline{\delta}.\underline{\delta}=\,\textrm{-}\,1$, secondly 
  	
  	$$2=\,\textrm{-}\, 4 + 2\sum_{j=1}^3 \alpha_j^2 + 2\beta^2 \quad\textrm{, hence}\quad \sum_{j=1}^3 \alpha_j^2 +\beta^2=3\;.$$
  	
  	We are left with two options. If $\beta\neq 0$, then $\beta=1$ and there exists a unique $l\neq 0$ such that $\alpha_l=0$, and $\alpha_j=1$ if $j\neq l$. In this case $D$ only goes through two out of the three points $\big\{g_\mu(\underline{s}_1),g_\mu(\underline{s}_2),g_\mu(\underline{s}_3)\}$, which is impossible because they are aligned. If on the contrary,  $\beta=0$, then $\alpha_j=1$ for any $j$, meaning $D$ goes through each $g_\mu(\underline{s}_j)$. In other words, $D= g_\mu(\underline{C}_0)$. Contradiction!\\
  	Let us consider at last the case $I_\mu(0) \geq 2$. There exists then at least one couple $j\neq k$ such that $\mu_j=0=\mu_k$, as well as an exceptional curve $\underline{\Gamma}_\nu$ fixed by $\iota_\mu$ (cf. §$3.4.5)$, hence contained in $\underline{R}_r=\underline{\Omega}+\underline{C}_0$, and in particular in $\underline{\Omega}$. $\blacksquare$

  \end{enumerate}
  We give hereafter more details about the irreducible components of $\underline{\Omega}\subset \underline{S}_\mu $ and the cubic $\mathcal{K}_0:=g_\mu(\underline{\Omega}) \subset \mathbb{P}^2(\mathbb{C})$, but also a geometric characterization of all exceptional curves of $\underline{S}_\mu $, in terms of the cardinal number $I_\mu(0).$
  
  \begin{theo}\hspace*{2mm}
  	
  	\begin{enumerate}
  		\item If $I_\mu(0)=0$, the curve $\underline{\Omega}$ is smooth, irreducible and has genus $1$. Its image $g_\mu(\underline{\Omega})$ is a smooth cubic of $j$-module
  		
  		$$j=2^8\frac{(\lambda^2\,\textrm{-}\,\lambda+1)^3}{\lambda^2(\lambda\,\textrm{-}\,1)^2}\; \quad\textrm{, with}\quad\lambda:= \Pi_0^3\frac{2\mu_i}{\big(\mu^{(1)}\,\textrm{-}\,2 \mu_i\big)} \in \mathbb{Q}^\ast\quad.$$
  		
  		\item If $I_\mu(0)=1$, $\underline{\Omega}$ is isomorphic to $\mathbb{P}^1$. Moreover, $g_\mu(\underline{\Omega})$ is an irreducible cubic with a node at $g_\mu(\underline{r}_{i_0})$, where $i_0$ is the index such that $\mu_{i_0}=0$.
  		
  		\item In the latter cases, the irreducible cubic $\mathcal{K}_0=g_\mu(\underline{\Omega})$ intersects the line $g_\mu(\underline{C}_0)$ at $\{g_\mu(\underline{s}_j)\}$, and for each $j\in \{1,2,3\}=\{j,k,l\}$, $g_\mu(\underline{s}_j)$ is an inflection point of $\mathcal{K}_0$ if and only if $\mu_k=\mu_l$.
  		
  		\item If $I_\mu(0)=2$, there exists $j\neq k$ such that $\mu_j= \mu_{k}=0\neq \mu_{l}$ $($where $\{j,k,l\}=\{1,2,3\})$. Moreover, $\underline{\Omega}$ has two irreducible components, one of which is the exceptional curve $\underline{\Gamma}_\nu$, with $\nu_i=1$ if $i \in \{j,k\}$, and $\nu_i=\mu_i$ if $i \in \{0,l\}$. Its projection $g_\mu(\underline{\Omega})$ is the union of $g_\mu(\underline{\Gamma}_\nu)$, the line through the $($aligned$!)$ points $g_\mu(\underline{r}_j),g_\mu(\underline{r}_{k})$ and $g_\mu(\underline{s}_{l})$, and an irreducible conic going through the four points $\big\{g_\mu(\underline{r}_j),g_\mu(\underline{r}_{k}),g_\mu(\underline{s}_j),g_\mu(\underline{s}_{k})\big\}.$
  		
  		\item If $I_\mu(0)=3$, i.e.: $\mu_j=0, \forall j\neq 0$, $\underline{\Omega}$ is the union of three exceptional curves $\underline{\Gamma}_{\mu+ \epsilon}$, with $\epsilon \in \{(0,0,1,1),(0,1,0,1),(0,1,1,0)\} $. Hence, its projection is the union of three lines going through the following three triplets of aligned points$:$\\
  		 $\big\{g_\mu(\underline{s}_{j}),g_\mu(\underline{r}_k),g_\mu(\underline{r}_{l})\big\}$, where $\{j,k,l\}=\{1,2,3\}$.
  		
  		\item The exceptional curves in $\underline{S}_\mu$, not fixed by $\iota_\mu$, correspond, via $g_\mu :\underline{S}_\mu \to \mathbb{P}^2(\mathbb{C})$, with the lines going through one of the three points $\{g_\mu(\underline{s}_j)\}$ and tangent to $($or going through a singular point of$)$ the cubic $\mathcal{K}_0$.
  	\end{enumerate}
  	
  \end{theo}
  \textbf{Proof.}
  \begin{enumerate}
  	\item Recall that $\underline{\Omega}$ is irreducible, numerically equivalent to $3\underline{C}_0+\sum_1^3 \underline{s}_j$, hence of arithmetic genus $1$. We will actually prove that its normalization has genus $1$, i.e.: $\underline{\Omega}$ is smooth; in fact, a double cover of $\mathbb{P}^1$, ramified over four points whose cross-ratio we also calculate.\\
  	For any $\epsilon \in \big\{(1,1,0,0),(0,0,1,1),(1,\textrm{-}1,0,0),(0,0,1,\textrm{-}1)\big\}$, we have $\mu + \epsilon \in \mathbb{T}_1$. We also know that the exceptional curves $\underline{\Gamma}_{\mu+ \epsilon}$ and $\underline{\Gamma}_{\mu\,\textrm{-}\,\epsilon}$, are exchanged by the involution $\iota_\mu$, intersect at a unique (hence smooth) point of $\underline{\Omega}$, and their sum $\underline{\Gamma}_{\mu+ \epsilon}+\underline{\Gamma}_{\mu\,\textrm{-}\,\epsilon}$ belongs to the pencil $|\underline{L}_\mu\,\textrm{-}\,\underline{s}_1|$. Any $\underline{D} \in |\underline{L}_\mu\,\textrm{-}\,\underline{s}_1|$, except the four latter divisors and $2\underline{C}_0+\underline{s}_2+\underline{s}_3$, is isomorphic to $\mathbb{P}^1$ and has $0$ self-intersection. Moreover, the projection $g_\mu:\underline{D}\to g_\mu(\underline{D})$ has degree $2$, hence two simple ramification points, corresponding to
  	$\underline{D}\cap \underline{\Omega}$.
  	The latter defines a degree-$2$ projection, of (the normalization of) $\underline{\Omega}$ onto $|\underline{L}_\mu\,\textrm{-}\,\underline{s}_1|\simeq \mathbb{P}^1$, ramified over the following four (smooth!) points:  $$\underline{\Gamma}_{\mu+ \epsilon}\cap \underline{\Omega}\quad\textrm{ with}\quad\epsilon \in \big\{(1,1,0,0),(0,0,1,1),(1,\textrm{-}1,0,0),(0,0,1,\textrm{-}1)\big\}.$$ Hence, the normalization of $\underline{\Omega}$, as well as  $\underline{\Omega}$, is a genus-$1$ smooth curve. \\
  	The above double cover $\underline{\Omega} \to \mathbb{P}^1$ has been defined by intersecting $\underline{\Omega}$ with each element of $|\underline{L}_\mu\,\textrm{-}\,\underline{s}_1|\simeq \mathbb{P}^1$. Analogously, given that $(\underline{L}_\mu\,\textrm{-}\,\underline{s}_1).\underline{s}_1 =2$, we define a degree-$2$ projection $\underline{s}_1 \to \mathbb{P}^1$, such that the points $\underline{\Gamma}_{\mu+ \epsilon} \cap \underline{s}_1$ and $\underline{\Gamma}_{\mu\,\textrm{-}\,\epsilon} \cap \underline{s}_1$ have same image, for any $\epsilon \in \big\{(1,1,0,0),(0,0,1,1),(1,\textrm{-}1,0,0),(0,0,1,\textrm{-}1)\big\}$. The corresponding involution of $\underline{s}_1$ is just the restriction of $\iota_\mu$ to $\underline{s}_1$. Recall the isomorphism $\underline{s}_1 \simeq \mathbb{C}\cup \{\infty\}$ (§$2.4)$), identifying $\underline{\Gamma}_{\mu \pm \epsilon} \cap \underline{s}_1$ with $n_{\mu \pm \epsilon}:=\frac{1}{2}\big((\mu \pm \epsilon)^{(2)}\,\textrm{-}\,1 \big)$, as well as the quotient of $\underline{s}_1$ by $\iota_\mu$ with $|\underline{L}_\mu\,\textrm{-}\,\underline{s}_1|\simeq \mathbb{P}^1$. It follows that the discriminant of the above double cover $\underline{\Omega} \to \mathbb{P}^1$ is equal to
  	$$\big\{(\mu_0+\mu_1)^2,(\mu_0 \,\textrm{-}\,\mu_1)^2 ,(\mu_2+\mu_3)^2,(\mu_2 \,\textrm{-}\,\mu_3)^2 \} \;,$$
  	
  	(and its cross-ratio is $\lambda:= \Pi_0^3\frac{2\mu_i}{\big(\mu^{(1)}\,\textrm{-}\,2 \mu_i\big)}$, as announced).
  	
  	\item Let $I_\mu(0)=1$ and denote $i$ the unique index such that $\mu_i=0$. According to §$3.7.2$)):
  	$$\underline{r}_i.\underline{r}_i=\,\textrm{-}\,2\quad\,\textrm{and}\quad\underline{\Omega}\equiv 3\underline{C}_0+\sum_1^3 \underline{s}_j\,\;\textrm{-}\,\;\underline{r}_i.$$
  	The involution $\iota_\mu$ acts non-trivially on $\underline{r}_i$ (cf. §$3.4.6)$), while fixing the two points in $\underline{\Omega}\cap \underline{r}_i$. It immediately follows that $g_\mu(\underline{\Omega})$ has a node at $g_\mu(\underline{r}_i)$.
  	
  	\item Let us assume $\mu_1=\mu_2\neq 0 $, and consider the exceptional curves $\underline{\Gamma}_{\mu+\epsilon}$ and $\underline{\Gamma}_{\mu\,\textrm{-}\,\epsilon}$, where $\epsilon := (0,1,\textrm{-}1,0)$. The latter are exchanged by the Geiser involution $\iota_\mu$, and $\underline{s}_3\cap \underline{\Omega}$ is their unique point of intersection, naturally identified with $n_{\mu+\epsilon}=n_{\mu\,\textrm{-}\,\epsilon}=\frac{1}{2}(\mu^{(2)}+1)$. Hence, $g_\mu(\underline{\Gamma}_{\mu+\epsilon})=g_\mu(\underline{\Gamma}_{\mu\,\textrm{-}\,\epsilon})$ only intersects $\mathcal{K}_0=g_\mu(\underline{\Omega})$ at $g_\mu(\underline{s}_3)$. In other words, the latter is an inflection point of the cubic $\mathcal{K}_0$.\\
  	Conversely, assume $g_\mu(\underline{s}_3)$ is an inflection point of $\mathcal{K}_0$ and let $\mathcal{H}_3$ denote the corresponding tangent line. Then, its inverse image $\underline{D}_3:=g_\mu^\ast(\mathcal{H}_3) \in |\underline{L}_\mu\,\textrm{-}\,\underline{s}_3|$ can only intersect $\underline{\Omega}$ at $\underline{s}_3 \cap \underline{\Omega}$, because $\mathcal{H}_3$ only intersects $g_\mu(\underline{\Omega})$ at $g_\mu(\underline{s}_3)$. On the other hand, any element of the pencil $|\underline{L}_\mu\,\textrm{-}\,\underline{s}_3|$ intersects $\underline{s}_3$ at two points exchanged by $\iota_\mu$ (cf. §$3.4.6$)). In particular, $\underline{\Omega}\cap \underline{s}_3$ being fixed by $\iota_\mu$, $\underline{D}_3$ intersects $\underline{s}_3$ and $\underline{\Omega}_3$ with multiplicity two (only) at $\underline{s}_3 \cap \underline{D}_3$. The latter being a smooth point of both curves, $\underline{D}_3$ must have a singular point there. Hence, is reducible and equal to $\underline{\Gamma}_\nu+ \iota_\mu(\underline{\Gamma}_\nu)$, for some $\nu \in \mathbb{T}_3$. It follows that $n_\nu:=\frac{1}{2}(\nu^{(2)}\,\textrm{-}\,1)=\frac{1}{2}(\mu^{(2)}+1)$ is fixed by ${\iota_\mu}_{|\underline{s}_3} $, which can only happen if $\nu=\mu \pm (0,1,\textrm{-}1,0)$ and $\mu_1=\mu_2$.
  	
  	\item Let us suppose $\mu_1=\mu_2=0\neq \mu_3$, and consider this time the exceptional curve $\underline{\Gamma}_{\mu+\epsilon}$, where $\epsilon := (0,1,1,0)$. According to §$3.4.5)$ it is fixed by $\iota_\mu$, hence an irreducible component of $\underline{\Omega}$ going through the points $\{ g_\mu(\underline{s}_3), g_\mu(\underline{r}_1), g_\mu(\underline{r}_2)\}$. Its complement, $\underline{\Omega}\,\textrm{-}\,\underline{\Gamma}_{\mu+\epsilon}$, intersects $\underline{s}_1$ and $\underline{s}_2$, as well as 
  	$\underline{r}_1$ and $\underline{r}_2$, because $\underline{\Omega}. \underline{r}_1=2=\underline{\Omega}. \underline{r}_2$. Therefore, the latter projects onto a conic going through the four points $\{ g_\mu(\underline{s}_1), g_\mu(\underline{s}_2), g_\mu(\underline{r}_1), g_\mu(\underline{r}_2)\}$. At last, in order to prove its irreducibility, it is enough to construct a line going through $g_\mu(\underline{s}_3)$ and tangent to the latter conic. Let us consider the exceptional curves $\underline{\Gamma}_{\mu+\delta}$ and $\underline{\Gamma}_{\mu\,\textrm{-}\,\delta}$, where $\delta := (1,0,0,1)$. They are disjoint with $\underline{\Gamma}_{\mu + \epsilon}$ and exchanged by $\iota_\mu$. We also know they intersect $\underline{s}_3$ at the points $n_{\mu+\delta}:=\frac{1}{2}(\mu^{(2)}+1) +\mu_0+\mu_3$ and $n_{\mu\,\textrm{-}\,\delta}:=\frac{1}{2}(\mu^{(2)}+1) \,\textrm{-}\,\mu_0\,\textrm{-}\,\mu_3$, respectively, and their intersection point must belong to $\underline{\Omega}\,\textrm{-}\,\underline{\Gamma}_{\mu+\epsilon}$, which is fixed by $\iota_\mu$. It follows that the line $g_\mu(\underline{\Gamma}_{\mu+\delta})=g_\mu(\underline{\Gamma}_{\mu\,\textrm{-}\,\delta})$ is (not contained, but) tangent to the conic $g_\mu(\underline{\Omega}\,\textrm{-}\,\underline{\Gamma}_{\mu+\epsilon})$.
  	
  	\item The last case follows directly from the above results. In fact, whenever $\mu_j=\mu_k=0$  for some $j\neq k$, there exists a line contained in the cubic $\mathcal{K}_0$ (cf. §$3.4.5$)). Hence, $\mathcal{K}_0$ is the union of the three lines going through the three points  $\{ g_\mu(\underline{r}_j)\}$.
  	
  	\item This caracterization follows from a direct calculation. Recall that an exceptional curve $\underline{\Gamma}_\nu$ satisfying $\iota_\mu(\underline{\Gamma}_\nu)\neq \underline{\Gamma}_\nu$ (and $\nu \in \mathbb{T}_j$), projects onto a line $g_\mu(\underline{\Gamma}_\nu) $ going through the point $g_\mu(\underline{s}_j)$ and: either $\iota_\mu(\underline{\Gamma}_\nu). \underline{\Gamma}_\nu=1=\underline{\Gamma}_\nu.\underline{\Omega}$, or $\iota_\mu(\underline{\Gamma}_\nu). \underline{\Gamma}_\nu=0=\underline{\Gamma}_\nu.\underline{\Omega}$. In the former case, let us denote $\underline{p}_\nu$ its unique point of intersection with $\underline{\Omega}$. We can easily check that the line $g_\mu(\underline{\Gamma}_\nu)$ is tangent to $g_\mu(\underline{\Omega}_j)$ at $g_\mu(\underline{p}_\nu)$. In the second case, there must exist an index $i\in \{0,\cdots,3\}$ such that ($\mu_i=0,\,\underline{r}_i^2=\,\textrm{-}2$ and) $\underline{r}_i$ is a positive dimensional fibre satisfying $\underline{r}_i.\underline{\Gamma}_\nu=1$. Hence it intersects $\underline{\Gamma}_\nu$. We can then check that the line $g_\mu(\underline{\Gamma}_\nu)$ goes through $g_\mu(\underline{r}_i)$, which is a node of $g_\mu(\underline{\Omega})$ (see the proof of §$3.8.2)$).  $\quad \quad \blacksquare$
  \end{enumerate}
  \begin{remark}\hspace*{2mm}
  	The generic element of the pencil $|\widetilde{\gamma}_X(\mu,1)|\simeq |\underline{L}_\mu|_{\underline{p}_X}$ is a smooth irreducible $\iota_\mu$-invariant curve of genus $1$. Therefore, it is a double cover of a line, ramified at $\underline{p}_X\in \underline{C}_0$ plus three other points lying on the ramification divisor $\underline{\Omega}$. Up to a finite number of points, $\underline{\Omega}$ is made of all such triplets. Hereafter we develop the analogous results for the linear system $\psi_\mu(|\widetilde{\gamma}_X(\mu,2)|)$ and corresponding Severi Variety $\mathcal{S}\mathcal{V}_X(\mu,2)$.
  	\end{remark}
  
   \begin{lem}\hspace*{2mm}
   	
	For any $X \in \mathfrak{X}$, $\psi_\mu(|\widetilde{\gamma}_X(\mu,2)|)$ has dimension $2$ and arithmetic genus $3$. Its generic element, say $\underline{\Gamma}$, has the following properties$:$
	\begin{enumerate}
		\item $\underline{\Gamma}$ and $\underline{\Omega}$ intersect transversally at six points$;$
		
		\item  $\underline{\Gamma}$ is a $\iota_\mu$-invariant irreducible curve, with a node at $\underline{p}_X\in \underline{C}_0$ and smooth everywhere else, i.e.$:$ of geometric genus $2;$
		
		\item $\underline{\Gamma}$ is a degree-$2$ cover of the conic $g_\mu(\underline{\Gamma})$, ramified at its intersection with $ \underline{\Omega}.$

		\end{enumerate}

\end{lem}
\textbf{Proof.}\\

1) - The generic element of the pencil $|\underline{L}|_{\underline{p}_X}$, say $\underline{\Gamma}_1$, is a smooth, $\iota_\mu$-invariant, irreducible curve of genus $1$, double cover of the projective line $g_\mu(\underline{\Gamma}_1 )$. Besides $\underline{p}_X$, it has three other ramification points, all lying at $\underline{\Omega}\setminus\cup_{j=1}^3\underline{s}_j$. Hence, $2\underline{C}_0+\sum_{j=1}^3\underline{s}_j+\underline{\Gamma}_1  \in\psi_\mu(|\widetilde{\gamma}_X(\mu,2)|)$ and intersects $\underline{\Omega}$ at six distinct points. It follows that the generic element $\underline{\Gamma}$ has the same property.

2) - The linear system $|\widetilde{\gamma}_X(\mu,2)|$ has dimension and arithmetic genus equal $2$, and its generic element, say $\widetilde{\Gamma}$, is a smooth irreducible curve. Recall also that $\widetilde{s}_0\setminus \widetilde{C}_0$ has been canonically identified with $\mathbb{C}$ (cf. Lemma §$2.4$). Accordingly, $\widetilde{\Gamma}\cap\widetilde{s}_0$ and $\widetilde{\Gamma}_\mu\cap \widetilde{s}_0$ correspond to the values $n:=\frac{1}{2}(\mu^{(2)}\,\textrm{-}\,1)$ and $n\,\textrm{-}\,4$ respectively. It follows that $\underline{\Gamma}:=\psi_\mu(\widetilde{\Gamma}_\mu)$ inherits a node at $\underline{p}_X\in \underline{C}_0$, with tangent directions $n$ and $n\,\textrm{-}\,4$.\\
\indent 
The latter results force $\iota_\mu(\underline{\Gamma})$ to have a node with same tangent cone as $\underline{\Gamma}$ at $\underline{p}_X$. Assuming $\iota_\mu(\underline{\Gamma}) \neq \underline{\Gamma}$ would imply they intersect with multiplicity $\geq 4$ at $\underline{p}_X$, as well as at the six points of $\underline{\Gamma}\cap \underline{\Omega}$. We should then have $8=\underline{\Gamma}.\underline{\Gamma}=\iota_\mu(\underline{\Gamma}) . \underline{\Gamma}\geq 4+6=10$. Contradiction! \\

\indent 3) - Once we know $\underline{\Gamma}$ is $\iota_\mu$-invariant of geometric genus $2$, we deduce via the projection and Riemann-Hurwitz formula, that $g_\mu(\underline{\Gamma})$ is a smooth conic and the natural double cover $\underline{\Gamma} \to g_\mu(\underline{\Gamma})$ is ramified at six smooth points. Moreover, the latter points must lie on $\underline{\Omega}$ and make all $\underline{\Gamma}\cap\underline{\Omega}$ because $\underline{\Gamma}.\underline{\Omega}=6$. $\blacksquare$\\

	Given any $ \underline{\Gamma} \in \psi_\mu(|\widetilde{\gamma}_X(\mu,2)|)$, let $p_g(\underline{\Gamma})$ denote its geometric genus and $I_{\underline{p} }(\underline{\Gamma},\underline{\Omega} )$ its intersection multiplicity with $ \underline{\Omega} $ at $\underline{p}$. The singularities of $\underline{\Gamma} $, other than its node at $\underline{p}_X\in \underline{C}_0$, lie along $\underline{\Omega}$ but outside $\big\{\underline{s}_j,j\neq 0\big\}\bigcup \big\{\underline{r}_i,\mu_i= 0\big\} $. They can be easily read off by looking at the decomposition $\sum_{\underline{p}}I_{\underline{p} }(\underline{\Gamma},\underline{\Omega} )=6$, as explained hereafter.
	
	\begin{prop}\hspace*{2mm}
		
	\begin{enumerate}
		
		\item The odd terms in $\sum_{\underline{p}}I_{\underline{p} }(\underline{\Gamma},\underline{\Omega} )$ correspond to the ramification divisor of $\underline{\Gamma} \to g_\mu(\underline{\Gamma});$
		
		\item $\underline{\Gamma}$ is a reduced divisor whenever $\sum_{\underline{p}}I_{\underline{p} }(\underline{\Gamma},\underline{\Omega} )=2+2+2$, $4+2$ or $6;$
			
		\item $\underline{\Gamma}$ is irreducible and $p_g(\underline{\Gamma})=2$ if and only if $\sum_{\underline{p}}I_{\underline{p} }(\underline{\Gamma},\underline{\Omega} )=1+1+1+1+1+1;$
		
		\item $\underline{\Gamma}$ is irreducible and $p_g(\underline{\Gamma})=1$ if and only if $\sum_{\underline{p}}I_{\underline{p} }(\underline{\Gamma},\underline{\Omega} )=2+1+1+1+1 $ or $3+1+1+1$, in which case $\underline{\Gamma} $ has a node or a cusp, respectively$;$
	
		\item  $ \underline{\Gamma}\in \mathcal{S}\mathcal{V}_X(\mu,2)$, if and only if $\sum_{\underline{p}}I_{\underline{p} }(\underline{\Gamma},\underline{\Omega})$ has two odd terms.\\
		 More precisely, it has$:$
		\begin{enumerate}
			\item two nodes if $\sum_{\underline{p}}I_{\underline{p} }(\underline{\Gamma} ,\underline{\Omega})=2+2+1+1;$
			\item a node and a cusp if $\sum_{\underline{p}}I_{\underline{p} }(\underline{\Gamma},\underline{\Omega} )=2+3+1;$
			
			\item two cusps if $\sum_{\underline{p}}I_{\underline{p} }(\underline{\Gamma},\underline{\Omega} )=3+3;$
			
			\item a tacnode if $\sum_{\underline{p}}I_{\underline{p} }(\underline{\Gamma} ,\underline{\Omega})=4+1+1;$
			
			\item a higher cusp if $\sum_{\underline{p}}I_{\underline{p} }(\underline{\Gamma},\underline{\Omega} )=5+1.$
			
		\end{enumerate}

		\end{enumerate}   
 
\end{prop}
 
 \begin{defi}\hspace*{2mm}
 	
 	Let  $\mathfrak{Dec}(6)=\{ \vec{d}=(d_k)\in \mathbb{N}^6, \Sigma_kd_k=6, d_i \geq d_{i+1}, \, \forall i<6\}$ denote the set of all decompositions of $6$ as a sum of non-negative integers and $\mathcal{S}\mathcal{V}(\mu, 2):=\bigcup_\mathfrak{X}\mathcal{S}\mathcal{V}_X(\mu,2)$. To any $\vec{d}\in \mathfrak{Dec}(6)$ we associate the subvariety $\underline{\Omega}^{\vec{d}}:=\{\sum_k d_k\underline{p}_k, \underline{p}_k\neq \underline{p}_j \,\textrm{if}\,k\neq j\} \subset \underline{\Omega}^{(6)}$. We deduce the partition $\underline{\Omega}^{(6)}=\cup_{\mathfrak{Dec}(6)}\underline{\Omega}^{\vec{d}}$, as well as $\mathcal{S}\mathcal{V}(\mu,2)=\cup_{\mathfrak{Dec}(6)}\mathcal{S}\mathcal{V}^{\vec{d}}$, where $\mathcal{S}\mathcal{V}^{\vec{d}}$ denotes the inverse image of 
 	$\underline{\Omega}^{\vec{d}}$ with respect to the natural morphism $\mathcal{S}\mathcal{V}(\mu,2) \to \underline{\Omega}^{(6)}$.
 \end{defi}

\begin{remark}\hspace*{2mm}
	
The partition $\mathcal{S}\mathcal{V}(\mu,2)=\cup_{\mathfrak{Dec}(6)}\mathcal{S}\mathcal{V}^{\vec{d}}$ runs over all vectors $\vec{d} \in \mathfrak{Dec}(6)$ with two odd terms. We prove in Appendix $\textbf{B}$ that $\mathcal{S}\mathcal{V}^{\vec{d}}$ is a finite set for any $\vec{d}\neq (2,2,1,1,0,0)$.\\
\indent  In other words, that for a generic $X \in \mathfrak{X}$, any $\underline{\Gamma} \in \mathcal{S}\mathcal{V}_X(\mu,2)$ has just two nodes.
\end{remark}

 \section{The complementary degree-2 projection $\widetilde{S}_X \to \mathbb{P}^2(\mathbb{C})$.}
 
 We will henceforth fix $\mu \in \mathbb{T}_0$ and consider, for any $X \in \mathfrak{X}$, the blowing-up of $b_n\in \widetilde{s}_0 $, with $n:=\frac{1}{2}(\mu^{(2)}\,\textrm{-}\,1)$, the unique base point of $|\widetilde{\gamma}_X(\mu,2)|$, say $u: \widehat{S}_X \to \widetilde{S}_X$ (see §$2.5.2)$). Let $\widehat{b} \subset \widehat{S}_X$ denote its exceptional fiber, $\widehat{K}$ its canonical divisor and  $\widehat{C}_0,\widehat{l},\widehat{s}_i,\widehat{r}_i, \widehat{S}_i$ the strict transforms of the corresponding divisors of $\widetilde{S}_X$. \\
 \indent We recall hereafter the main properties of the strict transform $\widehat{\gamma}_X  :=u^\ast(\widetilde{\gamma}_X(\mu,2))\,\textrm{-}\,\widehat{b}$ and its associated projection, say $t: \widehat{S}_X  \to |\widehat{\gamma}_X|^\vee \simeq \mathbb{P}^2(\mathbb{C})$. We also give the matrix of intersection product in its Picard group $Pic(\widehat{S}_X)$.

 \begin{prop}\hspace*{2mm} 
 	\begin{enumerate}
 		\item We have $u_\ast \circ\, u^\ast = 2Id$ and $Pic(\widehat{S}_X)=Pic(\widetilde{S}_X)\oplus \widehat{b}\mathbb{Z}$, the direct sum being orthogonal with respect to the intersection of divisors.
 		
 		\item We have $u^\ast(\widetilde{C}_0)=\widehat{C}_0$, $u^\ast(\widetilde{l})=\widehat{l}$, $u^\ast(\widetilde{s}_0)=\widehat{s}_0+ \widehat{b}$, as well as the following identities$:$
 		
 		$$ u^\ast(\widetilde{S}_i) =\widehat{S}_i,\quad u^\ast(\widetilde{r}_i)=\widehat{r}_i,\quad \textrm{and}\quad  u^\ast(\widetilde{s}_i)=\widehat{s}_i\quad \textrm{for any}\quad i=0,\cdots,3.$$
 		
 		\item The canonical divisor is numerically equivalent to $\widehat{K}\equiv \textrm{-}\, 2\widehat{C}_0\,\textrm{-}\sum_i \widehat{s}_i$ and satisfies, for any $i=0,\cdots,3$ and $j=1,2,3:$
 		
 		$$\widehat{K}.\widehat{S}_i=\widehat{K}.\widehat{b}=\widehat{K}.\widehat{K}=\textrm{-}\,1,\;\; \widehat{K}.\widehat{l}=\textrm{-}\,2,\;\;  \widehat{K}.\widehat{C}_0=\widehat{K}.\widehat{r}_i=\widehat{K}.\widehat{s}_j =0,\;\; \textrm{and}\;\; \widehat{K}.\widehat{s}_0=1\;.$$
 		
 		\item Any $\widetilde{\Gamma} \in |\widetilde{\gamma}_X(\mu,2)|$ goes through the base point $b_n$ and $u^\ast(\widetilde{\Gamma} )\textrm{-}\,\widehat{b} \in |\widehat{\gamma}_X |.$
 		\item 
 		For all $\widehat{D} \in |\widehat{\gamma}_X |$ we have $u_\ast(\widehat{D})\in |\widetilde{\gamma}_X(\mu,2)|$. The latter morphism $\widehat{D} \mapsto u_\ast(\widehat{D})$ identifies $|\widehat{\gamma}_X |$ with $|\widetilde{\gamma}_X(\mu,2)|$. We also have$:$
 		
 		$$\widehat{\gamma}_X .\widehat{\gamma}_X =2,\;\; \widehat{\gamma}_X .\widehat{b}_X =1,\;\;\textrm{and}\;\; \,\widehat{\gamma}_X .\widehat{s}_i =\widehat{\gamma}_X .\widehat{C}_0 =\widehat{K}_X .\widehat{\gamma}_X =0, \;\; \forall i=0,\cdots,3.$$
 		
 		\item The generic element of $ |\widehat{\gamma}_X |$ is a smooth irreducible curve of arithmetic genus $2$.
 	\end{enumerate}
 	As for the reducible divisors in $|\widehat{\gamma}_X |$, they correspond to the reducible ones in $|\widetilde{\gamma}_X(\mu,2)|$. The latter are well known and always include the pencil $2\widetilde{C}_0+\sum_i \widetilde{s}_i+|\widetilde{\gamma}_X(\mu,1)|$. There are no other ones unless $\mu_k \leq 1$ for some $k=0,\cdots,3$. Let indeed $\mu(k)\in \mathbb{T}_0$ denote the vector satisfying $\mu(k)_i\,\textrm{-}\,\mu_i=2\delta_{k,i}$ for any $i=0,\cdots,3$. Then, any other reducible divisor is, either equal to $\widetilde{\Gamma}_{\mu(k)}+\widetilde{r}_k$ or belongs to the pencil $2\widetilde{C}_0+\sum_i \widetilde{s}_i+|\widetilde{\gamma}_X(\mu(k),1)|$, for some index $k$ such that $\mu_k=1$ or $\mu_k=0$ respectively.
 \end{prop}
 
 Having identified $|\widehat{\gamma}_X |$ with $|\widetilde{\gamma}_X(\mu,2)|$, let us denote $m\in |\widehat{\gamma}_X |^\vee$ the dual of the pencil (corresponding to) $2\widetilde{C}_0+\sum_i\widetilde{s}_i+|\widetilde{\gamma}_X(\mu,1)|\subset |\widetilde{\gamma}_X(\mu,2)|$ and $\mathcal{H}\subset |\widehat{\gamma}_X |^\vee$ the dual of the point (corresponding to) $\widetilde{\Gamma}_\mu +4\widetilde{C}_0+2\sum_i \widetilde{s}_i$, so that $m\in \mathcal{H}$. Let as well $ \mathcal{H}_k \subset |\widehat{\gamma}_X |^\vee$ denote the line corresponding to the point $\widetilde{\Gamma}_{\mu(k)}+\widetilde{r}_k$ if $\mu_k=1$, and  $m_k \in \mathcal{H}_k $ the point corresponding to the pencil $2\widetilde{C}_0+\sum_i \widetilde{s}_i+|\widetilde{\gamma}_X(\mu(k),1)|$, whenever $\mu_k=0$.\\
 \indent The following properties of $t: \widehat{S}_X  \to |\widehat{\gamma}_X|^\vee \simeq \mathbb{P}^2(\mathbb{C})$ follow.
 
 \begin{prop} $($cf.\,\cite{T2}$)$\hspace*{2mm}
 	
 	\begin{enumerate}
 		
 		\item $t^{\textrm{-}1}(m)=\widehat{C}_0\cup \bigcup_i\widehat{s}_i$.
 		\item $t^{\textrm{-}1}(m_k)=\widehat{r}_k$, whenever $\mu_k=0$.
 		\item $t$ is a finite degree-$2$ cover over the complement of $\{m\}\cup \{m_k, \mu_k=0\} \subset |\widehat{\gamma}_X |^\vee$.
 		\item $t(\widehat{\Gamma}_\mu)= \mathcal{H} =t(\widehat{b})\quad$ and $\quad t^\ast(\mathcal{H})=4\widehat{C_0}+2\sum_i\widehat{s}_i+\widehat{\Gamma}_\mu+\widehat{b}$.
 	\end{enumerate}
 \end{prop}
 
 \begin{defi}\hspace*{2mm}
 	  Let $\widehat{\mathcal{R}}_X$ denote the ramification divisor of $t$, $\mathcal{D}_X:=t_\ast(\widehat{\mathcal{R}}_X)$ its discriminant and $\widehat{\rho}_X$ the reduced inverse image of $\mathcal{D}_X$.
 \end{defi}
\begin{remark}
	 The latter divisors will be henceforth denoted $\widehat{\mathcal{R}}, \mathcal{D}$ and $\widehat{\mathcal{\rho}}$ to simplify notations, although they are not independent of $X\in \mathfrak{X}$. The divisor $\widehat{\mathcal{R}}\,\textrm{-}\,\widehat{\rho}$ is effective, with support contained in the fibers of $t$. We dress hereafter the list of their basic properties.
 \end{remark}
 
  \begin{prop}\hspace*{2mm}
  	\begin{enumerate}
  		\item $\widehat{\mathcal{\rho}}$ is linearly equivalent to $\quad3( \widehat{\mathcal{C}}_0+\widehat{\Gamma}_\mu+\widehat{b}+ \widehat{s}_0) +\sum_{j\neq 0}\widehat{s}_j\;\textrm{-}\;\sum_{\mu_i =0}\widehat{r}_i\;;$
  		
  		\item $\widehat{\mathcal{R}}\,\textrm{-}\,\widehat{\rho}$ is equal to $\quad 7\widehat{\mathcal{C}}_0 +2 \widehat{s}_0 +4\sum_{j\neq 0}\widehat{s}_j+\sum_{\mu_i =0}\widehat{r}_i\;;$
  		
  		\item $\mathcal{D}$ is a sextic and has a singular point at $m$, with three smooth branches tangent to the line $\mathcal{H};$
  		
  		\item $\mathcal{D}$ has also a node at each point of $\{m_k, \mu_k=0\}$, is smooth everywhere else and decomposes as sum of $Max(1,I_\mu(0))$ irreducible components. More precisely$:$
  		
  		\begin{enumerate}
  			\item $\mathcal{D}$ is irreducible and has geometric genus $p_g(\mathcal{D})=1\,\textrm{-}\,I_\mu(0)$, whenever $I_\mu(0) \leq 1;$
  			\item $\mathcal{D}=\mathcal{C}\cup \mathcal{Q}$, with $\mathcal{C}$ and $\mathcal{Q}$ a smooth conic and a rational quartic, if $I_\mu(0) =2;$
  			\item $\mathcal{D}=\mathcal{C}_1\cup \mathcal{C}_2 \cup \mathcal{C}_3$, union of three smooth conics, if $I_\mu(0) =3;$
  		\end{enumerate}
  		\item The singular point $m\in \mathcal{D}$ has $\delta$-invariant $\delta_m=9$.
  	\end{enumerate}
  	
  	 \end{prop}
 \textbf{Proof.}\\

 	\begin{enumerate} 
 		
 		\item The generic $\widehat{\Gamma} \in |\widehat{\gamma}_X|$ is a smooth irreducible curve of genus $2$, double cover of the line $t(\widehat{\Gamma})$. Therefore, it intersects $\widehat{\rho}$ at its set of (six) Weierstrass points. According to §$3.9$, the latter are the inverse images of $\underline{\Gamma}\cap \underline{\Omega}$, where $\underline{\Gamma}:=\psi_\mu(u(\widehat{\Gamma} ))\subset \underline{S}_\mu$.
 		 Hence, $\widehat{\rho}$ can be characterized as the inverse image of $\underline{\Omega}$. The linear equivalence of $\underline{\Omega}$ being known, we deduce that of $\widehat{\rho}$.
 		 
 		 \item Recall that by definition, $\widehat{\mathcal{R}}\equiv  \widehat{K}+3\widehat{\gamma}_X$. It follows that $\widehat{\mathcal{R}}\,\textrm{-}\,\widehat{\rho}$ is linearly equivalent, hence equal to $\quad 7\widehat{\mathcal{C}}_0 +2 \widehat{s}_0 +4\sum_{j\neq 0}\widehat{s}_j+\sum_{\mu_i =0}\widehat{r}_i$.
 		 
 		 \item In particular, according to the projection formula, $\mathcal{D}.\mathcal{H}=\widehat{\rho}.t^\ast(\mathcal{H})=6$, hence $\mathcal{D}$ is a sextic. We also deduce that $\widehat{\rho}$ intersects each $\widehat{s_j}$ with multiplicity $1$ $(j=1,2,3)$, and each one of these branches intersects $t^\ast(\mathcal{H})$ with multiplicity $2$. The latter give three branches of $\mathcal{D}$ tangent to $\mathcal{H}$ at the singular point $m\in \mathcal{D}$. 
 		
 		\item 
 		It follows from Lemma §$3.4.6)$ that $\underline{\Omega}$ intersects each $(\textrm{-}2)$-curve $\{ \underline{r}_k, \mu_k=0\}$ transversally at two points. Analogously, $\widehat{\rho}$ intersects each $(\textrm{-}2)$-curve $\{ \widehat{r}_k, \mu_k=0\}$ transversally at two points, and its projection $\mathcal{D}$ has a node at such point $m_k$. Everywhere else, outside the fibers of $t$, $\mathcal{D}$ is isomorphic to $\widehat{\rho}$, hence smooth. The decomposition of $\mathcal{D}$ and the geometric geni of its irreducible components in terms of $I_\mu(0)$ follow from that of $\underline{\Omega}$ (see §$ 3.8$).
 		
 		\item Let us consider each case $I_\mu(0)=0,\cdots,3$ separatly. \begin{enumerate}
 			\item When $I_\mu(0)=0$ the sextic $\mathcal{D}$ is irreducible, has geometric genus $p_g=1$, arithmetic genus $p_a=\frac{5.4}{2}=10$ and is smooth outside its singular point $m$. Hence, $\delta_m=p_a\,\textrm{-}\,p_g =10\,\textrm{-}\,1=9$. 
 			
 			\item When $I_\mu(0)=1$, say $\mu_k=0$, $\mathcal{D}$ is a rational irreducible sextic, smooth outside $m$ and $m_k$ (where it has a node). Hence, $\delta_m=p_a\,\textrm{-}\,p_g \,\textrm{-}\,1=10\,\textrm{-}\,0\,\textrm{-}\,1=9$. 
 			
 			\item  When $I_\mu(0)=2$, say $\mu_k=0=\mu_j$, $\mathcal{D}=\mathcal{C}\cup\mathcal{Q}$ is the union of a smooth conic and a rational irreducible quartic smooth outside $m$. The latter intersect transversally at $m_k$ and $m_j$. Hence $\mathcal{D}$ has geometric genus $\,\textrm{-}\,1$ and $\delta_m=p_a\,\textrm{-}\,p_g\,\textrm{-}\,1\,\textrm{-}\,1= 10\,\textrm{-}\,(\textrm{-}\,1)\,\textrm{-}\,1\,\textrm{-}\,1=9$.
 			
 			\item At last, when $I_\mu(0)=3$, $\mathcal{D}$ is the union of three smooth conics, intersecting transversally each pair of them at one of the points $\{m_1, m_2,m_3\}$. Hence $\mathcal{D}$ has geometric genus $\,\textrm{-}\,2$, and three nodes. We deduce again that $\delta_m=9$.
 		\end{enumerate}
 	\end{enumerate} 
 
 \begin{defi}\hspace*{2mm}
 	
Given any $\widehat{\Gamma} \in |\widehat{\gamma}_X|$, we will denote $t(\widehat{\Gamma})^\vee \in \mathbb{P}^2(\mathbb{C})^\vee$ the dual of the line $t(\widehat{\Gamma}) \subset \mathbb{P}^2(\mathbb{C})$ and $\widehat{\Gamma}_r$ its reduced inverse image in $\widehat{S}_X$. We will also denote $\mathcal{H}^\vee \in m^\vee \subset \mathcal{D}^\vee $ the duals of the line $\mathcal{H}:=t(4\widehat{C_0}+2\sum_i\widehat{s}_i+\widehat{\Gamma}_\mu+\widehat{b})$, the point $m \in \mathcal{H}$ and the discriminant $\mathcal{D}$. 
\end{defi}

The following Lemma characterizes the divisors $\widehat{\Gamma} \in |\widehat{\gamma}_X|$ such that $t(\widehat{\Gamma})^\vee$ is a singular point (of pre-determined analytical type) of $\mathcal{D}^\vee$.

 \begin{lem}\hspace*{2mm}
 	
	For any $\widehat{\Gamma} \in |\widehat{\gamma}_X|$ the point $t(\widehat{\Gamma})^\vee$ does not belong to $\mathcal{D}^\vee$ if and only if $\widehat{\Gamma}$ is transverse to $\widehat{\rho}$ $($i.e.$:$ $\#\,\widehat{\Gamma}\cap \widehat{\rho}=6)$. Otherwise, $\widehat{\Gamma} _r$ has at least one singular point along $\widehat{\rho}$, and is related to $t(\widehat{\Gamma})^\vee$ as follows$:$
	\begin{enumerate}
		\item if $\widehat{\Gamma} _r$ is reducible then, either $(\widehat{\Gamma}=\widehat{\Gamma}_\mu+\widehat{b}+4\widehat{C_0}+2\sum_i\widehat{s}_i$ hence$)$ $\widehat{\Gamma}_r=\widehat{\Gamma}_\mu+\widehat{b}$, or there exists $k=0,\cdots,3$ such that $\mu_k=1$ and $\widehat{\Gamma}=\widehat{\Gamma}_r=\widehat{\Gamma}_{\mu(k)}+\widehat{r}_k$, where $\mu(k) \in \mathbb{T}_0$ satisfies $\mu(k)_i=\mu_i+2\delta_{i,k}$, for any $i=0,\cdots,3$. In both cases $t(\widehat{\Gamma})^\vee \in \mathcal{D}^\vee$ is a singular point of multiplicity $3$, and $\delta$-invariant equal to $9$ and $3$ respectively$;$
		
			\item if $\widehat{\Gamma} _r$ is irreducible but $\widehat{\Gamma}$ reducible then, $\mathcal{D}^\vee$ has a smooth point or a cusp at  $t(\widehat{\Gamma})^\vee;$
		
		\item if $\widehat{\Gamma} _r$ is irreducible and equal to $\widehat{\Gamma}$, it has geometric genus $0\leq p_g(\widehat{\Gamma})\leq 1$ and $:$
		
		\begin{enumerate}
			\item if $p_g(\widehat{\Gamma})=1$ the curve $\widehat{\Gamma}$ has a unique node or cusp, and $\mathcal{D}^\vee$ has a smooth point or a cusp at $t(\widehat{\Gamma})^\vee ;$ 
			
			\item if $p_g(\widehat{\Gamma})=0$, for a generic $X\in \mathfrak{X}$  the curve $\widehat{\Gamma}$ has two nodes, and $\mathcal{D}^\vee$ has a node at $t(\widehat{\Gamma})^\vee .$
		\end{enumerate}
		
	\end{enumerate}
	
\end{lem}
\textbf{Proof.}\\

Recall that any reducible curve $\widetilde{\Gamma} \in |\widetilde{\gamma}_X(\mu,2)|$ is, either equal to $\widetilde{\Gamma}_\mu+4\widetilde{C}_0+2\sum_i\widetilde{s}_i$, to $\widetilde{\Gamma}_{\mu(k)}+\widetilde{r}_k$ where $\mu_k=1$, or $\widetilde{C}+ 2\widetilde{C}_0+\sum_i\widetilde{s}_i +\widetilde{r}_k$ where $\mu_k=0$ and  the irreducible curve $\widetilde{C}$ belongs to the pencil $|\widetilde{\gamma}_X(\mu(k),1)|$. In the latter case, $\widetilde{C}$ has arithmetic genus $1$ and, for generic $X\in \mathfrak{X}$, is either smooth or rational with one node (cf. \cite{T2} §$7.8.4)$ \& $9.6.2))$.  In both cases $\mu(k) \in \mathbb{T}_0$ satisfies $\mu(k)_i=\mu_i+2\delta_{i,k}, \forall i=0,\cdots,3$.

\begin{enumerate}
	\item In case $\widehat{\Gamma}_r$ is reducible, it is either equal to $\widehat{\Gamma}_\mu+\widehat{b}$ or there exists an index $k$ such that $\mu_k=1$ and $\widehat{\Gamma}_r=\widehat{\Gamma}=\widehat{\Gamma}_{\mu(k)}+\widehat{r}_k$. In the first case  $t(\widehat{\Gamma})$ is the common tangent line to the three branches of $\mathcal{D}$ at its singular point $m$. It follows that $\mathcal{D}^\vee$ has a singular point at $\mathcal{H}^\vee=t(\widehat{\Gamma})^\vee $ of same analytical type, hence $\delta_{\mathcal{H}^\vee}=\delta_m=9$. In the second case, instead, the exceptional curve $\widetilde{\Gamma}_{\mu(k)}$ is transverse to $\widetilde{r}_k$, at least for a generic $X\in \mathfrak{X}$. It follows that $\widehat{\Gamma}_{\mu(k)}$ is transverse to $\widehat{r}_k$, hence (the line $t(\widehat{\Gamma}_{\mu(k)})$ is tritangent to $\mathcal{D}$ and) $\mathcal{D}^\vee$ has an ordinary singularity of multiplicity (and $\delta$-invariant) $3$ at $t(\widehat{\Gamma}_{\mu(k)})^\vee$.
	
	\item In case $\widehat{\Gamma}_r$ is irreducible but $\widehat{\Gamma}$ reducible, it is a rational curve of arithmetic genus $1$, with a node or a cusp, and $t(\widehat{\Gamma}^\vee)$ is either a smooth point or a cusp of $\mathcal{D}^\vee$.
	
	\item At last, in cases $\widehat{\Gamma}_r$ is irreducible and $\widehat{\Gamma}_r=\widehat{\Gamma}$, it has arithmetic genus $2$ and geometric genus $p_g(\widehat{\Gamma}) \leq 1$. If $p_g(\widehat{\Gamma}) = 1$ it must have a node or cusp, hence $t(\widehat{\Gamma}^\vee)$ is either a smooth point or a cusp of $\mathcal{D}^\vee$. Otherwise $p_g(\widehat{\Gamma}) = 0$, i.e.: $\widehat{\Gamma}$ is a rational irreducible curve, i.e.: $\widetilde{\Gamma} \in \mathcal{S}\mathcal{V}_X(\mu,2)$, and for generic $X\in \mathfrak{X}$ it has two nodes (cf. \textbf{B}.$11$). Therefore, $t(\widehat{\Gamma})$ is bitangent to $\mathcal{D}$, and $\mathcal{D}^\vee$ has a node at $t(\widehat{\Gamma})^\vee$. $\quad \quad \blacksquare$
	\end{enumerate}
We summarize hereafter the basic properties of $\mathcal{D}^\vee$ and its singularities, in terms of $I_\mu(0)$ and $I_\mu(1)$, the number of coordinates of $\mu$ equal to $0$ or $1$, respectively.

 \begin{theo}\hspace*{2mm}
 	
 	The curve $\mathcal{D}^\vee$ has
 	$Max\big(1,I_\mu(0)\big)$ irreducible components, degree $d^\vee =12\,\textrm{-}\,I_\mu(0)$ and a singularity at $\mathcal{H}^\vee$, with three smooth branches tangent to a common line and $\delta$-invariant $\delta_{\mathcal{H}^\vee}=9$.\\
 	Moreover, for generic $X\in \mathfrak{X}$, the remaining singular points are$:$
  \begin{enumerate}
  	\item  $I_\mu(1)$ ordinary singularities of multiplicity $3;$
  	
  	\item $18\,\textrm{-}\,6I_\mu(0)$ cusps and $27\,\textrm{-}\,14I_\mu(0)+2I_\mu(0)^2\,\textrm{-}\,3I_\mu(1)$ nodes.
  \end{enumerate}
 	\end{theo}
 \textbf{Proof.}\\
 
Consider a plane curve $\mathcal{C}$ of degree $d\geq 2$ and geometric genus $g$. For any $p\in \mathcal{C}$ we let $m_p$ and $\nu_p$ denote the multiplicity and number of branches of $\mathcal{C}$ at $p$, and $d^\vee$ the degree of its dual curve $\mathcal{C}^\vee$. According to Pl\"{u}cker's formulas$:$
 $$2\,\textrm{-}\,2g=2d\,\textrm{-}\,d^\vee \,\textrm{-}\,\sum_{p\in \mathcal{C}} (m_p\,\textrm{-}\,\nu_p) .$$
 
 For $\mathcal{C}=\mathcal{D}$ we have $d=6$, $g=1\,\textrm{-}\,I_\mu(0)$, and $m_p=\nu_p$ for any $p\in \mathcal{D}$. Hence $\mathcal{D}^\vee$ has degree $d^\vee=12\,\textrm{-}\,2I_\mu(0)$, as well as a singular point at $\mathcal{H}^\vee$ of $\delta$-invariant $\delta_{H^\vee}=9$ (as proven in the above lemma). \\
 	Since $\mathcal{D}$ is the dual of its dual curve, and has the same geometric genus, it follows from Pl\"{u}cker's formulas that $2I_\mu(0)=2d^\vee\,\textrm{-}\,d\,\textrm{-}\,\sum_{p^\vee\in \mathcal{D}^\vee} (m_{p^\vee}\,\textrm{-}\,\nu_{p^\vee}).$ In other words $$\sum_{p^\vee\in \mathcal{D}^\vee} (m_{p^\vee}\,\textrm{-}\,\nu_{p^\vee})=2d^\vee\,\textrm{-}\,d\,\textrm{-}\,2I_\mu(0)=18\,\textrm{-}\,6I_\mu(0).$$ 
 	
 	Recall now that for generic $X\in \mathfrak{X}$ and any $p^\vee\in \{t(\widehat{r}_k ), \mu_k=1 \} \subset \mathcal{D}^\vee$, we have $(m_{p^\vee}\,\textrm{-}\,\nu_{p^\vee}, \delta_{p^\vee})=(0,3)$. Moreover, the remaining singularities are nodes and cusps, where $(m_{p^\vee}\,\textrm{-}\,\nu_{p^\vee}, \delta_{p^\vee})=(0,1)$ or $(1,1)$ respectively. It follows that 
 	$$\sum_{p^\vee\in \mathcal{D}^\vee} (m_{p^\vee}\,\textrm{-}\,\nu_{p^\vee})=18\,\textrm{-}\,6I_\mu(0)$$ counts the number of cusps.\\
 	
 	As for the number of nodes, equal to $\sum_{p^\vee} \delta_{p^\vee}\,\textrm{-}\,3I_\mu(1)\,\textrm{-}\,\delta_{\mathcal{H}^\vee}\,\textrm{-}\,18+6I_\mu(0)$, we deduce it from $\delta_{\mathcal{H}^\vee}=9$ and 
 	
 	$$p_a(\mathcal{D}^\vee)\,\textrm{-}\,p_g(\mathcal{D}^\vee)=\sum_{p^\vee \in \mathcal{D}^\vee} \delta_{p^\vee},$$
 	
 	where $p_a(\mathcal{D}^\vee)=\frac{1}{2}(d^\vee\,\textrm{-}\,1)(d^\vee\,\textrm{-}\,2)$ denotes its arithmetic genus and $p_g(\mathcal{D}^\vee)=1\,\textrm{-}\,I_\mu(0)$ its geometric genus.        $\quad \quad \quad\blacksquare$
 \begin{cor}\hspace*{2mm}
 	
 	For any $\mu\in \mathbb{T}_0$ and generic $X\in \mathfrak{X}$, the Severi set $\mathcal{S}\mathcal{V}_X(\mu,2)$ is in a $1\,\textrm{-}\,1$ correspondence with the subset of nodes of $\mathcal{D}$ and has, in particular, cardinal $$\#\mathcal{S}\mathcal{V}_X(\mu,2)=27\,\textrm{-}\,14I_\mu(0)+2I_\mu(0)^2\,\textrm{-}\,3I_\mu(1).$$
 	
 \end{cor}
	 
	 	  \section{On spectral data of even hyper-elliptic potentials.}
	 	  
	 	  Given any $X \in \mathfrak{X}$, identified with an elliptic curve, we let $\wp$ denote its Weierstrass function, with a double pole at its origin $\omega_0$, $\wp'$ its derivative, vanishing at the three other half-periods $\{\omega_1, \omega_2, \omega_3\} $ and such that $\wp'(x)^2=4\wp^3(x)\,\textrm{-}\,g_2\wp(x)\,\textrm{-}\,g_3= 4\Pi_1^3(\wp(x)\,\textrm{-}\,e_j)$, with $e_j=\wp(\omega_j)$. For any $\alpha \in \mathbb{N}^4$ we then denote $\mathcal{P}ot_X(\alpha,2)$ the finite set of even hyper-elliptic potentials decomposing as 
	 	  
	 	  $$\sum_0^3 \alpha_i(\alpha_i+1)\wp(x\,\textrm{-}\,\omega_i)+2\big[\wp(x\,\textrm{-}\,\rho_1)+\wp(x+\rho_1)+\wp(x\,\textrm{-}\,\rho_2)+\wp(x\,+\rho_2)\big]$$
	 	  for some $\rho_1\neq \rho_2$ (not in $\{\omega_i\})$. Any such even doubly periodic function is known to be a finite-gap potential, if and only if the pair $(\rho_1,\rho_2)$ satisfies the so-called Duistermaat-Gr\"{u}nbaum set of equations\\
	 	  
	 	  $(D\,\textrm{-}\,G) \left\{ \begin{array}{rcl}
	 	  	8[\wp'(\rho_1\,\textrm{-}\,\rho_2)+\wp'(\rho_1+\rho_2)]+ \sum_0^3 (2\alpha_i+1)^2\wp'(\rho_1\,\textrm{-}\,\omega_i)=0 \\\hspace*{2mm} \\
	 	  	8[\wp'(\rho_2\,\textrm{-}\,\rho_1)+\wp'(\rho_2+\rho_1)]+ \sum_0^3 (2\alpha_i+1)^2\wp'(\rho_2\,\textrm{-}\,\omega_i)=0
	 	  \end{array}\right.$
 	  
 	  \begin{lem}\hspace*{2mm}
 	  	
 	  	Let us denote $\Pi(x)=\Pi_1^3(x\,\textrm{-}\,e_j)$, $P_j(x)=(x\,\textrm{-}\,e_k)^2(x\,\textrm{-}\,e_l)^2$ and $E_j=(e_j\,\textrm{-}\,e_k)(e_j\,\textrm{-}\,e_l)$ $($with $\{j,k,l\}= \{1,2,3\})$. The above (D-G) system is equivalent to the vanishing of two polynomials in the variables $(x,y)= \left(\wp(\rho_1), \wp(\rho_2)\right)$, namely, $F(x,y)$ and $F(y,x)$, where 
 	  	
 	  	$$F(x,y)=4\Pi(x)^2[(12y^2\,\textrm{-}\,g_2)(x\,\textrm{-}\,y)+16\Pi(y)]+(x\,\textrm{-}\,y)^3[(2\alpha_0+1)^2\Pi(x)^2 \,\textrm{-}\,\sum_1^3(2\alpha_j+1)^2E_jP_j(x)]$$

 	  	\end{lem}
   	
 	  	\textbf{Proof.}\\
 	  	 The following (well known) addition formulas for $\wp$ and $\wp'$ can be checked directly:\\
 	  	 
 	  	 $\quad \quad \wp'(z\,\textrm{-}\,\omega_j)[\wp(z)\,\textrm{-}\,\wp(e_j)]^2=\,\textrm{-}\,\wp'(z)E_j \quad \quad \textrm{and}$ $$2[\wp'(z\,\textrm{-}\,w)+\wp'(z+w)][\wp(z)\,\textrm{-}\,\wp(w)]^3=\,\textrm{-}\,\wp'(z)[(12\wp(w)^2\,\textrm{-}\,g_2)(\wp(z)\,\textrm{-}\,\wp(w))+16\Pi(\wp(w))].$$
 	  	
 	  	Replacing in (D-G) with $(z,w)=(\rho_1,\rho_2)$ and simplifying by $\wp'(\rho_1)$ and $\wp'(\rho_2)$, which do not vanish because $\rho_1$ and $\rho_2$ are not half-periods, we immediately obtain the equivalent system $F(x,y)=0=F(y,x)$.  $\blacksquare$
 	  	
 	  	 \begin{lem}\hspace*{2mm}
 	  	 	
 	  	 	The above affine plane curves $\{F(x,y)=0\}$ and $\{F(y,x)=0\}$ intersect the diagonal $\Delta \subset \mathbb{C}^2$ at the three points $\{(e_j,e_j)\}$, and have singularities of multiplicity three at each one of them.

 	  	 	\end{lem}
 	  	
 	  	\textbf{Proof.}\\
 	  	
 	  	Since $F(x,x)=4.16\Pi(x)^2\Pi(x)$, we deduce that $\Delta \cap \{F(x,y)=0\}=\{(e_j,e_j)\}$. On the other hand, at any point $(e_j,e_j)$, with local coordinates $x_j:=x\,\textrm{-}\,e_j$ and $y_j:=y\,\textrm{-}\,e_j$, we have the following development for $F(x,y):$ $$x_j^2[\,\textrm{-}\,4E_j^2+ O(x_j)][(12e_j^2\,\textrm{-}\,g_2)(x_j\,\textrm{-}\,y_j)+16E_jy_j+O(y_j^2)]+(x_j\,\textrm{-}\,y_j)^3[(2\alpha_j+1)^2E_j^3 +O(x_j)].$$
 	  	
 	  	Hence, its zero-locus has a singularity of multiplicity $3$ at $(e_j,e_j)$, with tangent cone 
 	  	
 	  	$$\left\{\,4E_j^2x_j^2[(12e_j^2\,\textrm{-}\,g_2)(x_j\,\textrm{-}\,y_j)+16E_jy_j]\,\,\textrm{-}\,(2\alpha_j+1)^2E_j^3(x_j\,\textrm{-}\,y_j)^3=0\right\}.\quad \blacksquare$$
 	  	
 	  	\begin{prop}\hspace*{2mm}
 	  		
 	  		For any elliptic curve $X\in \mathfrak{X}$ and $\alpha \in \mathbb{N}^4$, the cardinal of $\mathcal{P}ot_X(\alpha,2)$ is bounded by $27$.
 	  		\end{prop}
 	  	
 	  		\textbf{Proof.}\\
 	  		
 	  		According to the above lemma, the degree-$9$ curves $\{F(x,y)=0\}$ and $\{F(y,x)=0)\}$ intersect at the three points $\{(e_j,e_j)\}$ with total multiplicity $3.9=27$, and nowhere else along the diagonal $\Delta$. Hence, outside $\Delta$ they intersect at no more than $9.9\,\textrm{-}\,27=54$ points, i.e.: no more than $27$ pairs of non-ordered pairs of points. In other words, $\#\mathcal{P}ot_X(\alpha,2) \leq 27$ as asserted. $\quad \blacksquare$
 	  		
 	  		\begin{defi}\hspace*{2mm}
 	  			
 	  			For any $(j,k)\in \mathbb{Z}_2\times \mathbb{Z}_4$, we define the map $\mathcal{C}^{j,k}: \mathbb{T}_0 \to \mathbb{N}^4,\quad \mu \mapsto \alpha:=\mathcal{C}^{j,k}(\mu):$
 	  			\begin{enumerate}
 	  				\item for $(j,k)=(0,0)$, $\alpha:=\mathcal{C}^{0,0}(\mu)$ has coefficients$:$
 	\begin{enumerate}
 	  \item $\alpha_0:=\frac{1}{2}(\mu^{(1)}\,\textrm{-}\,1)$ and
 	\item $\alpha_j:=\frac{1}{2}(|2(g_\mu\,\textrm{-}\,\mu_0\,\textrm{-}\,\mu_j)+1|\,\textrm{-}\,1)$, hence $\alpha_j(\alpha_j+1)=(g_\mu\,\textrm{-}\,\mu_0\,\textrm{-}\,\mu_j)(g_\mu\,\textrm{-}\,\mu_0\,\textrm{-}\,\mu_j+1);$
 	
 	  			\end{enumerate}
   			\item for $(j,k)=(1,0)$, the i-th coefficient of  $\alpha:=\mathcal{C}^{1,0}(\mu)$ is $\alpha_i:=\frac{1}{2}(|2(g_\mu\,\textrm{-}\,\mu_i)+1|\,\textrm{-}\,1)$.
 	  				
   			\item For any other $(j,k) \in \mathbb{Z}_2\times \mathbb{Z}_4$,  $\mathcal{C}^{j,k}(\mu)_i=\mathcal{C}^{j,0}(\mu)_{i+k}.$
   		\end{enumerate}
   	Let us also define the maps $M, S, m: \mathbb{N}^4 \to \mathbb{N}$ as follows$:$\\
   	
   	for any $\alpha \in \mathbb{N}^4$, $M(\alpha):=Max\{\alpha_i\}$, $m(\alpha):=min\{\alpha_i\}$ and $S(\alpha):=\sum_i \alpha_i$,\\
   	
   	and denote
   	 $ \left\{ \begin{array}{rcl}
   	\mathbb{N}^{4,\geq}:=\left\{\alpha \in \mathbb{N}^4, 2M(\alpha)\geq S(\alpha)\,\textrm{-}\,(1+(\textrm{-}\,1)^{S(\alpha)})m(\alpha)\right\} \\  \\
   	\mathbb{N}^{4,\leq}:=\left\{\alpha \in \mathbb{N}^4, 2M(\alpha)\leq S(\alpha)\,\textrm{-}\,(1+(\textrm{-}\,1)^{S(\alpha)})m(\alpha)\right\}
   	\end{array}\right.$

 	  			\end{defi}

 	  	\begin{prop}\hspace*{2mm} $($\cite{TV2} Prop. §$4.2$ \& Th. §$5.2)$\\
 	  		
 	  	 \noindent The map $\mathcal{C}^{j,k}$ is injective for any $(j,k) \in \mathbb{Z}_2\times \mathbb{Z}_4$. The images of $ \mathcal{C}^{0,k}$ and $\mathcal{C}^{1,k}$ are, respectively, $ \{\alpha \in \mathbb{N}^{4,\geq},  M(\alpha)=\alpha_k\}$ and
 	  	$$\{\alpha \in \mathbb{N}^{4,\leq},S\equiv 0 \;\textrm{mod}.2 \;\textrm{and}\; \alpha_k=m(\alpha)\} \cup \{\alpha \in \mathbb{N}^{4,\leq},S\equiv \alpha_k \;\textrm{-}\;\alpha_j\equiv 1\;\textrm{mod}.2 \;\forall j\neq k\}.
 	  	$$
 	  	
 	  	In particular, for any $\alpha \in \mathbb{N}^4$ there exists $(j,k)\in \mathbb{Z}_2\times \mathbb{Z}_4$ and a unique $\mu \in \mathbb{T}_0$ satisfying $\mathcal{C}^{j,k}(\mu)=\alpha$.
 	  	\end{prop}
   	
   	\begin{remark}\hspace*{2mm}
   		
   		We skip the formulas for each inverse map $\left(\mathcal{C}^{j,k}\right)^{\textrm{-1}}:\mathcal{C}^{j,k}\left( \mathbb{T}_0\right) \to \mathbb{T}_0$ $($to be found in \cite{TV2}, §$4.4\, \&\, $§$5.3$ with $\mathcal{C}^{0,k}=\mathcal{A}^{k}$ and $\mathcal{C}^{1,k}=\mathcal{B}^{k})$. They imply the following $($in-$)$equalities.
   	\end{remark}
   
   \begin{lem}\hspace*{2mm}
   	
   	For any $\big(\mu,(j,k)\big) \in \mathbb{T}_0\times\mathbb{Z}_2\times \mathbb{Z}_4$ and $\alpha:=\mathcal{C}^{j,k}(\mu)$, we have
   	\begin{enumerate}
   		\item $\mu^{(2)}\,\textrm{-}\,1=\sum_i\alpha_i(\alpha_i+1);$
   		\item $\mu^{(1)}\,\textrm{-}\,1=Max\left( 2M(\alpha), S(\alpha)+1\,\textrm{-}\,(1+(\textrm{-}\,1)^{S(\alpha)})(m(\alpha)+\frac{1}{2})\right);$
   		\item $min\{\mu_i\} \geq 2$ if and only if $ \left| 2M(\alpha) \,\textrm{-}\,S(\alpha)+(1+(\textrm{-}\,1)^{S(\alpha)})m(\alpha)\right| \geq 4$.
   	\end{enumerate} 
   	
   \end{lem}

\begin{remark}\hspace*{2mm}
	
	Recall that any spectral data $(\pi,\xi)$, with $\pi$ a degree-$n$ $\it{ht}$ cover and $\xi$ a theta-characteristic, gives rise to an even hyper-elliptic potential in $\mathcal{P}ot_X(\alpha,m)$, for some $(\alpha,m)\in \mathbb{N}^4 \times \mathbb{N}$ such that $2n =\sum_i\alpha_i(\alpha_i+1) +4m$. Conversely, any element in $\mathcal{P}ot_X(\alpha,m)$ corresponds to such spectral data. However, while all spectral curves are $\it{ht}$ covers of same degree $n= \frac{1}{2}\left(\sum_i\alpha_i(\alpha_i+1) +4m\right)$, they may have different types, and in particular, different arithmetic geni $($as it is already the case for $m=1)$.\\
	
	We will define hereafter eight theta-characteristics naturally associated to any $\it{ht}$ cover.
	\end{remark}
	
	\begin{defi}\hspace*{2mm}
	
	Let $\pi: (\Gamma,p) \to (X,\omega_0)$ be a degree-$n$ $\it{ht}$ cover of arithmetic genus $g$. Then, for any $k\in \mathbb{Z}^4$ we define the following theta-characteristic$:$	
	$$\xi_{0,k}:=\mathcal{O}_\Gamma \left((g\,\textrm{-}\,1\,\textrm{-}\,2n)p+\pi^*(\omega_k+\omega_0)\right),$$
	
	in particular $\xi_{0,0}:=\mathcal{O}_\Gamma \left((g\,\textrm{-}1)p\right)$, and 
	
	$$\xi_{1,k}:=\mathcal{O}_\Gamma \left((g\,\textrm{-}\,1\,\textrm{-}\,n)p+\pi^*(\omega_k)\right).$$
	\end{defi}
	
	\begin{prop} (\cite{TV2}, Th.$3.4$) \hspace*{2mm}
	
	For any $\it{ht}$ cover $\pi$ coming from $\mathcal{S}\mathcal{V}_X(\mu,2)$ and $\alpha= \mathcal{C}^{j,k}(\mu)$, for any $(j,k) \in \mathbb{Z}_2\times \mathbb{Z}_4$, the spectral data $(\pi,\xi_{j,k})$ gives rise to a hyper-elliptic potential in $\mathcal{P}ot_X(\alpha,2)$.
	
	\end{prop}
	
	\begin{remark}  \hspace*{2mm}
		
		According to §$5.3$ and §$5.8$, $\#\mathcal{P}ot_X(\alpha,2)$ ranges between $\#\mathcal{S}\mathcal{V}_X(\mu,2)$ and $27$. On the other hand, $\#\mathcal{S}\mathcal{V}_X(\mu,2)=27$ for generic $X\in \mathfrak{X}$, whenever $min\{\mu_i\}\geq 2$ $($cf. §$4)$. It follows that all spectral data corresponding to $\mathcal{P}ot_X(\alpha,2)$ have same arithmetic genus $g_\alpha$ and theta-characteristic $\xi_{j,k}$. In the remaining case $\#\mathcal{S}\mathcal{V}_X(\mu,2)<27$, we construct herafter $27 \,\textrm{-}\,\#\mathcal{S}\mathcal{V}_X(\mu,2)$ more spectral data giving rise to elements in $\mathcal{P}ot_X(\alpha,2)$, as follows directly from the same formula in \cite{TV2}, Th.$5.3$.
		
		\end{remark}
	
	\begin{lem}\hspace*{2mm}
		
		For any $(\mu,i_0) \in \mathbb{T}_0\times\mathbb{Z}_4$ satisfying $\mu_{i_0}=0$, let $\nu \in \mathbb{T}_0$ be such that $\nu_i=\mu_i +2\delta_{i,i_0}$ for any index $i$, hence $\nu_{i_0}=2$ and $\nu^{(2)}=\mu^{(2)}+4$. Let $\pi : (\Gamma,p) \to (X,\omega_0)$ be any $\it{ht}$ cover of type $\nu$ and arithmetic genus $\frac{1}{2}(\nu^{(1)}\,\textrm{-}\,1)=\frac{1}{2}(\mu^{(1)}\,\textrm{-}\,1)+1$ coming from $\mathcal{S}\mathcal{V}_X(\nu,1)$, and $\{p_1,p_2\}$ the unique pair of Weierstrass points in $\pi^*(\omega_{i_0}) \setminus \{p\}$. Then, the even hyper-elliptic potentials corresponding to $\xi_{j,k}(p_1\,\textrm{-}\,p)$ and $\xi_{j,k}(p_2\,\textrm{-}\,p)$ belong to $\mathcal{P}ot_X(\alpha,2)$.

		\end{lem}
	
	\begin{lem}\hspace*{2mm}
		
		For any $(\mu,i_0) \in \mathbb{T}_0\times\mathbb{Z}_4$ satisfying $\mu_{i_0}=1$, let $\nu \in \mathbb{T}_0$ be such that $\nu_i=\mu_i +2\delta_{i,i_0}$ for any index $i$, hence $\nu_{i_0}=3$ and $\nu^{(2)}=\mu^{(2)}+8$. Let $\pi : (\Gamma,p) \to (X,\omega_0)$ be the exceptional $\it{ht}$ cover of type $\nu$ and arithmetic genus $\frac{1}{2}(\nu^{(1)}\,\textrm{-}\,1)= \frac{1}{2}(\mu^{(1)}\,\textrm{-}\,1)+1$, and $\{p_1,p_2, p_3\}$ the unique Weierstrass points in $\pi^*(\omega_{i_0}) \setminus \{p\}$. Then, the three even hyper-elliptic potentials corresponding to $\xi_{j,k}(p_l\,\textrm{-}\,p)$ $(l=1,2,3)$ belong to $\mathcal{P}ot_X(\alpha,2)$.
	
\end{lem}

	\begin{lem}\hspace*{2mm}
		
	For any $\mu \in \mathbb{T}_0$ satisfying $\mu_{j'}=\mu_{k'}=0$ $(\{j',k',l'\}=\{1,2,3\})$, consider $\tau \in \mathbb{T}_3$, where $\tau_0=\mu_0, \tau_{l'}=\mu_{l'}$ and $\tau_{j'}=\tau_{k'}=2$. The exceptional $\it{ht}$ cover of type $\tau$ has degree $\frac{1}{2}(\tau^{(2)}\,\textrm{-}\,1)=\frac{1}{2}(\mu^{(2)}\,\textrm{-}\,1)+4=n$ and arithmetic genus $\frac{1}{2}(\tau^{(1)}\,\textrm{-}\,1)=\frac{1}{2}(\mu^{(1)}\,\textrm{-}\,1)+2$. Moreover, let $\{p_{j',1},p_{j',2}\}$ and $\{p_{k',1},p_{k',2}\}$ denote the two pairs of Weierstrass points over $\{\omega_{j'},\omega_{k'}\}$. Then, the four even hyper-elliptic potentials corresponding to $\xi_{j,k}(q_{j'}+ q_{k'}\,\textrm{-}\,2p)\}$ with $ q_{j'}\in \{p_{j',1},p_{j',2}\}$ and $ q_{k'}\in \{p_{k',1},p_{k',2}\}$ belong to $\mathcal{P}ot_X(\alpha,2)$.
	
	\end{lem}

Given $\alpha \in \mathbb{N}^4$, let $\big(\mu,j,k\big) \in \mathbb{T}_0\times\mathbb{Z}_2\times\mathbb{Z}_4$ be such that $\mathcal{C}^{j,k}(\mu)= \alpha$, and denote $g_\alpha:=\frac{1}{2}(\mu^{(1)}\,\textrm{-}\,1)$ and $(M,S,m)= \left(M(\alpha),S(\alpha),m(\alpha)\right)$ (see §$5.4$ and §$5.5$). Recall that $g_\alpha= \frac{1}{2}Max\left(2M,S+1\,\textrm{-}\,(1+(\,\textrm{-}\,1)^S)(m+\frac{1}{2})\right)$, and that for any $d\in \mathbb{N}$, the spectral curve of any potential in $\mathcal{P}ot_X(\alpha,d)$ is a degree-$n$ $\it{ht}$ cover, where $2n=\sum_i \alpha_i(\alpha_i+1)+4d$. \\
\indent The following theorem, together with the preceding lemmas, imply the next corollary. They generalize to the case $d=2$, the results already proved for $d=0,1$.

\begin{theo}\hspace*{2mm}
	
	 For generic $X\in \mathfrak{X}$ and any $\alpha \in \mathbb{N}^4$, the cardinal of $\mathcal{P}ot_X(\alpha,2)$ is equal to $27$.
 
\end{theo}
 \textbf{Proof.}\\
 We divide the proof according to the values of $\left(I_0(\mu),I_1(\mu)\right)$. Recall that for any $\mu \in \mathbb{T}_0$ and generic $X\in \mathfrak{X}$ the Severi sets $\mathcal{S}\mathcal{V}_X(\mu,1)$ and $\mathcal{S}\mathcal{V}_X(\mu,2)$ have cardinals $$\#\mathcal{S}\mathcal{V}_X(\mu,1)=6\,\textrm{-}\,2I_0(\mu)\quad\textrm{and}\quad
 \#\mathcal{S}\mathcal{V}_X(\mu,2)= 27\,\textrm{-}\,14I_0(\mu)+2I_0(\mu)^2\,\textrm{-}\,3I_1(\mu).$$
 \begin{enumerate}
 	\item The condition $|2M\,\textrm{-}\,S+(1+(\,\textrm{-}\,1)^S m| \geq 4$ is equivalent to $min\{\mu_i\}\geq 2$, i.e.: $\left(I_0(\mu),I_1(\mu)\right)=(0,0)$. For such type $\#\mathcal{S}\mathcal{V}_X(\mu,2)=\#\mathcal{P}ot_X(\alpha,2)$ and the proof follows the lines of §$5$ Remark $7$.
 	\item In case $\left(I_0(\mu),I_1(\mu)\right)\neq (0,0)$ but $I_0(\mu)\leq 1$, for generic $X\in \mathfrak{X}$ we only have $\#\mathcal{S}\mathcal{V}_X(\mu,2)=27\,\textrm{-}\,12I_0(\nu)\,\textrm{-}\,3I_1(\nu)$ spectral data as above. On the other hand, for each index $i_0$ such that $\mu_{i_0}=1$, there is an exceptional $\it{ht}$ cover with three theta-characteristics giving rise to three elements in $\mathcal{P}ot_X(\alpha,2)$ (according to Lemma §$5.10$). At last, let us assume $\mu_{j_0}=0$ for some index $j_0$, and let $\nu\in \mathbb{T}_0$ denote the type with $\nu_i=\mu_i+2\delta_{i,j_0}$, hence $\nu^{(2)}=\mu^{(2)}+4$, $I_0(\nu)= 0$ and for generic $X\in \mathfrak{X}$, $\#\mathcal{S}\mathcal{V}_X(\nu,1)=6$. Moreover, according to Lemma §$5.9$, any $\it{ht}$ cover coming from $\mathcal{S}\mathcal{V}_X(\nu,1)$ can be equipped with two theta-characteristics giving rise to a total of $2.6=12$ more elements in $\mathcal{P}ot_X(\alpha,2)$. We thus complete to a total of $27$ elements in $\mathcal{P}ot_X(\alpha,2)$, the maximum possible. Their spectral curves have arithmetic genus $g_\alpha$, the first ones, and $g_\alpha+1$ the last ones.
 	\item If $I_0(\mu)=2$, say $\mu_2=\mu_3=0$ hence $\mu_0$ odd, but $\mu_1\in2\mathbb{N}^*$, for generic $X\in \mathfrak{X}$ we know $\#\mathcal{S}\mathcal{V}_X(\mu,2)=7\,\textrm{-}\,3I_1(\mu)$. The exceptional $\it{ht}$ cover of type $\tau=(\mu_0,\mu_1,2,2)$ has degree $n$, arithmetic genus $g_\alpha+2$ and is naturally equipped (see Lemma §$5.11$) with four theta-characteristics, giving rise to four elements of $\mathcal{P}ot_X(\alpha,2)$. On the other hand, if $\nu=(\mu_0,\mu_1,2,0)$ or $\nu=(\mu_0,\mu_1,0,2)$, for generic $X \in \mathfrak{X}$ we know $\#\mathcal{S}\mathcal{V}_X(\nu,1)=4$ and any $\it{ht}$ cover coming from this Severi set can be equipped with two theta-characteristics giving rise to two elements in $\mathcal{P}ot_X(\alpha,2)$. Summing up, we have $7\,\textrm{-}\,3I_1(\mu) +4+2.4+2.4=27\,\textrm{-}\,3I_1(\mu)$ elements in $\mathcal{P}ot_X(\alpha,2)$, (for generic $X \in \mathfrak{X}$), to which we add as precedently, $3I_1(\mu)$ more (in fact $3$) if $\mu_0=1$. Total $27$. The corresponding arithmetic geni range from $g_\alpha$ to $g_\alpha +2$.
 	\item At last the case $I_0=3$ (i.e.: $\mu=(\mu_0,0,0,0)$ with $\mu_0$ odd). For generic $X\in \mathfrak{X}$ we just have $\#\mathcal{S}\mathcal{V}_X(\mu,2)=3\,\textrm{-}\,3I_1(\mu)$. There are as precedently, three exceptional $\it{ht}$ covers of types $(\mu_0,0,2,2),(\mu_0,2,0,2)$ and $(\mu_0,2,2,0)$, each one equipped with four theta-characteristics giving rise to four elements in $\mathcal{P}ot_X(\alpha,2)$. We have as well that any $\nu\in \{(\mu_0,2,0,0),(\mu_0,0,2,0),(\mu_0,0,0,2)\}$ satisfies $\nu^{(2)}=\mu^{(2)}+4$, hence for generic $X\in \mathfrak{X}$ we have $\#\mathcal{S}\mathcal{V}_X(\nu,1)=6\,\textrm{-}\,2I_0(\nu)=2$. Moreover, each $\it{ht}$ cover coming from the three latter Severi sets is equipped with two theta-characteristics giving rise to two elements in $\mathcal{P}ot_X(\alpha,2)$. Summing up, for generic $X\in \mathfrak{X}$ we have obtained $3\,\textrm{-}\,3I_1(\mu)+3.4+3.2.2=27\,\textrm{-}\,3I_1(\mu)$ potentials in $\mathcal{P}ot_X(\alpha,2)$, to which we add $3I_1(\mu)$ (in fact $3$) if $\mu_0=1$. Total $=27$. We can easily check that the arithmetic geni  range again from $g_\alpha$ to $g_\alpha +2$.  $\blacksquare$
 	\end{enumerate}
 
 \begin{cor}\hspace*{2mm}
 	 Any spectral data coming from $\mathcal{P}ot_X(\alpha,2)$, say $(\pi: (\Gamma,p) \to (X,\omega_0),\xi)$, has type $\nu=\mu+ 2\gamma$ and arithmetic genus $g$, where $\gamma \in \mathbb{N}^4$ satisfies $\mu.\gamma+\gamma^{(2)}\leq 2$ and $g_\alpha  \leq g \leq g_\alpha +2$, while $\xi=\xi_{j,k}(p'+p''\,\textrm{-}\,2p)$ for a suitable pair $\{p',p''\}$ of Weierstrass points. More precisely$:$
 		\begin{enumerate} 
 			\item if $|2M\,\textrm{-}\,S+(1+(\,\textrm{-}\,1)^S m| \geq 4$,  each $\pi$ has type $\mu$, arithmetic genus $g_\alpha$ and $\xi=\xi_{j,k};$
 				\item if $|2M\,\textrm{-}\,S+(1+(\,\textrm{-}\,1)^S m| < 4$ instead, $\pi$ has type $\nu=\mu+2\gamma$ and arithmetic genus $g_\alpha \leq g \leq g_\alpha +2$ satisfying$:$
 			\begin{enumerate}
 			\item $g=g_\alpha$ if and only if $\gamma^{(2)}=0$ $($i.e.$:\gamma=\overrightarrow{0})$, if and only if $\xi=\xi_{j,k};$
 			
 			\item $g=g_\alpha +1$ if and only if $\gamma^{(2)}=1;$
 			
 			\item $g=g_\alpha +2$ if and only if $\gamma^{(2)}=2.$
 			\end{enumerate}
 			
 		\end{enumerate}

 	\end{cor}
 We conclude this final chapter with partial results about the spectral data of $ \mathcal{P}ot_X(\alpha,d)$, for any $(\alpha,d) \in \mathbb{N}^4 \times \mathbb{N}$. They are given in terms of the unique type $\mu \in \mathbb{T}_0$ and $(j,k) \in \mathbb{Z}_2 \times \mathbb{Z}_4$ such that $\alpha = \mathcal{C}^{j,k}(\mu)$. Their proof, to be given elsewhere, is based on formula \cite{TV2}, Prop.$3.2$ (see also \cite{F}), which calculates the coefficient of each term $\wp(x\,\textrm{-}\,\omega_i)$ in $u_\xi(x)$. We also deduce partial information about the cardinals $\#\mathcal{P}ot_X(\alpha,d)$ and $\#\mathcal{S}\mathcal{V}_X(\mu,d)$, for any $(\alpha,\mu,d) \in \mathbb{N}^4 \times \mathbb{T}_0 \times \mathbb{N}$ and generic $X\in \mathfrak{X}$, leading to a conjectural recursive formula for both cardinals.
 
 \begin{prop}\hspace*{2mm}
 	
  The spectral curve of any potential in $\mathcal{P}ot_X(\alpha,d)$ is a $\it{ht}$ cover of degree $n$, type $\nu$ and arithmetic genus $g$ bounded as follows: $2n=\sum_i\alpha_i(\alpha_i+1)+4d$, for some $d\in \mathbb{N}$, $\nu\,\textrm{-}\,\mu \in 2\mathbb{N}^4$, $\nu^{(2)}\leq \mu^{(2)}+ 4d$ and $g_\alpha\leq g \leq g_\alpha+2\sqrt{d}$. Moreover, $\mathcal{P}ot_X(\alpha,d)$ has cardinal
  
  $$\#\mathcal{P}ot_X(\alpha,d)= \#\mathcal{S}\mathcal{V}_X(\mu,d) +\sum_{\nu,\gamma,l} \#\mathcal{S}\mathcal{V}_X(\nu,l)\Pi_{i=0}^3 C_{\gamma_i}^{\nu_i},$$
  
  where $(\nu,\gamma,l)$ ranges over all triplets in $\mathbb{T}_0\times \mathbb{N}^4 \times \mathbb{N}$, satisfying $\nu=\mu +2\gamma$, $\gamma \neq \overrightarrow{0}$ and $\mu.\gamma +\gamma^{(2)} + l= d$. In particular, $\mu= (1,0,0,0)$ implies $\alpha=\overrightarrow{0}$, in which case $\#\mathcal{P}ot_X(\overrightarrow{0},d)$ is given in terms of cardinals $\#\mathcal{S}\mathcal{V}_X(\nu,l)$ with $l<d$.
 
 	\end{prop}
 
 	\begin{cor}\hspace*{2mm}
 		
 	Let $\alpha=\mathcal{C}^{j,k}(\mu)$ and assume $|2M\,\textrm{-}\,S+(1+(\,\textrm{-}\,1)^S m| \geq 2d$ $($i.e.: $min\{\mu_i\} \geq d)$. Then $\#\mathcal{P}ot_X(\alpha,d)= \#\mathcal{S}\mathcal{V}_X(\mu,d)$, and for any spectral data $(\pi,\xi)$ coming from $\mathcal{P}ot_X(\alpha,d)$, the  $\it{ht}$ cover $\pi$ has degree $n$, type $\mu$ and arithmetic genus $g_\alpha$, while $\xi=\xi_{j,k}$.
 		
\end{cor}
 	
 	\begin{remark} $($Conjectural recursive formula for both cardinals$)$\hspace*{2mm}
 		
 		Could we also prove $($just as for $d\leq 2)$ that $\#\mathcal{P}ot_X(\alpha,d)$ is independent of $\alpha\in \mathbb{N}^4$, hence equal to $$\#\mathcal{P}ot_X(\overrightarrow{0},d)= \sum_{\nu,\gamma,l} \#\mathcal{S}\mathcal{V}_X(\nu,l)\Pi_{i=0}^3 C_{\gamma_i}^{\nu_i},$$
 		
 	with $\nu:=(2\gamma_0+1,2\gamma_1,2\gamma_2,2\gamma_3)$, $\gamma \neq \overrightarrow{0}$ and $\gamma_0 +\gamma^{(2)} + l= d$, we would deduce its value from the knowledge of all cardinals $\#\mathcal{S}\mathcal{V}_X(\nu,l)$ with $l<d$.
 	At last, we would deduce from Prop. §$5.14$, that
 	
 	$$\#\mathcal{S}\mathcal{V}_X(\mu,d)=\#\mathcal{P}ot_X(\overrightarrow{0},d)\,\textrm{-}\,\sum_{\nu,\gamma,l} \#\mathcal{S}\mathcal{V}_X(\nu,l)\Pi_{i=0}^3 C_{\gamma_i}^{\nu_i},$$
 	
 	where $(\nu,\gamma,l)$ ranges over all triplets in $\mathbb{T}_0\times \mathbb{N}^4 \times \mathbb{N}$, satisfying $\nu=\mu +2\gamma$, $\gamma \neq \overrightarrow{0}$ and $\mu.\gamma +\gamma^{(2)} + l= d$.

 	\end{remark}
	\appendix
	\section{Contraction data and families of blowing-ups}
	\begin{defi}\hspace*{2mm}

Let $(\overrightarrow{t})	:=(\overrightarrow{t_0},\overrightarrow{t_1},\overrightarrow{t_2},\overrightarrow{t_3}) \in \Pi_{i=0}^3 \mathbb{T}_i$ and denote $\widetilde{\Gamma}_{(\overrightarrow{t})}:=\left(\widetilde{\Gamma}_{\overrightarrow{t}_0},\widetilde{\Gamma}_{\overrightarrow{t}_1},\widetilde{\Gamma}_{\overrightarrow{t}_2},\widetilde{\Gamma}_{\overrightarrow{t}_3}\right)$ the corresponding vector of exceptional curves of $\widetilde{S}_X$. We will call it contraction data if and only if they are disjoint one from the other, and each $\widetilde{r}_k$ ($0\leq k\leq 3)$ intersects at least one of them. In other words$:$
	 		
		$$\forall i\neq j,\quad \big(\overrightarrow{t_i}\,\textrm{-}\,\overrightarrow{t_j} \big)^{(2)}=2\;\textrm{, and }\sum_{i=0}^3 \overrightarrow{t_i}\quad\textrm{has positive coordinates}\;.$$
	 		We then let, $\psi_{ (\overrightarrow{t})}: \widetilde{S}_X \to \underline{S}_{(\overrightarrow{t})}$ denote the contraction of the four exceptional curves $\big\{\widetilde{\Gamma}_{\overrightarrow{t_i}}\big\}$, and $\psi_{ (\overrightarrow{s})}: \underline{S}_{(\overrightarrow{t})} \to \underline{S}_{(\overrightarrow{s},\overrightarrow{t})}$ the contraction $($of the images$)$ of $\big\{\widetilde{s}_i\big\}$.

	 	\end{defi}
	 	
	 	\begin{remark}\hspace*{2mm}
	 		
	 		We will henceforth fix $\mu \in \mathbb{T}_0$ and complete $\widetilde{\Gamma}_\mu$ to a contraction data as follows$:$ 
	 		$$\widetilde{\Gamma}_{(\overrightarrow{t})}=\big(\widetilde{\Gamma}_\mu,\widetilde{\Gamma}_\mu^{(1,1,0,0)},\widetilde{\Gamma}_\mu^{(1,0,1,0)},\widetilde{\Gamma}_\mu^{(1,0,0,1)}\big) \;.$$
	 		We will also consider the following curve vectors$:$
	 		$$\big(\widetilde{F}_i\big):=\big(\widetilde{\Gamma}_\mu^{(1,1,1,1)},\widetilde{\Gamma}_\mu^{(0,0,1,1)},\widetilde{\Gamma}_\mu^{(0,1,0,1)},\widetilde{\Gamma}_\mu^{(0,1,1,0)}\big)$$
	 		and  $$\big(\widetilde{G}_i\big):=\big(\widetilde{\Gamma}_\mu^{(2,0,0,0)},\widetilde{\Gamma}_\mu^{(1,\,\textrm{-}1,0,0)},\widetilde{\Gamma}_\mu^{(1,0,\,\textrm{-}1,0)},\widetilde{\Gamma}_\mu^{(1,0,0,\,\textrm{-}1)}\big)\;.$$
	 		
	 		In case $\mu_j=0$ for some $j\neq 0$, the coordinate $\widetilde{G}_j$ should be replaced by $\widetilde{r_j}+\widetilde{\Gamma}_{\overrightarrow{t_j}}$.
	 		
	 	\end{remark}
	 	\begin{lem}\hspace*{2mm}
	 		
	 		Given $\mu \in \mathbb{T}_0 $, we choose the above contraction data $\widetilde{\Gamma}_{(\overrightarrow{t})}$ and curve vectors $\big(\widetilde{F}_i\big)$ and $\big(\widetilde{G}_i\big)$. Let $\Delta$ denote the diagonal of $\mathbb{P}^1\times \mathbb{P}^1$ and $p_i$ the unique point in $\widetilde{s}_i\cap \widetilde{C}_0$ $(i=0,\cdots,3)$. Then, up to an automorphism of $\mathbb{P}^1\times \mathbb{P}^1$, we can identify$:$
	 		\begin{enumerate}
	 			\item $\underline{S}_{(\overrightarrow{s},\overrightarrow{t})}$ with $\mathbb{P}^1\times \mathbb{P}^1$, and the image of the pair $(\widetilde{C}_0,\,p_i)$, with $(\Delta,\,(p_i,p_i));$
	 			
	 			\item the morphism $\psi_{(\overrightarrow{s})}:\underline{S}_{(\overrightarrow{t})} \to \underline{S}_{(\overrightarrow{s},\overrightarrow{t})} \cong \mathbb{P}^1\times \mathbb{P}^1$ with the blowing-up of four points $\big\{(p_i,p_i)\} \subset \mathbb{P}^1\times \mathbb{P}^1;$
	 			\item the exceptional fibers of $\psi_{(\overrightarrow{s})}$ with the projections $\psi_{(\overrightarrow{t})}(\widetilde{s}_i) \subset \underline{S}_{(\overrightarrow{t})}$, $i=0,\cdots,3$.
	 		\end{enumerate}

	 	\end{lem}
	 	\textbf{Proof.}\\
	 	\textbf{1)} - The curves $\big\{\widetilde{s}_i\big\}$ project onto exceptional ones in the surface $\underline{S}_{(\overrightarrow{t})}$, but remain disjoint and may be contracted. Any other exceptional curve intersects $\{\widetilde{s}_i\}$, for some index $i$, hence its image in $\underline{S}_{(\overrightarrow{s},\overrightarrow{t})}$ is no longer exceptional. The same happens to the $(\textrm{-}2)$-curves $\big\{\widetilde{r_k}\big\}$: for any $k$ there exists $i$ such that $\widetilde{r}_k $ intersects $\widetilde{\Gamma}_{\overrightarrow{t_i}}$ in $\widetilde{S}_X$ (because $\sum_i \overrightarrow{t_i}$ has positive coordinates), hence intersects $\widetilde{s}_i$ in $\underline{S}_{(\overrightarrow{t})}$. It follows that $\widetilde{r}_k.\widetilde{r}_k \geq 0$ in $\underline{S}_{(\overrightarrow{s},\overrightarrow{t})}$. The rational surface $\underline{S}_{(\overrightarrow{s},\overrightarrow{t})}$ has a rank-$2$ Picard group and no negative curve. The latter implies $\underline{S}_{(\overrightarrow{s},\overrightarrow{t})}$ is isomorphic to $\mathbb{P}^1\times \mathbb{P}^1$. Knowing that the image of $\widetilde{C}_0$ in $\underline{S}_{(\overrightarrow{s},\overrightarrow{t})} \simeq \mathbb{P}^1\times \mathbb{P}^1$ is numerically equivalent to $\Delta $, we deduce the existence of a unique automorphism of $\mathbb{P}^1\times \mathbb{P}^1$ identifying the image of each $p_i$ with $(p_i,p_i)$, as well as the image of $\widetilde{C}_0$, numerically equivalent to $\Delta$, with $\Delta$ (because $\Delta^2=2$).\\
	 	\textbf{2)} et \textbf{3)} - Recall that $\psi_{(\overrightarrow{s},\overrightarrow{t})}$ is obtained as the contraction of curves $\big\{\widetilde{\Gamma}_{\overrightarrow{t_i}}\big\} $, and the images of $\big\{\widetilde{s}_i\big\}$.  The last two identifications follow directly from the universal property of the blowing-up (cf. \cite{B}, Proposition II.$8$).
	 	$\quad \quad \quad \blacksquare$
	 	
	 	\begin{remark}\hspace*{2mm}
	 		
	 		According to the above identifications, for each index $i$, the strict transforms in $\underline{S}_{(\overrightarrow{t})}$ of the fibers $\mathbb{P}^1\times \{p_i\}$ and $ \{p_i\} \times \mathbb{P}^1$, denoted $F_i$ and $G_i$, are exceptional curves. We check via the projection formula for $\psi_{(\overrightarrow{t})}$, that their inverse images in $\widetilde{S}_X$,   $\widetilde{F_i}:=\psi_{(\overrightarrow{t})}^\ast(F_i) \;\textrm{et}\;\widetilde{G}_i:=\psi_{(\overrightarrow{t})}^\ast(G_i)$, satisfy the following properties$:$
	 		
	 		$$\forall i,j, \quad \widetilde{F}_i.\widetilde{\Gamma}_{\overrightarrow{t_j}}=0=\widetilde{G}_i.\widetilde{\Gamma}_{\overrightarrow{t_j}},\quad \widetilde{F}_i.\widetilde{G}_j=1\,\textrm{-}\,\delta_{i,j} \quad\textrm{et}\quad \widetilde{F}_i.\widetilde{s}_j=\delta_{i,j}=\widetilde{G}_i.\widetilde{s}_j\, .$$
	 	\end{remark}
	 	
	 	\begin{defi}\hspace*{2mm}
	 		
	 		Consider the curve vectors $\widetilde{F}_{(\overrightarrow{t})}:=\big(\widetilde{F}_0,\cdots,\widetilde{F}_3\big)$ and $\widetilde{G}_{(\overrightarrow{t})}:=\big(\widetilde{G}_0,\cdots,\widetilde{G}_3\big)$ defined above. We will say $\Big(\widetilde{\Gamma}_{(\overrightarrow{t})},\widetilde{F}_{(\overrightarrow{t})},\widetilde{G}_{(\overrightarrow{t})}\Big)$ is an exceptional configuration associated to the contraction data $\widetilde{\Gamma}_{(\overrightarrow{t})}$.
	 		
	 	\end{defi}

	 	Fix $\mu \in \mathbb{T}_0$ and let $ (X,\omega_0)$ be the elliptic curve defined by choosing a fourth point $p_0 \in \mathfrak{X}:=\mathbb{P}^1\setminus \{p_1,p_2,p_3\}$. The canonical projection $\widetilde{S}_X \to \mathbb{P}^1$ (see §$2.1.3.$) defines an isomorphism between $\widetilde{C}_0$ and $\mathbb{P}^1$, identifying the unique point in $\widetilde{s}_i \cap \widetilde{C}_0$ with $p_i$. Let $\psi_{\mu}:\widetilde{S}_X \to \underline{S}_{\mu}$ denote the contraction of the divisor $\widetilde{\Gamma}_{\mu}+\widetilde{s}_0$, and $\underline{C}_{0}$ the projection of $\widetilde{C}_0 $ in $ \underline{S}_{\mu}$. The following result is in order.
	 	
	 	\begin{theo}\hspace*{2mm}
	 		
	 		The couple $\big(\underline{S}_{\mu},\underline{C}_{0}\big)$ is independent of $p_0 \in \mathfrak{X}$ $($hence of $X)$.\\
	 		In other words, $\big\{\widetilde{S}_X , X \in \mathfrak{X}\big\}$ is a family of blowing-ups of $\underline{S}_{\mu}$ along $\underline{C}_{0}$.	
	 	\end{theo}
	 	\textbf{Proof.}\\
	 	Consider the above exceptional configuration $\Big(\widetilde{\Gamma}_{(\overrightarrow{t})},\widetilde{F}_{(\overrightarrow{t})},\widetilde{G}_{(\overrightarrow{t})}\Big)$, associated to the contraction data $\widetilde{\Gamma}_{(\overrightarrow{t})}=\big(\widetilde{\Gamma}_\mu,\widetilde{\Gamma}_\mu^{(1,1,0,0)},\widetilde{\Gamma}_\mu^{(1,0,1,0)},\widetilde{\Gamma}_\mu^{(1,0,0,1)}\big)$.\\
	 	Let $\psi_{(\overrightarrow{t})}: \widetilde{S}_X \to \underline{S}_{(\overrightarrow{t})}$ and $\psi_{(\overrightarrow{s})}:  \underline{S}_{(\overrightarrow{t})} \to \underline{S}_{(\overrightarrow{s},\overrightarrow{t})}$, denote the contractions of $\big\{\widetilde{\Gamma}_{\overrightarrow{t_i}}\big\}$ and (the images of) $\big\{\widetilde{s}_i\big\}$, respectively. According to the preceding lemma, $\psi_{(\overrightarrow{s})}$ is isomorphic to the blowing up of points $\big\{(p_i,p_i\big)\} \subset \mathbb{P}^1 \times \mathbb{P}^1$. In particular, the exceptional fiber $\psi_{(\overrightarrow{s})}^\ast((p_i,p_i)) \subset \underline{S}_{(\overrightarrow{t})}$, denoted $E_i$, is equal to the image of $\widetilde{s}_i$ by $\psi_{(\overrightarrow{t})}$. On the other hand, the images $\psi_{(\overrightarrow{t})}(\widetilde{F}_i)$ and $\psi_{(\overrightarrow{t})}(\widetilde{G}_i)$ coincide with the strict transforms (by $\psi_{(\overrightarrow{s})}$) of $\mathbb{P}^1 \times \{p_i\}$ and $\{p_i\} \times \mathbb{P}^1$ respectively, while the image of $\widetilde{C}_0$ coincides with the strict transform of the diagonal $\Delta$. The diagrams below represent the corresponding morphisms and curves$:$ 
	 	
	 	\begin{displaymath}
	 		\xymatrix{
	 			\widetilde{C}_0 \ar[d]^{\psi_{(\overrightarrow{t})}}	& \widetilde{F}_i \ar[d]^{\psi_{(\overrightarrow{t})}} &   \widetilde{G}_i \ar[d]^{\psi_{(\overrightarrow{t})}} & \widetilde{s}_i \subset \widetilde{S}_X \ar[d]^{\psi_{(\overrightarrow{t})}}\\
	 			\psi_{(\overrightarrow{t})}(\widetilde{C}_0) \ar[d]^{\psi_{(\overrightarrow{s})}}	& \psi_{(\overrightarrow{t})}(\widetilde{F}_i) \ar[d]^{\psi_{(\overrightarrow{s})}} & \psi_{(\overrightarrow{t})}(\widetilde{G}_i) \ar[d]^{\psi_{(\overrightarrow{s})}} &  E_i:=\psi_{(\overrightarrow{t})}^\ast((p_i,p_i)) \subset \underline{S}_{(\overrightarrow{t})} \ar[d]^{\psi_{(\overrightarrow{s})}}\\
	 			\Delta &  \mathbb{P}^1 \times \{p_i\}   &  \{p_i\}\times  \mathbb{P}^1 &(p_i,p_i) \in \mathbb{P}^1 \times \mathbb{P}^1  
	 		}
	 	\end{displaymath}

	 	The exceptional fiber $E_i$ is marked at the three points of its intersection with the images of the triplet of curves $\big(\widetilde{C}_0,\widetilde{F}_i,\widetilde{G}_i\big)$. Recalling the identifications of each $E_i$ with $\widetilde{s}_i$, and of $\widetilde{s}_i$ with $\mathbb{C}\cup \{\infty\}$, we deduce an isomorphism $E_i\simeq \mathbb{C}\cup \{\infty\}$, identifying  the triplet of marked points of $E_i$, with the following triplet. In case $i=0,1,2$ or $3$, it is equal to:
	 	\begin{enumerate}
	 		\item $i=0:$ $\big(\infty,n_{\mu+(1,1,1,1)},n_{\mu+(2,0,0,0)}\big)\;;$
	 		 
	 		\item $i=1:$ $\big(\infty,n_{\mu+(0,0,1,1)},n_{\mu+(1,\,\textrm{-}1,0,0)}\big)\quad$ if $\mu_1\neq 0$, or $\big(\infty,n_{\mu+(0,0,1,1)},0\big)\quad$ if $\mu_1=0\;;$ 
	 		
	 		\item $i=2:$ $\big(\infty,n_{\mu+(0,1,0,1)},n_{\mu+(1,0,\,\textrm{-}1,0)}\big)\quad$ if $\mu_2\neq 0$, or $\big(\infty,n_{\mu+(0,1,0,1)},0\big)\quad$ if $\mu_2=0\;;$
	 		
	 		\item $i=3:$ $\big(\infty,n_{\mu+(0,1,1,0)},n_{\mu+(1,0,0,\,\textrm{-}1)}\big)\quad$ if $\mu_3\neq 0$, or $\big(\infty,n_{\mu+(0,1,1,0)},0\big)\quad$ if $\mu_3=0\;$.
	 	\end{enumerate} 
	 	On the other hand, for any $i=0,\cdots,3$, $\psi_{(\overrightarrow{t})}$ contracts $\widetilde{\Gamma}_{\overrightarrow{t_i}}$ to the point of $E_i$ corresponding, respectively, to the values $n_\mu, n_{\mu+(1,1,0,0)},n_{\mu+(1,0,1,0)},n_{\mu+(1,0,0,1)}$.\\
	 	In brief, up to isomorphism the family $\{\widetilde{S}_X,X\in \mathfrak{X}\}$ is obtained as follows.
	 	\begin{enumerate}
	 		\item After blowing-up $\{(p_i,p_i)\} \subset \mathbb{P}^1 \times \mathbb{P}^1$ we mark the exceptional fiber over each $(p_i,p_i)$, denoted $E_i$, at its intersection with the strict transforms of the triplet of curves $\big(\mathbb{P}^1 \times \{p_i\}, \{p_i\} \times  \mathbb{P}^1, \Delta\big).$
	 		
	 		\item We then identify $E_i$ with $ \mathbb{C}\cup \{\infty\}$, via the unique isomorphism identifying the above triplet of marked points of $E_i$ with:
	 		\begin{enumerate}
	 			\item $(\infty,n_{\mu+(1,1,1,1)},n_{\mu+(2,0,0,0)})$, if $i=0\;;$
	 			\item $(\infty,n_{\mu+(0,0,1,1)}, n_{\mu+(1,\,\textrm{-}1,0,0)})$, if $i=1$ and $\mu_1\neq 0\;;$
	 			\item $(\infty,n_{\mu+(0,0,1,1)}, 0)$, if $i=1$ and $\mu_1= 0\;;$
	 			\item $(\infty,n_{\mu+(0,1,0,1)}, n_{\mu+(1,0,\,\textrm{-}1,0)})$, if $i=2$ and $\mu_2\neq 0\;;$
	 			\item $(\infty,n_{\mu+(0,1,0,1)}, 0)$, if $i=2$ and $\mu_2= 0\;;$ 
	 			\item $(\infty,n_{\mu+(0,1,1,0)}, n_{\mu+(1,0,0,\,\textrm{-}1)})$, if $i=3$ and $\mu_3\neq 0\;;$
	 			\item $(\infty,n_{\mu+(0,1,1,0)}, 0)$, if $i=3$ and $\mu_3= 0\;.$
	 		\end{enumerate}
	 		\item At last, for each $i=0,\cdots,3$ we blow up the points of $E_i$ corresponding to the values $n_\mu,n_{\mu+(1,1,0,0)},n_{\mu+(1,0,1,0)}$ and $n_{\mu+(1,0,0,1)}$, respectively.

	 	\end{enumerate} 
	 	We thus obtain a surface isomorphic to $\widetilde{S}_X$ (for $X$ associated to $p_0 \in \mathfrak{X}$).
	 	It follows that $\underline{S}_\mu \to \mathbb{P}^1\times \mathbb{P}^1$ is isomorphic to the blowing-ups of $\big\{(p_j,p_j),j=1,2,3\big\} \subset \mathbb{P}^1\times \mathbb{P}^1$, and the above marked points of $\big\{E_j,j=1,2,3\big\}$. In particular, $\underline{S}_\mu$ is independent of $p_0$, hence of $X\in \mathfrak{X}$. Moreover, $\underline{C}_0$ coincides with the strict transform of the diagonal $\Delta$. The blowing-up of $p_0 \in \underline{C}_{0}$
	 	and $n_\mu \in E_0$ is then isomorphic to $\psi_\mu : \widetilde{S}_X \to \underline{S}_\mu$.
	 	$\quad \quad \quad \blacksquare$

	 	\section{The Severi Variety $\mathcal{S}\mathcal{V}_X(\mu,2)$ for generic $X\in \mathfrak{X}$.}
	 	We will identify hereafter the family $\mathcal{S}\mathcal{V}:=\bigcup_\mathfrak{X}\mathcal{S}\mathcal{V}_X(\mu,2) \subset \bigcup_\mathfrak{X}|\gamma_X(\mu,2)|$ of Severi Varieties, with its projection in the complete linear system $|2\underline{L}_\mu|$, and consider its closure in twice the anticanonical system, denoted $\widehat{\mathcal{S}\mathcal{V}} \subset |2\underline{L}_\mu|\simeq \mathbb{P}^6(\mathbb{C})$.\\
	 	\indent Up to the natural identification between $\mathfrak{X} \subset \mathbb{P}^1$ and $\underline{C}_0\setminus \bigcup_j\underline{s}_j \subset \underline{C}_0$, the subvariety $\widehat{\mathcal{S}\mathcal{V}}$ comes with a natural projection onto $\underline{C}_0$, as well as an embedding $\widehat{\mathcal{S}\mathcal{V}} \to \underline{\Omega}^{(6)}$, sending $\underline{\Gamma} \in \widehat{\mathcal{S}\mathcal{V}}$ to the degree-$6$ effective divisor $\underline{\Gamma}\cap \underline{\Omega}$. We deduce in particular, a stratification $\widehat{\mathcal{S}\mathcal{V}}=\cup_{\mathfrak{Dec}(6)}\widehat{\mathcal{S}\mathcal{V}}^{\vec{d}}$, where $\widehat{\mathcal{S}\mathcal{V}}^{\vec{d}}$ denotes the inverse image of $\underline{\Omega}^{\vec{d}}$ with respect to the natural morphism $\widehat{\mathcal{S}\mathcal{V}} \to \underline{\Omega}^{(6)}$ (see §$3.11$). We will prove that $dim(\widehat{\mathcal{S}\mathcal{V}}^{\vec{d}})=0$ unless $ \vec{d}=(2,2,1,1,0,0)$. In other words, for generic $X\in \mathfrak{X}$ the Severi Variety $\mathcal{S}\mathcal{V}_X(\mu,2)$ can only contain (rational irreducible) curves with two simple nodes.
	 	 We gather hereafter, known basic results needed to do it, and let $\underline{p}_X \in \underline{C}_0 \setminus \bigcup_j \underline{s}_j$ denote the point corresponding to $X\in \mathfrak{X} \subset \mathbb{P}^1$.
	 	
	 	\begin{lem}\hspace*{2mm}
	 		\begin{enumerate}
	 			
	 			 			\item A divisor $\underline{\Gamma} \in \psi_\mu(|\gamma_X(\mu,2)|)$ is reducible if and only if, either $\underline{C}_0$ is one of its components, or there exists an index $i_0$ such that $\mu_{i_0}=0$ and $\underline{\Gamma} = \underline{r}_{i_0} +\underline{\Gamma}_\nu $, where $\nu \in \mathbb{T}_0$ and $\nu_i=\mu_i + 2\delta_{i_0,i}$ for any $i=0,\cdots,3$. In the latter case, the conic $g_\mu(\underline{\Gamma} )$ is tangent to the cubic $\mathcal{K}_0:=g_\mu(\underline{\Omega})$ at any point of their intersection.
	 			
	 			\item Otherwise $\underline{\Gamma} \in \psi_\mu(|\gamma_X(\mu,2)|)$ is irreducible, and  $\underline{\Gamma} \cap \underline{\Omega} \in \Omega^{\vec d}$ for some partition ${\vec d}\in \mathfrak{Dec}(6)$ with two odd terms.
	 			
	 			\item Any divisor $\underline{\Gamma}\in \widehat{\mathcal{S}\mathcal{V}}$ in the fiber over $\underline{p}_X$ is contained in $\psi_\mu(|\gamma_X(\mu,2)|)$, and the conic $g_\mu(\underline{\Gamma})$ is either, reducible with $\mathcal{H}_0:=g_\mu(\underline{C}_0)$ as one of its components, or smooth and tangent to $\mathcal{H}_0$ at $g_\mu(\underline{p}_X)$.
 			
	 				\item For any $\underline{p}_X\in \underline{C}_0\setminus \bigcup_j\underline{s}_j$ and $\underline{\Gamma}_1, \underline{\Gamma}_2 \in \psi_\mu(|\gamma_X(\mu,2)|)$, their projections are both tangent to $\mathcal{H}_0$ at $g_\mu(\underline{p}_X)$, where they intersect with multiplicity $\geq 3$.
	 				
	 				\item Let $\mathcal{K}\subset \mathbb{P}^2(\mathbb{C})$ be an irreducible cubic, $q\in \mathcal{K}$ a smooth point, $\mathcal{H}_q$ the tangent line to $\mathcal{K}$ at $q$ and $\mathcal{H}$ any line not going through $q$. Then, 
	 				\begin{enumerate}
	 					\item if $q$ is an inflection point of $\mathcal{K}$, any smooth conic $\mathcal{C}$ has multiplicity of intersection $I_q(\mathcal{K},\mathcal{C})\leq 2$ with $\mathcal{K}$ at $q;$ 
	 					
	 					\item otherwise, the conics satisfying $I_q(\mathcal{K},\mathcal{C})\geq 4$ make a pencil, with only one irreducible element tangent to $\mathcal{H}$. Moreover, $2\mathcal{H}_q$ is its unique reducible element.
	 					
	 				\end{enumerate}
	 		
	 		\end{enumerate}
	 		\end{lem}
 		
 		In order to prove that $\widehat{\mathcal{S}\mathcal{V}}^{\vec{d}}$ is finite unless $\vec{d}=(2,2,1,1,0,0)$, we will consider different cases, according to whether or not the cubic $\mathcal{K}_0$ is irreducible and has inflection points outside its intersection with $\mathcal{H}_0$. The latter happens if and only if $I_\mu(0):= \# \{i,\mu_i=0\}\leq 1$ and $\mu$ is not a multiple of $(0,1,1,1)$. In the particular case with $I_\mu(0)=1$ and $\mu$ multiple of $(0,1,1,1)$, up to the action of $\mathbb{P}GL(3,\mathbb{C})$, we have explicit equations for the couple $(\mathcal{K}_0,\mathcal{H}_0)$ allowing to check our claim.  When $I_\mu(0)=2$ we also have explicit equations but yet we need to distinguish between the cases $\mu_0 \geq 3$ and $\mu_0=1$. At last, when $I_\mu(0)=3$ the cubic $\mathcal{K}_0$ decomposes as a union of lines and our claim is trivial.\\

 		\begin{lem}\hspace*{2mm}
 			
 			If $I_\mu(0) =0$, or $I_\mu(0) =1$ but $\mu \notin \left\{(0,2h+1,2h+1,2h+1), h\in \mathbb{N}\right\}$, the cubic $\mathcal{K}_0$ is irreducible and has at most one inflection point in $\mathcal{K}_0\cap \mathcal{H}_0$.
 			\end{lem}
 		
 		\textbf{Proof.} \\
 		\indent 
 		When $I_\mu(0)=0$ the cubic $\mathcal{K}_0$ is smooth and has therefore $9$ inflection points.\\
 		\indent If $I_\mu(0)=1$ instead, it is known to be irreducible but with a node, hence it has $3$ inflection points, one of which at most is in $\mathcal{H}_0$, due to the assumption \\
 		$\mu \notin\{(0,2h+1,2h+1,2h+1), h\in \mathbb{N}\}$ (see §$3.8\,$-$\,3)$. $\blacksquare$
 		
 		\begin{prop} $(I_\mu(0) \leq 1$ and $\mu \notin \left\{(0,2h+1,2h+1,2h+1), h\in \mathbb{N}\right\})$\hspace*{2mm}

 			For any $\vec{d} \in \mathfrak{Dec}(6)$ other than $\vec{d}=(2,2,1,1,0,0)$, the set $\widehat{\mathcal{S}\mathcal{V}}^{\vec{d}}$ is finite.\\
 			
 			\end{prop}
 		\textbf{Proof.}
 		
 		According to §$3$ Remark $3$, the stratification $\widehat{\mathcal{S}\mathcal{V}}=\cup_{\mathfrak{Dec}(6)}\widehat{\mathcal{S}\mathcal{V}}^{\vec{d}}$ runs over all $\vec{d} \in \mathfrak{Dec}(6)$ with two odd terms. We first consider the cases $\vec{d}=(3,2,1,0,0,0), (4,1,1,0,0,0)$ and $(5,1,0,0,0,0)$. For any $X\in \mathfrak{X}$ and $\underline{\Gamma} \in \mathcal{S}\mathcal{V}_X^{\vec{d}} (\mu,2)$, the divisor $\underline{\Gamma}\cap \underline{\Omega}$ of $\underline{\Omega}$ is equal to $3p_1+2p_2+p_3$, $4p_1+p_2+p_3$ or $5p_1+p_2$, respectively (with $p_j \neq p_k$ if $j\neq k$). We deduce in particular a natural morphism  $\underline{\Gamma} \in \widehat{\mathcal{S}\mathcal{V}}^{\vec{d}} \mapsto p_1 \in \underline{\Omega}$.
 		 Let us denote $q_j:=g_\mu(p_j)$ for any $j$, and assume that the latter morphism is constant and $\widehat{\mathcal{S}\mathcal{V}}^{\vec{d}}$ is not finite.\\
 		 \indent For $\vec{d}=(4,1,1,0,0,0)$ or $(5,1,0,0,0,0)$ the set of conics $\{g_\mu(\underline{\Gamma}), \underline{\Gamma} \in \widehat{\mathcal{S}\mathcal{V}}^{\vec{d}}\}$ should be contained in the pencil of conics intersecting $\mathcal{K}_0$ at $q_1$ with multiplicity $\geq 4$. Such pencil has only two elements tangent to the line $\mathcal{H}_0$. Contradiction!\\
 		\indent In case $d=(3,2,1,0,0,0)$ instead, we have a second natural projection, namely $\underline{\Gamma}\in \widehat{\mathcal{S}\mathcal{V}}^{\vec{d}} \mapsto p_2 \in \underline{\Omega}$, which can not be constant either; hence, it must be surjective. We can therefore pick $\underline{\Gamma} \in \widehat{\mathcal{S}\mathcal{V}}^{\vec{d}}$ in the fiber over $p_2\in \underline{\Omega}\cap \underline{C}_0$. Its projection $g_\mu(\underline{\Gamma})$ should then be a conic tangent at $q_2$ to both, $\mathcal{H}_0$ and $\mathcal{K}_0$, while satisfying $I_{q_1}(g_\mu(\underline{\Gamma}), \mathcal{K}_0)\geq 3$. This can only happen if $g_\mu(\underline{\Gamma})$ decomposes as $2\mathcal{H}_{q_1}$, twice the tangent line to $\mathcal{K}_0$ at $q_1$, an option we can rule out by changing $p_2$ to another point in $\underline{\Omega}\cap\underline{C}_0$.\\
 		We have thus proved that the first projection $\underline{\Gamma}\in \widehat{\mathcal{S}\mathcal{V}}^{\vec{d}} \mapsto p_1 \in \underline{\Omega}$ is surjective in the above three cases.\\
 		 Now choose $\underline{p}_1\in \underline{\Omega}\setminus \underline{C}_0$, such that $q_1:=g_\mu(\underline{p}_1)$ is an inflection point of $\mathcal{K}_0$, and $\underline{\Gamma} \in \widehat{\mathcal{S}\mathcal{V}}^{\vec{d}}$ in the fiber over $\underline{p}_1$. The conic $g_\mu(\underline{\Gamma})$ intersects $\mathcal{K}_0$ at $q_1$ with multiplicity $\geq 3$, while being tangent to $\mathcal{H}_0$ (meaning $\#\big(g_\mu(\underline{\Gamma})\cap \mathcal{H}_0)=1\big)$. Hence, it must decompose, and be equal to $2\mathcal{H}_{q_1}$, twice the tangent line to $\mathcal{K}_0$ at $q_1$. In particular, it intersects $\mathcal{H}_0$ outside $\mathcal{H}_0 \cap \mathcal{K}_0$, say at $g_\mu(\underline{p}_X)$, for some $X\in \mathfrak{X}$. It follows that $\underline{\Gamma} \in \psi_\mu(|\gamma_X(\mu,2)|)$ and $g_\mu(\underline{\Gamma})$ is reducible. Contradiction!: such conic should have $\mathcal{H}_0$ as one of its components.\\
 		\indent Consider at last the remaining case, i.e.: $ \vec{d}=(3,3,0,0,0,0)$, for which there is a natural morphism $\widehat{\mathcal{S}\mathcal{V}}^{\vec{d}}\to \underline{\Omega}^{(2)},\; \underline{\Gamma}\mapsto \{ \underline{p}_1,\underline{p}_2\}$, where $3\underline{p}_1 + 3\underline{p}_2$ is the divisor $\underline{\Gamma}\cap \underline{\Omega}$ of $\underline{\Omega}$. Assume $\widehat{\mathcal{S}\mathcal{V}}^{\vec{d}}$ not finite, then choose an inflection point $q^*\in \mathcal{K}_0\setminus \mathcal{H}_0$ and let $\underline{p}^*:=g_\mu(q^*)\in \underline{\Omega}\setminus \underline{C}_0 $. The divisor $\underline{p}^*+\underline{\Omega}=\{\underline{p}^*+\underline{p}, \underline{p}\in \underline{\Omega}\}\subset \Omega^{(2)}$ being ample, it should intersect the image of $\widehat{\mathcal{S}\mathcal{V}}^{\vec{d}}$. We would deduce the existence of $\underline{\Gamma} \in \widehat{\mathcal{S}\mathcal{V}}^{\vec{d}}$ such that $I_{q^*}(g_\mu( \underline{\Gamma}),\mathcal{K}_0)\geq 3$, forcing $g_\mu( \underline{\Gamma})$ to decompose as $2\mathcal{H}_{q^*}$. Arguing as before leads to a contradiction. $\blacksquare$ \\
 		
 		  \underline{The last (sub-)case: $I_\mu(0) =1$ and $\mu \in\left\{(0,2h+1,2h+1,2h+1), h\in \mathbb{N}\right\}$}.\\

 		  In this sub-case, $\underline{r}_0$ is a fiber of the projection $g_\mu:\underline{S}_\mu \to |\underline{L}_\mu|^\vee\simeq \mathbb{P}^2(\mathbb{C})$, transverse to $\underline{\Omega}$ and $\underline{r}_0.\underline{\Omega}=2$. Hence, $\mathcal{K}_0$ is irreducible and has a node at $g_\mu(\underline{r}_0)$. Moreover, the coordinates $\mu_j,j=1,2,3$ being equal to each other, $\{g_\mu(\underline{s}_j),j=1,2,3\}$ are the three inflection points of $\mathcal{K}_0$. \\
 	\indent	  Modulo action of $\mathbb{P}GL(3,\mathbb{C})$, the line $\mathcal{H}_0$ and the cubic $\mathcal{K}_0$ are given as the zero-loci $$\mathcal{H}_0= \{z=0\}\,\;\textrm{and}\,\;\mathcal{K}_0= \{xyz\,\textrm{-}\,x^3\,\textrm{-}\,y^3=0\}.$$
 		  
 		  Let us also recall the following useful parameterization of $\mathcal{K}_0$:
 		  $$ [u:v]\in \mathbb{P}^1\; \mapsto\; [uv^2:u^2v:u^3+v^3] \in \mathcal{K}_0\subset \mathbb{P}^2(\mathbb{C}).$$
 		  
 		  \begin{prop} \underline{$(I_\mu(0) =1$ and $\mu \in\left\{(0,2h+1,2h+1,2h+1), h\in \mathbb{N}\right\})$}.\hspace*{2mm}
 		  	
 		  		For any $\vec{d} \in \mathfrak{Dec}(6)$ other than $\vec{d}=(2,2,1,1,0,0)$, the set $\widehat{\mathcal{S}\mathcal{V}}^{\vec{d}}$ is finite.\\
 		  	
 	  	\end{prop}
 		  	
 		  	\textbf{Proof.}
 		  	
 		  	For any $\vec{d}\in \mathfrak{Dec}(6)$ with two odd coordinates, consider again the first morphism $\underline{\Gamma} \in \widehat{\mathcal{S}\mathcal{V}}^{\vec{d}} \mapsto \underline{p}_1 \in \underline{\Omega}$, with $\underline{\Gamma} \cap \underline{\Omega}$ equal to $4\underline{p}_1+\underline{p}_2+\underline{p}_3$,  $5\underline{p}_1+\underline{p}_2$, $3\underline{p}_1+2\underline{p}_2+\underline{p}_3$ or $3\underline{p}_1+3\underline{p}_2$ (depending on $\vec{d}$). One can check, via the above parameterization, that such morphism is finite, e.g.: for $\vec{d}=(4,1,0,0,0)$, any $\underline{p}_1 \in \underline{\Omega}$ and $\underline{\Gamma}\in \widehat{\mathcal{S}\mathcal{V}}^{\vec{d}}$ in the fiber over $\underline{p}_1$, its projection $g_\mu(\underline{\Gamma})$ is a conic having multiplicity of intersection $I_{q_1}(\mathcal{K}_0,g_\mu(\underline{\Gamma}))\geq 4$ with $\mathcal{K}_0$ at $q_1:=g_\mu(\underline{p}_1)$. Hence, it belongs to the corresponding pencil of conics, a finite number of which is tangent somewhere to $\mathcal{H}_0:=\{z=0\}$.\\
 	\indent	  	 It follows, assuming $\widehat{\mathcal{S}\mathcal{V}}^{\vec{d}}$ not finite and arguing as in the preceding cases, that the natural morphism $\underline{\Gamma} \in \widehat{\mathcal{S}\mathcal{V}}_X^{\vec{d}} \mapsto \underline{p}_1 \in \underline{\Omega}$ is surjective.\\
 		  	 \indent We can also check that over any $\underline{p}_1 \in \underline{\Omega}\cap \underline{r}_0$ the multiplicity jumps by at least one, e.g.: $I_{q_1}(\mathcal{K}_0,g_\mu(\underline{\Gamma}))\geq 5$, when $\vec{d}=(4,1,1,0,0,0)$. On the other hand, $g_\mu(\underline{\Gamma})$ must stay tangent to $\mathcal{H}_0$, meaning $\mathcal{H}_0$ is not a component and they intersect at just one point. There are only two such conics, namely, the zero-loci 
 		  	 $$\{zy\,\textrm{-}\,x^2=0\} \quad\textrm{and}\quad \{zx\,\textrm{-}\,y^2=0\},$$ 
 	  	 both intersecting $\mathcal{K}_0$ at its singular point $[0:0:1]$, and tangent to $\mathcal{H}_0=\{z=0\}$ outside $\mathcal{K}_0$, at $q_X=[0:1:0]$ and $[1:0:0]$ respectively. It follows that such divisor $\underline{\Gamma}$ should decompose as $\underline{\Gamma}'+\iota_\mu(\underline{\Gamma}')+\underline{r}_0$. Contradiction!\\
 	  	 \indent Following the same line of arguments we deal hereafter with the other cases. For $\vec{d} = (3,2,1,0,0,0)$, one can check that for any $\underline{p}_1 \in\underline{\Omega}\setminus \underline{r}_0$, there is a finite number of conics with multiplicity of intersection $3$ with $\mathcal{K}_0$ at $q_1:=g_\mu(\underline{p}_1)\in \mathcal{K}_0\setminus \{[0:0:1]\}$, tangent to $\mathcal{K}_0$ somewhere else and tangent to $\mathcal{H}_0$ as well. When we specialize at $\underline{p}_1 \in\underline{\Omega}\cap\underline{r}_0$ instead, the multiplicity of intersection at $q_1:= [0:0:1]$ jumps, as already mentioned, i.e.: the divisor $g_\mu(\underline{\Gamma})\cap \mathcal{K}_0$ is equal, either to $5q_1 +q_3$ or to $4q_1+2q_2$. The case $5q_1+q_3$ with $q_1 \neq q_3$ can not happen, as already proved, but we can obtain (four possible conics, which intersect $\mathcal{K}_0$ at):
 	  	 \begin{enumerate}
 	  	 	\item either a divisor of the form $6q_1$;
 	  	 	\item or $4q_1+2q_2$, with $q_2 =[2\zeta^2:4\zeta:9]$ and $\zeta^3=1$. 
 	  	 \end{enumerate}  
   	 The associated conics are $\{zy\,\textrm{-}\,x^2=0\}$ and $\{zy\,\textrm{-}\,x^2+4\zeta xy\,\textrm{-}\,4\zeta^2y^2=0\}$, which intersect $\mathcal{H}_0=\{z=0\}$ outside $\mathcal{K}_0$, namely, at $q_X=[0:1:0]$ and $q_X=[2\zeta:1:0]$, respectively.
   	 
   	 The two former cases correspond to decompositions $\vec{d} \in \mathfrak{Dec}(6)$ with no odd coordinate. Hence, their inverse images should belong to $\psi_\mu(|\gamma_X(\mu,2)|)$ but yet be reducible. Acccording to §$3.10.2)$ they must decompose as $\underline{\Gamma}'+\iota_\mu(\underline{\Gamma}')+\underline{r}_0$. Contradiction!$\quad \blacksquare$\\
   	 
   	We now consider the remaining cases, i.e.: $I_\mu(0)=2$ and $I_\mu(0)=3$. Recall that $I_\mu(0)=3$ if and only if $\mu=(2h+1,0,0,0)$ for some $h \in \mathbb{N}^*$, in which case the cubic $\mathcal{K}_0$ decomposes as sum of three distinct lines (see §$3.8.5$). In particular, at any smooth point $q\in \mathcal{K}_0$ and for any conic $\mathcal{C}$, we have $I_q(\mathcal{K}_0, \mathcal{C})\leq 2$. Hence, $\mathcal{S}\mathcal{V}^{\vec{d}}=\emptyset$ for any decomposition $\vec{d}\neq (2,2,1,1,0,0)$. This settles the case $I_\mu(0)=3$.\\
   	
   	\begin{prop} \underline{$(I_\mu(0) =2$ and $\mu_0\neq 1).$}\hspace*{2mm}\\
   	
   	For any $\vec{d} \in \mathfrak{Dec}(6)$ other than $\vec{d}=(2,2,1,1,0,0)$, the set $\widehat{\mathcal{S}\mathcal{V}}^{\vec{d}}$ is finite.\\
   	
   \end{prop}

\textbf{Proof.}

   	When $I_\mu(0)=2$, we can assume $\mu_1=0
   	=\mu_2=0\neq \mu_3$, i.e.: $\mu=(2h+1,0,0,2l\,\textrm{-}\,2)$ for some $h,l\in \mathbb{N}^*$. In this case the exceptional curve $\underline{\Gamma}_\nu$, with $\nu=(2h+1,1,1,2l\,\textrm{-}\,2) \in \mathbb{T}_3$, also denoted $\underline{C}_3$, is a component of $\underline{\Omega}$, hence fixed by the Geiser involution $\iota_\mu$. In particular the cubic $\mathcal{K}_0$ decomposes as the union of a conic $\mathcal{C}_0$ and the line $\mathcal{H}_3:=g_\mu(\underline{C}_3)$. There exists another involution, say $\sigma \in Aut(\underline{S}_\mu)$, such that $\sigma(\underline{C}_0)=\underline{C}_3$, $\sigma(\underline{s}_1)=\underline{r}_1$, $ \sigma(\underline{s}_2)=\underline{r}_2$ and $ \sigma(\underline{s}_3)=\underline{s}_3$. Its existence follows from the equalities $\mu_1=\mu_2=0$ and can be constructed as pull-back of an involution by the blowing-up $\underline{S}_\mu \to \mathbb{P}^1\times \mathbb{P}^1$ defined in \textbf{A.$4$}. We can also prove that $\sigma$ commutes with $\iota_\mu$. Hence, it induces an involution on $\mathbb{P}^2(\mathbb{C})$, also denoted $\sigma$.\\
   	\indent  Modulo action of 
   	$\mathbb{P}GL(3,\mathbb{C})$ we can make the following identifications: $$\mathcal{H}_0:=g_\mu(\underline{C}_0)=\{y\,\textrm{-}\,z=0\}\;\textrm{,}\quad \mathcal{H}_3:=g_\mu(\underline{C}_3)=\{y+z=0\}\;\textrm{,}\quad\mathcal{C}_0:=\{x^2+y^2\,\textrm{-}\,2z^2=0\}$$ and $\sigma :[x:y:z] \in \mathbb{P}^2(\mathbb{C}) \mapsto [x:\,\textrm{-}\,y:z] \in \mathbb{P}^2(\mathbb{C})$. \\
   	\indent We start ruling out the decompositions $\vec{d}=(3,3,0,0,0,0)$ and $\vec{d}=(5,1,0,0,0,0)$: no conic can have multiplicity of intersection $\geq 3$ with the line $\mathcal{H}_3$.\\
   	\indent  We restrict therefore to the cases $\vec{d}=(4,1,1,0,0,0)$ and $\vec{d}=(3,2,1,0,0,0)$. \\
   	\indent For any $\underline{\Gamma} \in \mathcal{S}\mathcal{V}^{\vec{d}}$, the divisor $g_\mu(\underline{\Gamma})\cap  \mathcal{K}_0$ of $\mathcal{K}_0$ is either equal to $4q_1+q_2+q_3$ or to $3q_1+2q_2+q_3$, with $q_1 \in \mathcal{C}_0$, $q_2\in \mathcal{H}_3$ and $q_j\neq q_k$ if $j\neq k$.
   	 On the other hand, for any $q\in \mathcal{C}_0$ outside $\mathcal{H}_0\cap \mathcal{H}_3$ there is only one irreducible conic tangent to $\mathcal{H}_0$, having multiplicity of intersection $\geq 4$ with $\mathcal{C}_0$ at $q$. Analogously, there are at most four conics tangent to both, $\mathcal{H}_0$ and $\mathcal{H}_3$, having multiplicity of intersection $\geq 3$ with $\mathcal{C}_0$ at $q$. It follows that the natural morphism $\underline{\Gamma} \in \mathcal{S}\mathcal{V}^{\vec{d}} \mapsto q_1 \in \mathcal{C}_0$ has finite fibers. Hence, assuming $\mathcal{S}\mathcal{V}^{\vec{d}}$ not a finite set would imply that its closure $\widehat{\mathcal{S}\mathcal{V}}^{\vec{d}}$ projects onto $\mathcal{C}_0$. Choosing $\underline{\Gamma}\in \mathcal{S}\mathcal{V}^{\vec{d}}$ in a particular fiber, as shown hereafter, we finally obtain a contradiction.
   	\begin{enumerate}
   		\item For $\vec{d}=(4,1,1,0,0,0)$, any element $\underline{\Gamma}$ in the fiber over $q_1=[\sqrt{2}:0:1]$ should project onto the irreducible conic $\{2x^2+y^2+2\sqrt{2}xz=0\}$ or onto $2\mathcal{H}_{q_1}$, twice the tangent to $\mathcal{C}_0$ at $q_1$. The latter are $\sigma$-invariant and intersect $\mathcal{H}_0$ and $\mathcal{H}_3=\sigma(\mathcal{H}_0)$ at $q_X$ and $\sigma(q_X)$, where $q_X$ is equal to $[\textrm{-}\,\frac{\sqrt{2}}{2}:1:1]$ and $q_X:=[\sqrt{2}:1:1]$ respectively. Either way, $\underline{\Gamma}\in \psi_\mu(|\gamma_X(\mu,2)|)$ and $g_\mu(\underline{\Gamma})\cap \mathcal{K}_0$ is equal to $4q_1+2\sigma(q_X)$, with no odd coordinate, forcing $\underline{\Gamma}$ to be reducible. Contradiction!: any reducible divisor in $\psi_\mu(|\gamma_X(\mu,2)|)$, with $\mu_i\neq 1$ for any $i=0,\cdots, 3$, should have $\underline{C}_0$ as one of its components.
   		
   		\item If $\vec{d}=(3,2,1,0,0,0)$ instead, we choose $\underline{\Gamma}$ in the fiber over $q_1=[1:\textrm{-}\,1:1]\in \mathcal{C}_0\cap \mathcal{H}_3$. Its projection should be tangent to both $\mathcal{C}_0$ and $\mathcal{H}_3$ at $q_1$, hence reducible, and equal to $g_\mu(\underline{\Gamma})=2\mathcal{H}_{q_1} $. In particular $\underline{\Gamma}\in \psi_\mu(|\gamma_X(\mu,2)|)$, with $q_X:=[3:1:1]$. Contradiction!: should $g_\mu(\underline{\Gamma})$ be reducible, $\mathcal{H}_0$ would be one of its components.$\quad\blacksquare$

   		\end{enumerate}
 		 
 		 \begin{remark}\hspace*{2mm}\\
 		 	
 		 The proof given above for $(I_\mu(0)=2$ and$)$  $\vec{d}=(3,2,1,0,0,0)$ does not assume $\mu_0\neq 1$. We are thus left with the $($sub-$)$case $\mu\in \{(1,0,0,2l), l\in \mathbb{N}^*\}$ and $\vec{d}=(4,1,1,0,0,0)$. \\
 		 \end{remark}
 		
 		\indent In our particular case, up to the action of $\mathbb{P}(GL(3,\mathbb{C}))$, we have the same characterization as above for $\mathcal{K}_0=\mathcal{C}_0\cup \mathcal{H}_3$, $\mathcal{H}_3$ and the involution of $\mathbb{P}^2(\mathbb{C})$ induced by $\sigma$. The new feature is that for any $X\in \mathfrak{X}$, the curves  $\tilde{r}_0,\tilde{\Gamma}_{(3,0,0,2l)} \subset \tilde{S}_X$ satisfy the equalities $\tilde{r}_0.\tilde{\Gamma}_\mu=1$ and $\tilde{\Gamma}_{(3,0,0,2l)}.\tilde{\Gamma}_\mu=0$, and their sum belongs to $|\gamma_X(\mu,2)|$. Hence, they project onto zero self-intersection curves in $\underline{S}_\mu$, denoted $\underline{r}_0$ and $\underline{\Gamma}_{(3,0,0,2l)}$, although they depend upon $X\in\mathfrak{X}$. Recall again the natural identification of $\mathfrak{X}$ with $\mathcal{H}_0 \setminus \mathcal{K}_0$. In particular, for any $X\in \mathfrak{X}$ there exists, say $c_X\in \mathbb{C}\setminus\{1,\textrm{-}\,1\}$, such that the conic $g_\mu(\underline{r}_0)=g_\mu(\underline{\Gamma}_{(3,0,0,2l)})$ is tangent to $\mathcal{H}_0$ at $[c_X:1:1]$. We start characterizing the latter conic in terms of its tangency properties, and give its explicit equation in $\mathbb{P}^2(\mathbb{C})$.\\
 		
 		\begin{prop} $(\mu\in \{(1,0,0,2l), l\in \mathbb{N}^*\})$\hspace*{2mm}\\
 			
 		For any $X \in \mathfrak{X}$ the curves $\underline{r}_0$ and $\underline{\Gamma}_{(3,0,0,2l)}$ are $\sigma$-invariant.
 		
 		\end{prop}
 	\textbf{Proof.}\\
 	We start checking that the three $\sigma$-invariant divisors
 	
 	$$2\,\underline{\Gamma}_{(0,1,0,2l)}+\underline{s}_1+\underline{r}_1,\quad2\,\underline{\Gamma}_{(0,0,1,2l)}+\underline{s}_2+\underline{r}_2,\quad \underline{\Gamma}_{(0,0,0,2l+1)}+\underline{\Gamma}_{(0,0,0,2l\,\textrm{-}1)}+\underline{s}_3\,,$$
 	
as well as $\underline{r}_0$, for any $X\in \mathfrak{X}$, have same linear equivalence class, say $\underline{P}_\mu$. Moreover, $\underline{P}_\mu.\underline{P}_\mu=0$ and $\underline{P}_\mu.\underline{L}_\mu=2$. Hence $|\underline{P}_\mu|$ is a pencil (cf. \cite{H}) and for any $\underline{p}\in \underline{S}_\mu$ there is a unique divisor in $|\underline{P}_\mu|$ going through it. In particular, since $\underline{P}_\mu.\underline{C}_0=1=\underline{P}_\mu.\underline{C}_3$, we deduce an isomorphism $\underline{C}_0 \to \underline{C}_3$, sending $\underline{s}_j\cap\underline{C}_0$ to $\underline{r}_j\cap\underline{C}_3$ (for $j=1,2)$ and fixing $\underline{C}_0\cap\underline{C}_3$, just as $\sigma$ does on $\underline{C}_0$. Hence they coincide, implying each $\underline{r}_0$ is $\sigma$-invariant.\\
 	\indent Applying $\iota_\mu$ we obtain that $\underline{\Gamma}_{(3,0,0,2l)}$ is $\sigma$-invariant as well, for any $X\in \mathfrak{X}$.
 	 $\blacksquare$\\
 	\begin{remark}\hspace*{2mm}
 		
 		Should $\mathcal{S}\mathcal{V}^{\vec{d} }$ have positive dimension, the natural morphism $\widehat{\mathcal{S}\mathcal{V}}^{\vec{d} } \to \mathcal{H}_0,\; \underline{\Gamma} \mapsto \underline{\Gamma} \cap \mathcal{H}_0$ would be surjective. In that case, given any $q_X=[c:1:1]\;$ $(c\in \mathbb{C}\setminus\{1,\textrm{-}1\})$ there would exist $\underline{\Gamma} \in \widehat{\mathcal{S}\mathcal{V}}^{\vec{d} }$ such that $g_\mu(\underline{\Gamma})$ is tangent to $\mathcal{H}_0$ at $q_X$, meaning $g_\mu(\underline{\Gamma})\cap \mathcal{H}_0=\{q_X\}$, and also intersecting $\mathcal{C}_0$ at just one point. The next lemma dresses the list of all such conics, irreducible or not.\\
 	
 	\end{remark}
 
 \begin{lem}\hspace*{2mm}\\
 	
 	For any $q\in \mathcal{C}_0 \setminus \mathcal{H}_0$, say $q=[\alpha:\beta:\gamma]\in \mathcal{C}_0$ with $\beta\,\textrm{-}\,\gamma\neq0$, let $\mathcal{H}_q$ denote the tangent to $\mathcal{C}_0$ at $q$, i.e.:  the zero-locus $\{\alpha x+\beta y\,\textrm{-}\,2\gamma z=0\}$, and
 	$$\mathcal{C}_q^{(4,1,1)}= \Bigl\{(\beta\,\textrm{-}\,2\gamma)^2x^2+(3\beta^2\,\textrm{-}\,4\beta\gamma+2\gamma^2)y^2+2\alpha\beta xy\,\textrm{-}\,4\alpha\gamma xz\,\textrm{-}\,4\beta\gamma yz+4\beta(2\gamma\,\textrm{-}\,\beta)z^2=0\Bigr\}$$ 
 	
 	The conics $2\mathcal{H}_q$ and $\mathcal{C}_q^{(4,1,1)}$ are the unique ones intersecting $\mathcal{C}_0$ only at $q$, while tangent somewhere to $\mathcal{H}_0$.  More precisely, at $[2\gamma\,\textrm{-}\,\beta:\alpha:\alpha]$ and $[\alpha:2\gamma\,\textrm{-}\,\beta:2\gamma\,\textrm{-}\,\beta]$, respectively.\\
 	
\end{lem}
 \textbf{Proof.}\\
 
 Given $q=[\alpha:\beta:\gamma]\in \mathcal{C}_0$, the conic we are looking for belongs to the pencil generated by $\mathcal{C}_0$ and $2\mathcal{H}_q$. Hence, it is the zero-locus of $u(x^2+y^2\,\textrm{-}\,2y^2)+v(\alpha x+\beta y\,\textrm{-}\,2\gamma z)^2$, for some $[u:v] \in \mathbb{P}^1$. Moreover, it must intersect $\mathcal{H}_0=\{y=z\}$ at a unique point, forcing the discriminant of the equation $(u+\alpha^2v)x^2+2\alpha  (\beta\,\textrm{-}\,2\gamma)vxz+(v(\beta\,\textrm{-}\,2\gamma)^2 \,\textrm{-}\,u)z^2=0$ to vanish, i.e.: $\Delta=4u(u\,\textrm{-}\,2(\beta\,\textrm{-}\,\gamma)^2 v)$.\\
 \indent The solutions $u=0$ and $u=2(\beta\,\textrm{-}\,\gamma)^2 v$ correspond to $2\mathcal{H}_q$ and $\mathcal{C}_q^{(4,1,1)}$. We easily check at last, that they cut $\mathcal{H}_0$ at $[2\gamma\,\textrm{-}\,\beta:\alpha:\alpha]$ and $[\alpha:2\gamma\,\textrm{-}\,\beta:2\gamma\,\textrm{-}\,\beta]$, respectively. $\quad \blacksquare$\\

  \begin{defi} \hspace*{2mm}
  
  Recall that for any $X\in \mathfrak{X}$, the conic $g_\mu(\underline{\Gamma}_{(3,0,0,2l)})$ is tangent to $\mathcal{H}_0$ at a unique point, say $[c_X:1:1]$, where $c_X\in \mathbb{C}\setminus\{1,\textrm{-}\,1\}$. The latter conic will be denoted $\mathcal{C}_{c_X}^{(2,2,2)}$.
  	\end{defi}
  
 \begin{prop} \hspace*{2mm}
 
 	For any $c\in \mathbb{C}\setminus\{1,\textrm{-}\,1\}$ the curve $\mathcal{C}_{c}^{(2,2,2)}$ is the unique irreducible $\sigma$-invariant conic, tangent to $\mathcal{H}_0$ at $[c:1:1]$ and bitangent to $\mathcal{C}_0$. More precisely, $$\mathcal{C}_{c}^{(2,2,2)}:= \Bigl\{x^2+(1\,\textrm{-}\,c^2)y^2\,\textrm{-}\,2cxz+(2c^2\,\textrm{-}\,1)z^2=0\Bigl\}$$ and is tangent to $\mathcal{C}_0$ at $[1:\pm\sqrt{2c^2\,\textrm{-}\,1}:1]$, as well as to $\mathcal{H}_3$ at $[c:\textrm{-}\,1:1]$.
 	\end{prop}
 
 \textbf{Proof.}\\
 
 Any $\sigma$-invariant conic $\mathcal{C}$, tangent to $\mathcal{H}_0$ at $[c:1:1]$ and bitangent to $\mathcal{C}_0$, say at a pair of points $\{q, \sigma{(q)}\}\in \mathcal{C}_0^{(2)}$, belongs to the pencil generated by $\mathcal{C}_0$ and $\mathcal{H}_q+\mathcal{H}_{\sigma{(q)}}$. \\
  \indent In other words, if $q=[\alpha:\beta:\gamma]\in \mathcal{C}_0$, then $\sigma(q)= [\alpha:\textrm{-}\,\beta:\gamma]$ and $\mathcal{C}$ must be the zero-locus of $u(x^2+y^2\,\textrm{-}\,2z^2)+v(\alpha x +\beta y \,\textrm{-}\,2\gamma z)(\alpha x \,\textrm{-}\,\beta y\,\textrm{-}\,2\gamma z)$ for some $[u:v] \in \mathbb{P}^1(\mathbb{C})$.\\
  \indent Intersecting $\mathcal{C}$ with $\mathcal{H}_0$ amounts to setting $y=z$ and forcing the discriminant $\Delta$ of the corresponding equation to vanish, i.e.
   $$(u+\alpha^2v)x^2\,\textrm{-}\,4\alpha  \gamma vxz+(v(4\gamma^2 \,\textrm{-}\,\beta^2) \,\textrm{-}\,u)z^2=0\quad \textrm{and}\quad\,\Delta=\textrm{-}\,4(u\,\textrm{-}\,\alpha^2 v)(u\,\textrm{-}\,\beta^2 v)=0.$$ 
    The solution $u=\beta^2 v$ gives the (obvious) reducible conic $2\mathcal{H}_{q,\sigma(q)}$, twice the line joining $q$ and $\sigma(q)$. Choosing $u=\alpha^2 v$ instead, gives the(generically irreducible) conic $$\Bigl\{\alpha^2(x^2+y^2\,\textrm{-}\,2z^2)+(\alpha x\,\textrm{-}\,2\gamma z)^2\,\textrm{-}\,\beta^2y^2=0\Bigr\}\quad\textrm{i.e.:}\quad\Bigl\{\alpha^2 x^2 +(\alpha^2\,\textrm{-}\,\gamma^2)y^2\,\textrm{-}\,2\alpha \gamma xz+\beta^2z^2=0\Bigr\}$$
    as announced. 
    We easily check at last that it intersects $\mathcal{H}_0$ at any $[\alpha:\beta:\gamma]\in \mathcal{C}_0$ satisfying $[\gamma:\alpha]=[c:1]$.$\quad \blacksquare$\\

 		\begin{prop}$($\underline{$\,(I_\mu(0)=2$ and $\mu_0=1)$}$)$\hspace*{2mm}
 		 			
 		For any $\vec{d} \in \mathfrak{Dec}(6)$ other than $\vec{d}=(2,2,1,1,0,0)$, the set $\widehat{\mathcal{S}\mathcal{V}}^{\vec{d}}$ is finite.\\
 			
 		\end{prop}
 		
 		\textbf{Proof.}\\
 		
 	As explained in \textbf{B.}$11$, we only need to treat the case $\vec{d}=(4,1,1,0,0,0)$.\\
 \indent Assume therefore that $\widehat{\mathcal{S}\mathcal{V}}^{(4,1,1,0,0,0)}$ is not finite, or equivalently, that the natural morphism
 	$ \widehat{\mathcal{S}\mathcal{V}}^{(4,1,1,0,0,0)} \to \mathcal{H}_0,\,\underline{\Gamma} \mapsto g_\mu(\underline{\Gamma})\cap \mathcal{H}_0$ is surjective. Consider $\underline{\Gamma}$ in the fiber over $q_X:=[0:1:1]$. Its projection is a conic tangent to $\mathcal{H}_0$, just as $\mathcal{C}_0^{(2,2,2)}$. We will find hereafter the equation of $g_\mu(\underline{\Gamma})$ and then check that $I_{q_X}(g_\mu(\underline{\Gamma}),\mathcal{C}_0^{(2,2,2)})=2$, obtaning a contradiction. Both conics being in $g_\mu(\psi(|\gamma_X(\mu,2))|)$ they should intersect at $q_X$ with multiplicity $\geq 3$, according to \textbf{B.}$1.4)$. \\
 	\indent Reducible or not, the conic $g_\mu(\underline{\Gamma})$ intersects $\mathcal{C}_0$ at only one point, say $q\in \mathcal{C}_0$, and $I_{q_X}(g_\mu(\underline{\Gamma}),\mathcal{H}_0)=2$. According to \textbf{B.}$7$, the conic $g_\mu(\underline{\Gamma})$ is equal, either to $2\mathcal{H}_q$ or to $\mathcal{C}_q^{(4,1,1)}$, with $q=[\alpha:\beta:\gamma]\in \mathcal{C}_0$ satisfying $$[2\gamma\,\textrm{-}\,\beta:\alpha:\alpha]=[0:1:1]\quad\textrm{or}\quad[\alpha:2\gamma\,\textrm{-}\,\beta:2\gamma\,\textrm{-}\,\beta]=[0:1:1],$$ respectively.\\
 	\indent The reducible case is not an option: according to \textbf{B.}$1.1)$ the line $\mathcal{H}_0$ should be one of its components. It follows that $g_\mu(\underline{\Gamma})=\mathcal{C}_q^{(4,1,1)}$, with $q= [0:\sqrt{2}:1]$ or $q=[0:\textrm{-}\,\sqrt{2}:1]$. \\
 	\indent The latter conics are the zero-loci of

$$\Bigl\{(3\pm2\sqrt{2})x^2+(4\mp 2\sqrt{2})(y\,\textrm{-}\,z)^2+(8\mp6\sqrt{2})(y\,\textrm{-}\,z)=0\Bigr\}.$$

 	 We immediately check that $\mathcal{C}_0^{(2,2,2)}=\Bigl\{x^2+(y\,\textrm{-}\,z)^2+2(y\,\textrm{-}\,z)=0\Bigr\}$ intersects them at $q_X=[0:1:1]$ with multiplicity $2$.           $\quad  \blacksquare$
 		
 		\begin{cor}\hspace*{2mm}
 						
 			For any $\mu\in \mathbb{T}_0$ and $\vec{d}\neq (2,2,1,1,0,0)$, the set $\widehat{\mathcal{S}\mathcal{V}}^{\vec{d}}$ is finite, i.e.: for generic $X\in \mathfrak{X}$ the set $\mathcal{S}\mathcal{V}_X^{\vec{d}}$ is empty.\\
 			\end{cor}
 	
	 		\newpage
	 \underline{\textbf{Funding and/or Conflict of interests/Competing interests}}\\
	 
	 \textbf{Declaration of interests}: the author declares he has no competing financial interests or personal relationships that could have appeared to influence the work reported in this article.\\


\begin{thebibliography}{99}
	 		
	 		\bibitem{AKM}{Airault H., Mc Kean H.P. \& Moser J.} 
	 		{\textit{Rational and elliptic solutions of the Korteweg-de Vries
	 			equation and a related many body problem.}} 
	 		{Comm. Pure Appl. Math., vol 30 (1977), 95-148}{ (https://doi.org/10.1002/cpa.3160300106).}
	 		
	 		\bibitem{B}{A. Beauville, }{\textit{Surfaces Alg\'ebriques Complexes}, }{Astérisque, \textbf{54}, 1978.}
	 		
		\bibitem{BeBoEM}{Belokolos E.D., Bobenko A.I., Enolskii V.Z. \& Matveev V.B.}
	{ \textit{Algebraic-geometric principles of superposition of finite-zone solutions of integrable non-linear equations.}}{(Russian) Uspekhi Matem. Nauk \textbf{41}, $n^\circ 2$, 2-42,}{ English translation: Math. Surv. \textbf{41} (1986) 1-49}{  (https://doi.org/10.1070/RM1986v041n02ABEH003241 .}
	
	\bibitem{Do}{Dolgachev I.}{\textit{Topics in Classical Algebraic Geometry}}{ mathweb.ucsd.edu/~eizadi/207A-14/Dolgachev-topics.pdf (2010)} 
	
	
	\bibitem{D1}{ Dubrovin B.A.}
	{\textit{ Periodic problems for the Korteweg-de Vries equation in the class of finite-gap potentials.} }{Funk. Anal. Priloz. \textbf{9}:3 (1975) 41-51 }{(https://doi.org/10.1007/BF01075598).}
	
	\bibitem{D2}{ Dubrovin B.A.}
	{\textit{Theta functions and non-linear equations.}}
	{Russ. Math. Surveys,Uspekhi Mat. Nauk vol 2 (1981).}{English translation: Russian Math. Surv. \textbf{36} (1981), 11-92.}{ (https://doi.org/10.1070/RM1981v036n02ABEH002596).}
	
	\bibitem{DN}{ Dubrovin B.A  \& Novikov S.P.}{\textit{Periodic and conditionally periodic analogs of the multisoliton solutions of the KDV equation.}}{ (Russian) Z. Esper Teoret. Fiz. \textbf{67} (1974) n6 2131-2144,}{ English translation: Soviet Physics JETF 40 (1974) $n^\circ 6$ 1058-1063.}      
	
	\bibitem{DMN}{ Dubrovin B.A, Matveev V.B. \& Novikov S.P.}
	{\textit{Non linear equations of $KdV$ type, finite zone linear operators, and 
			abelian varieties.}}
	{English translation : Russian math. Surv., vol. 31 (1976), 59-146}{ (https://doi.org/10.1070/RM1976v031n01ABEH001446).}
	
	\bibitem{F}{ Fay J.D.}
	{\textit{On the even-order vanishing of Jacobian theta functions.}}
		{Duke Math. J., vol 51 (1984), 109-132}{ (https://doi.org/10.1215/S0012-7094-84-05106-8).}
		
		\bibitem{GW}{ Gesztesy F. \& Weikard R.}
		{\textit{Treibich-Verdier potentials and the stationary $(m)KdV$ hierarchy.}} 
		{Math. Zeitsch. vol 219 (1995), 451-476}{ (https://doi.org/10.1007/BF02572375).}
		
		\bibitem{H}{Harbourne B.}{\textit{ Anticanonical rational surfaces,}}{Trans.AMS, \textbf{349}-3 (1997), 1191-1208}{ (https://doi.org/10.1090/S0002-9947-97-01722-4).}
		
		\bibitem{I}{Ince E.L.}
		{\textit{Further investigations into the periodic Lam\'e equation.}}{Proc. Roy. Soc. Edin., \textbf{60} (1940), 88-99.}
			
				
		
				
				\bibitem{IM} {Its A.R.\& Matveev V.B.}
				{\textit{Schr\"{o}dinger operator with finite-band spectrum and $N$-soliton solutions of the Korteweg-de Vries equation.}}{Theoret. Mat. Fiz. 23 $n^\circ 1$  51-68 (Russian)}{English translation: Theor. Math. Phys. (1975) 345-355}{ (https://doi.org/10.1007/BF01038218).}
					
					\bibitem{K1}{ Krichever I.M. }
					{\textit{Algebraic-geometric construction of the Zakharov-Shabat equations
							and their periodic solutions.}}
					{Dokl. Akad. Nauk  227 $n^\circ 2$ (1976),
					English transl. : Soviet Math. Dokl., vol 17 (No 2)
					(1976), 394-397.}
		
					
					\bibitem{K2}{ Krichever I.M.}
					{\textit{Methods of algebraic geometry in the theory of non-linear equations.}} {English transl. : Russ. Math. Surv. vol 32 (1977)}{ (https://doi.org/10.1070/RM1977v032n06ABEH003862).}
					
					\bibitem{K3}{ Krichever I.M.}
					{\textit{Elliptic solutions of the $KP$ equation and integrable systems of	particles.}}
					{English transl. : Funct. An. vol 14 ($n^\circ 4$ ) (1980)}{ (https://doi.org/10.1007/BF01078304).}
					
					
					\bibitem{N}{Novikov S.P.}
					{\textit{The periodic problem for the KDV equation.}} {(Russian)}{Funct. Anal. Appl. \textbf{8} (1974) 236-246}{ (https://doi.org/10.1007/BF01075697).}
					
					\bibitem{SW}{ Segal G. \& Wilson G.}
					{\textit{Loop groups and equations of $KdV$ type.}}Publ. Math. de l'I.H.E.S., vol 61 (1985), 5-65,
					Publ. Math. de l'I.H.E.S., vol 61 (1985), 5-65{ (https://doi.org/10.1007/BF02698802).}
					
					\bibitem{S1}{Smirnov A.O.}
					{\textit{Elliptic solutions of the Korteweg-de Vries equation.}}
					{Mat. Zametki \textbf{45} : 6 (1989), 66-73 }{ English transl.: Math. Notes 45 
					(1989), $n^\circ 5$-$6$, 476-481}{ (https://doi.org/10.1007/BF01158237).}
					
					\bibitem{S2}{ Smirnov A.O.}
					{\textit{Finite gap elliptic solutions of the $KdV$ equation.}}
					{Acta Appl. Math., vol 36 (1994), 125-166}{ (https://doi.org/10.1007/BF01001546).}
					
					
					\bibitem{TV1}{ Treibich A. \& Verdier J.-L.
						with an appendix by Oesterl\'e J. : \textit{Solitons elliptiques.}}	{Prog. in Math. vol 88, The Grothendieck Festschrift, Vol III, Ed. :
						Birkhauser-Boston, 437-480 (1990)}{ (https://doi.org/10.1007/978081764576211).}
						
						\bibitem{TV2}{ Treibich A. \& Verdier J.-L.}
						{\textit{Rev\^etements  exceptionnels et sommes de 4 nombres triangulaires.}}
						{Duke Math. J., vol 68, $n^\circ 2$  (1992)}{ (https://doi.org/10.1215/S0012-7094-92-06809-8).}
						
						\bibitem{TV3}{ Treibich A. \& Verdier J.-L.}
						{\textit{Au del\`a des potentiels et rev\^etements tangentiels hyperelliptiques
								exceptionnels.}}
						{C.R. Acad. Sc. Paris, t. 325, S\'erie I, 1101-1106, (1997)}{ (https://doi.org/10.1016/S0764-4442(97)88713-7).}
						
						\bibitem{TV4}{ Treibich A. \& Verdier J.-L.}{\textit{Variétés de Kritchever des solitons elliptiques de KP.}}{ Proceedings of the Indo-French Conference on Geometry (Bombay, 1989), 187–232.
								Hindustan Book Agency, Delhi, 1993}
							{ISBN:81-85931-03-8}
							
							
							\bibitem{T1}{Treibich A.}{\textit{Tangential Polynomials and Elliptic Solitons. }}{Duke Math. J. 59(3): 611-627 (December 1989)}{ (https://doi.org/10.1215/S0012-7094-89-05928-0).}
						
						
						
						
						\bibitem{T2}{ Treibich A.}{\textit{On hyperelliptic tangential covers and elliptic finite-gap potentials.}} {(2001) Russ. Math. Surv. \textbf{56} 1107.}{ (https://doi.org/10.1070/rm2001v056n06abeh000454).}
						
						
						\bibitem{T3}{ Treibich A.}{\textit{Compactified Jacobians of Tangential Covers.}}{ Integrable Systems. The Verdier Memorial Conference.}{ Prog. in Math. Vol.$115$}.
						
						
							 
						
						
			\end{thebibliography}
\end{document}